\documentclass[11pt]{amsart}
\usepackage{latexsym,graphicx}
\numberwithin{equation}{section}
\theoremstyle{plain}

\theoremstyle{remark}

\theoremstyle{definition}

\newcommand{\D}{{\mathcal D}}
\newcommand{\E}{\mathcal E}

\newcommand{\G}{{\mathcal G}}

\newcommand{\K}{{\mathcal K}}
\renewcommand{\L}{{\mathcal L}}
\newcommand{\M}{{\mathcal M}}
\newcommand{\N}{\mathbb N}

\newcommand{\R}{\mathbb R}

\newcommand{\V}{{\mathcal V}}

\newcommand{\dist}{\operatorname{dist}}

\newcommand{\fp}{\operatorname{FP}}

\newcommand{\Int}{\operatorname{Int}}

\renewcommand{\span}{\operatorname{span}}

\newcommand{\supp}{\operatorname{Supp}}

\def\Ga{\Gamma}
\def\half{{1 \over 2}}

\newcommand{\oa}{\overrightarrow}

\newcommand{\ol}{\overline}

\def\XXint#1#2#3{{\setbox0=\hbox{$#1{#2#3}{\int}$}
      \vcenter{\hbox{$#2#3$}}\kern-.5\wd0}}

\newcommand{\note}[1]{\marginpar{\tiny\emph{#1}}}

\begin{document}

\def\cal{\mathcal}

\font\tpt=cmr10 at 12 pt
\font\fpt=cmr10 at 14 pt

\font \fr = eufm10

%\font\AAA=Times.dfont  at 12pt
 %\font\BBB=Times.dfont at 8pt

%\font\AAA=cmr10 at 12pt
%\font\BBB=cmr10 at 8pt

%\def\AAA{\bf}
%\def\BBB{\bf}

\overfullrule=0in

\def\boxit#1{\hbox{\vrule
 \vtop{%
  \vbox{\hrule\kern 2pt %
     \hbox{\kern 2pt #1\kern 2pt}}%
   \kern 2pt \hrule }%
  \vrule}}

  \def\harr#1#2{\ \smash{\mathop{\hbox to .3in{\rightarrowfill}}\limits^{\scriptstyle#1}_{\scriptstyle#2}}\ }

\def\AAA{1}
\def\BB{2}
\def\CC{3}
\def\DD{4}
\def\EE{5}
\def\FF{6}
\def\GGG{7}
\def\HH{8}
\def\II{9}
\def\JJ{10}
\def\KK{11}
\def\LL{12}
\def\MM{13}

\def\ALL{1}
\def\BTA{2}
\def\BL{3}
\def\BRE{4}
\def\CNS{5}
\def\CIL{6}
\def\CRA{7}
\def\DDD{8}
\def\DDR{9}
\def\GEO{10}
\def\HYP{11}
\def\BEL{12}
\def\AC{13}
\def\SURVEY{14}
\def\NOTES{15}
\def\AET{16}
\def\LAG{17}
\def\KRY{18}
\def\PLI{19}
\def\RT{20}
\def\SLO{21}
\def\TRU{22}
\def\TWC{23}
\def\WAL{24}

 \def\GG{{{\bf G} \!\!\!\! {\rm l}}\ }

\def\GL{{\rm GL}}

\def\bll{I \!\! L}

\def\IFF{\qquad\iff\qquad}
\def\bra#1#2{\langle #1, #2\rangle}
\def\bbf{{\bf F}}
\def\bbj{{\bf J}}
\def\Jtn{{\bbj}^2_n}  \def\JtN{{\bbj}^2_N}  \def\JoN{{\bbj}^1_N}
\def\jt{j^2}
\def\jtx{\jt_x}
\def\Jt{J^2}
\def\Jtx{\Jt_x}
\def\bpp{{\bf P}^+}
\def\bpt{{\wt{\bf P}}}
\def\fsh{$F$-subharmonic }
\def\mo{monotonicity }
\def\jet{(r,p,A)}
\def\ss{\subset}
\def\sse{\subseteq}
\def\half{\hbox{${1\over 2}$}}
\def\smfrac#1#2{\hbox{${#1\over #2}$}}
\def\oa#1{\overrightarrow #1}
\def\dim{{\rm dim}}
\def\dist{{\rm dist}}
\def\codim{{\rm codim}}
\def\deg{{\rm deg}}
\def\rank{{\rm rank}}
\def\log{{\rm log}}
\def\Hess{{\rm Hess}}
\def\Hessyp{{\rm Hess}_{\rm SYP}}
\def\trace{{\rm trace}}
\def\tr{{\rm tr}}
\def\max{{\rm max}}
\def\min{{\rm min}}
\def\span{{\rm span\,}}
\def\Hom{{\rm Hom\,}}
\def\det{{\rm det}}
\def\End{{\rm End}}
\def\Sym{{\rm Sym}^2}
\def\diag{{\rm diag}}
\def\pt{{\rm pt}}
\def\Spec{{\rm Spec}}
\def\pr{{\rm pr}}
\def\Id{{\rm Id}}
\def\Grass{{\rm Grass}}
\def\Herm#1{{\rm Herm}_{#1}(V)}
\def\arr{\longrightarrow}
\def\supp{{\rm supp}}
\def\Link{{\rm Link}}
\def\Wind{{\rm Wind}}
\def\Div{{\rm Div}}
\def\vol{{\rm vol}}
\def\foral{\qquad {\rm for\ all\ \ }}
\def\fpsh{{\cal PSH}(X,\f)}
\def\Core{{\rm Core}}
\def\dis{f_M}
\def\Re{{\rm Re}}
\def\rn{\bbr^n}
\def\pp{\cp^+}
\def\plp{\cp_+}
\def\Int{{\rm Int}}
\def\cix{C^{\infty}(X)}
\def\Gr#1{G(#1,\rn)}
\def\Symn{{\Sym(\rn)}}
\def\SymN{{\Sym(\bbr^N)}}
\def\Gpn{G(p,\rn)}
\def\fd{{\rm free-dim}}
\def\SA{{\rm SA}}
 \def\cd{{\cal C}}
 \def\cdt{{\widetilde \cd}}
 \def\cm{{\cal M}}
 \def\cmt{{\widetilde \cm}}

\def\Theorem#1{\medskip\noindent {\bf THEOREM \bf #1.}}
\def\Prop#1{\medskip\noindent {\bf Proposition #1.}}
\def\Cor#1{\medskip\noindent {\bf Corollary #1.}}
\def\Lemma#1{\medskip\noindent {\bf Lemma #1.}}
\def\Remark#1{\medskip\noindent {\bf Remark #1.}}
\def\Note#1{\medskip\noindent {\bf Note #1.}}
\def\Def#1{\medskip\noindent {\bf Definition #1.}}
\def\Claim#1{\medskip\noindent {\bf Claim #1.}}
\def\Conj#1{\medskip\noindent {\bf Conjecture \bf    #1.}}
\def\Ex#1{\medskip\noindent {\bf Example \bf    #1.}}
\def\Qu#1{\medskip\noindent {\bf Question \bf    #1.}}
\def\Exercise#1{\medskip\noindent {\bf Exercise \bf    #1.}}

\def\HoQu#1{ {\AAA T\BBB HE\ \AAA H\BBB ODGE\ \AAA Q\BBB UESTION \bf    #1.}}

\def\pf{\medskip\noindent {\bf Proof.}\ }
\def\qed{\hfill  $\vrule width5pt height5pt depth0pt$}
\def\equdef{\buildrel {\rm def} \over  =}
\def\qedqed{\hfill  $\vrule width5pt height5pt depth0pt$ $\vrule width5pt height5pt depth0pt$}
\def\mathqed{  \vrule width5pt height5pt depth0pt}

\def\V{W}

\def\df{d^{\phi}}
\def\hk{\_{\rm l}\,}
\def\n{\nabla}
\def\w{\wedge}

\def\cu{{\cal U}}   \def\cc{{\cal C}}   \def\cb{{\cal B}}  \def\cz{{\cal Z}}
\def\cv{{\cal V}}   \def\cp{{\cal P}}   \def\ca{{\cal A}}
\def\cw{{\cal W}}   \def\co{{\cal O}}
\def\ce{{\cal E}}   \def\ck{{\cal K}}
\def\ch{{\cal H}}   \def\cm{{\cal M}}
\def\cs{{\cal S}}   \def\cn{{\cal N}}
\def\cd{{\cal D}}
\def\cl{{\cal L}}
\def\cp{{\cal P}}
\def\cf{{\cal F}}
\def\ccr{{\cal  R}}

\def\gerG{{\fr{\hbox{g}}}}
\def\gerB{{\fr{\hbox{B}}}}
\def\gerR{{\fr{\hbox{R}}}}
\def\p#1{{\bf P}^{#1}}
\def\vf{\varphi}

\def\wt{\widetilde}
\def\wh{\widehat}

\def\and{\qquad {\rm and} \qquad}
\def\arr{\longrightarrow}
\def\ol{\overline}
\def\bbr{{\mathbb R}}\def\bbh{{\mathbb H}}\def\bbo{{\mathbb O}}
\def\bbc{{\mathbb C}}
\def\bbq{{\mathbb Q}}
\def\bbz{{\mathbb Z}}
\def\bbp{{\mathbb P}}
\def\bbd{{\mathbb D}}

\def\a{\alpha}
\def\b{\beta}
\def\d{\delta}
\def\e{\epsilon}
\def\f{\phi}
\def\g{\gamma}
\def\k{\kappa}
\def\l{\lambda}
\def\o{\omega}

\def\s{\sigma}
\def\x{\xi}
\def\z{\zeta}

\def\D{\Delta}
\def\L{\Lambda}
\def\G{\Gamma}
\def\O{\Omega}

\def\bd{\partial}
\def\bdf{\partial_{\f}}
\def\lag{Lagrangian}
\def\psh{plurisubharmonic }
\def\ph{pluriharmonic }
\def\pph{partially pluriharmonic }
\def\omp{$\omega$-plurisubharmonic \ }
\def\ffl{$\f$-flat}
\def\PH#1{\widehat {#1}}
\def\lloc{L^1_{\rm loc}}
\def\dbar{\ol{\partial}}
\def\lp{\Lambda_+(\f)}
\def\lpp{\Lambda^+(\f)}
\def\bo{\partial \Omega}
\def\Ob{\overline{\O}}
\def\fc{$\phi$-convex }
\def\PSH{{ \rm PSH}}
\def\SH{{\rm SH}}
\def\totr{ $\phi$-free }
\def\BM{\lambda}
\def\Der{D}
\def\CH{{\cal H}}
\def\RH{\overline{\ch}^\f }
\def\pconv{$p$-convex}
\def\MA{MA}
\def\lagpsh{Lagrangian plurisubharmonic}
\def\hermsk{{\rm Herm}_{\rm skew}}
\def\PSHl{\PSH_{\rm Lag}}
 \def\ppsh{$\pp$-plurisubharmonic}
\def\fp{$\pp$-plurisubharmonic }
\def\fh{$\pp$-pluriharmonic }
\def\Symn{\Sym(\rn)}
 \def\ci{C^{\infty}}
\def\USC{{\rm USC}}
\def\LSC{{\rm LSC}}
\def\fa{{\rm\ \  for\ all\ }}
\def\ppc{$\pp$-convex}
\def\cpt{\wt{\cp}}
\def\ft{\wt F}
\def\ob{\overline{\O}}
\def\Be{B_\e}
\def\K{{\rm K}}

\def\M{{\bf M}}
\def\N#1{C_{#1}}
\def\ds{Dirichlet set }
\def\dir{Dirichlet }
\def\Fa{{\oa F}}
\def\TR{{\cal T}}
 \def\ISO{{\rm ISO_p}}
 \def\Span{{\rm Span}}

\def\ALL{1}
\def\AV{2}
\def\BTA{3}
\def\BL{4}
\def\BRE{5}
\def\CNS{6}
\def\CP{7}
\def\CPW{8}
\def\CIL{9}
\def\CRA{10}
\def\DTT{11}
\def\DON{12}
\def\CG{13}
\def\DDD{14}
\def\DDR{15}
\def\GEO{16}
\def\HYP{17}
\def\BEL{18}
\def\SURVEY{19}
\def\AC{20}
\def\NOTES{21}
\def\AET{22}
\def\TANG{23}
\def\TANGG{24}
\def\LAG{25}
\def\SLE{26}
\def\JTY{27}
\def\KRY{28}
\def\PLI{29}
\def\RT{30}
\def\SLO{31}
\def\SPR{32}
\def\TRUU{33}
\def\TRU{34}
\def\TWC{35}
\def\TWCC{36}
\def\TWCCC{37}
\def\WAL{38}

\def\AAA{1}
\def\BB{2}
\def\CC{3}
\def\DD{4}
\def\EE{5}
\def\FF{6}
\def\HH{7}
\def\II{8}
\def\JJ{9}
\def\DU{10}
\def\GM{11}
\def\LO{12}
\def\BL{13}
\def\EX{14}
 \def\Ph{15}

\def\Aa{A}

\vskip .4in

\def\E{E}
\def\fpsi{{F_f(\psi)}}
\def\bL{{\bf \Lambda}}
\def\bdf{{\bf f}}
\def\UU{U}
\def\bbm{{\bf M}}
\def\gg{{\mathfrak g}}
\def\gra{\delta}
\def\Gm{G}
\def\LF{f}
\def\ON{^{1\over N}}
\def\On{^{1\over n}}
\def\Ga{G\aa rding\ }
\def\GD{G\aa rding-Dirichlet\ }
\def\Sn{{\mathcal S}(n)}
\def\cG{{\overline \G}}
\def\ggg{\gg}
\def\bbM{{\mathbb M}}
\def\Mbf{{\mathbf M}}
\def\Nul{{\mathbf N}}
\def\iv{^{-1}}
\def\ggb{\overline {\gg}}
\def\dlog{D_A \log \, \gg}

\def \xx{x} 
\def\yy{y}
\def\zz{z}
\def\aaa{a}
\def\bbb{b}
\def\ccc{c}

\def\Aa{A}
\def\Ab{B}
\def\Im{{\rm Im}}

\font\headfont=cmr10 at 14 pt
\font\aufont=cmr10 at 11 pt

%\centerline{ \headfont    Determinant Majorization}
%\medskip
%\centerline{\headfont  and the Work of Guo, Phong and Tong}

\title[ DUALITY FOR G\AA RDING OPERATORS AND A NEW INEQUALITY]
{\headfont DUALITY FOR G\AA RDING OPERATORS AND A NEW INEQUALITY OF GM/AM TYPE}

\date{\today}
\author{ \aufont F. Reese Harvey and H. Blaine Lawson, Jr.}
\thanks
{Second author was partially supported by the Simons Foundation}

\maketitle

\vskip .5in

\centerline{\bf Abstract}

This article concerns \Ga Dirichlet operators, and also the theory of \Ga polynomials.
\Ga Dirichlet  operators are those which arise from homogeneous polynomials on the space $S(n)$ of symmetric $n\times n$ polynomials,
which are hyperbolic with respect to any positive definite matrix. (In particular, any polynomial on $\bbr^n$, which is hyperbolic with respect to elements in $(\bbr_+)^n$, gives rise to a \Ga operator on $S(n)$,  in many ways.)

Here a  new duality theorem is developed with many special cases and examples. The dual $\gg^*$ of a \Ga
operator $\gg$ is used to prove a new basic inequality for  $\gg$.  An important Lemma of Guo-Phong-Tong is
improved  upon as a final application.

\medskip

 \vfill \eject

%\section{Introduction}
%\label{intro}

%{\small\tableofcontents}
\ 
\vskip1in
\centerline{\bf Table of Contents} \bigskip

%{{\parindent= .1in\narrower 

 \hskip .5 in  \AAA. Introduction

 \hskip .5 in  \BB. G\aa rding-Dirichlet Polynomial Operators

 \hskip .5 in  \CC.  Useful Characterizations of Completeness

 \hskip .5 in  \DD. A Uniform Ellipticity Formula

 \hskip .5 in  \EE. Derivatives

 \hskip .5 in  \FF.  The Relative Interior  of the Polar Cone

 \hskip .5 in  \HH.   The Gradient Diffeomorphism $D \log\, \gg$

 \hskip .5 in \II.   The Exhaustion Theorem

 \hskip .5in  \JJ. Uniform Ellipticity and Auxiliary Operators

 \hskip .5 in  \DU. A  G\aa rding Duality Theory

 \hskip .5in  \GM.    The Gradient Maps $D\gg$ and $D\gg\ON$

 \hskip .5 in  \LO.   The Duality Version of \Ga's Basic Inequality for 
 
 \hskip .8in  the Operator $\gg(A)^{1\over N}$.

   \hskip .5in \BL. Duality Under Pull-Back of \Ga Polynomials

 \hskip .5 in  \EX.   Important Special Cases of the Pull Back Theorem 13.2

  \hskip .5 in  15. Elliptic Regularization and Duality

    \hskip .5 in  16. Duality for General Cubic G\aa rding Polynomials in Two 
    
   \hskip .5 in    \qquad Variables

 \hskip .5in 17. Completing the  Guo-Phong-Tong Lemma.

 \hskip 1in Appendix A. Some Standard Ellipticity Remarks.

  \hskip 1in Appendix B. Natural   Linear Transforms which Turn

\hskip  1.3 in $\gg$ and  its Dual $\gg^*$ into Elliptic Differential Operators.

\vfill\eject

\centerline{\bf \headfont  \AAA. Introduction}

A {\bf  \Ga  polynomial} is a degree $N$ homogeneous polynomial  $\gg$ on a finite dimensional vector space $V$ along with a connected component $\G$ of $\{\gg>0\}$ such that for some $e\in \G$ and each $v\in V$ the one-variable 
polynomial $t \mapsto \gg(te+v)$ has exactly $N$  real roots.  (See  \S 2 for details.)  Such an $e$ is referred to as a hyperbolic direction and $\G$ is called the {\bf \Ga cone}.
These polynomials were studied in the classical
paper of Lars \Ga [G\aa r]  1957, where he showed that one direction $e\in\G$ is hyperbolic iff all directions in $\G$ are hyperbolic, and that the cone $\G$ is convex.
Such $\gg, \G$ have become important in many areas of mathematics including optimization and programming.
We are particularly interested in differential operators defined by polynomials $\gg$ on $V\equiv \Sn$ = real symmetric $n\times n$-matrices,
 and $\G\supset \Int\,  \cp =  \{A>0\}$.  These are called {\sl \Ga-Dirichlet polynomials} or {\sl \Ga-Dirichlet operators}.
 They include $\det(D^2 u)$ or generally  $\s_k(D^2 u)$, the  so called Hessian operators, as well as the newly introduced 
 $p$-fold sum operator and the Lagrangian Monge-Amp\`ere operator, along with many others.
 However, the main results in this paper have interest for any \Ga polynomial.
 
 The new results include:
 
 $\bullet$\ \ A generalization of the Arithmetic-Mean-Geometric-Mean Inequality to \Ga polynomials.
 
 \noindent
 This inequality is based on:
 
$\bullet$ \ \ A new  {\bf duality theory} for \Ga polynomials.

$\bullet$ \ \ Useful uniformly elliptic results in the \Ga-Dirichlet case.

 Let us start by assuming that 
 $(\gg, \G)$ is a \Ga polynomial of degree $N$ with \Ga cone $\G$ on  a vector space $V$, and mentioning that
 the open convex  cone $\G$ can have an {\bf edge} $E$, defined as the largest linear subspace in $\overline \G$. 
  It turns out that
 the polynomial $\gg$ is invariant under translations by elements of $E$, and is therefore defined on $V/E$.  It will frequently make
 the exposition clearer to assume that $E=\{0\}$.  The reader can easily do the general case by applying the results to 
 $\G/E \ss V/E$.   Interesting examples have edges.
 
  $\bullet$\ {\bf The Laplacian.} \ Consider $\gg(A)=\tr(A)$ on $\Sn$ where the edge  $E=\{\tr = 0\}$ is a hyperplane.
  
   $\bullet$\ {\bf The Complex Monge-Amp\`ere Operator.} \  Consider $\bbc^n = (\bbr^{2n}, J)$ and $\gg(A) = \det_\bbc(A_\bbc)$ with $A_\bbc = \half(A-JAJ)$ where the edge
   $E= \{A: JA=-JA\}$ in $\cs(2n)$, the space of skew-hermitian matrices.  In this case one could also apply the other elementary symmetric functions $\s_k$
   to $A_\bbc$, or in fact one could  apply any homogeneous symmetric polynomial in the eigenvalues $\l_1, ... , \l_n$, which is \Ga on $\rn$ and  whose \Ga cone contains $\bbr_{>0}^n$.  The edge remains the same except for the Laplacian $\s_1$.
   
   %With $\bbh^n= (\bbr^{4n}, I, J, K)$ and $A_\bbh = {1\over 4} (A-IAI-JAJ-KAK)$,
    %analogous remarks apply.

  $\bullet$\ {\bf The Quaternionic Monge-Amp\`ere Operator.} \  
   Consider $\bbh^n= (\bbr^{4n}, I, J, K)$ with $A_\bbh = {1\over 4} (A-IAI-JAJ-KAK)$.  Then
    analogous remarks apply.
    
    \bigskip
    
    \centerline{\bf Duality Theory}
    
    \medskip

    First,  we define the {\bf open polar} $\G^*$ to be the relative interior of the polar cone to $\G$, where 
    the polar cone is defined to be
    $$
    \G^0\  \equiv\  \{A^* \in V^* : A^*\bigr|\G \geq 0\},
    $$
    and {\sl interior} is taken in the span of $\G^0$.
 Then in \S\S  \   \HH \ and \II \ we prove that:
$$
\text{{\bf The function $f\equiv \log\,\gg$ is strictly concave,} and}
$$
 
 {\bf The gradient of $\log\, \gg$ gives a diffeomorphism}
  $$
 F \ \equiv \ D\log\, \gg : \G \arr\G^*
 $$
 of degree  -1.  Moreover, for $y\in\G^*$, 
 $$
 \begin{aligned}
&\text{\bf  The function $\bra xy -\log\, \gg(x)$ is a strictly convex exhaustion} \\
& \hskip  1in \text{\bf of the open convex cone $\G$}.
\end{aligned}
$$

 \Def{\DU.2} The {\bf dual \Ga function} $\gg^*$ on $\G^*$  is defined by 
$$
\gg^*(y) \ \equiv \ {1\over  \gg(F\iv (y))} \quad\text{for  $y\in \G^*$.}
%\eqno{(\DU.1)}
$$
The {\bf dual potential} $f^*$ of the potential $f= \log\, \gg$ is defined by
$$
f^*(y) \ \equiv \log\, \gg^*(y).
$$

 We have the following three identities.  The first two follow from Definition \DU.2

\Prop{\DU.3}  {\sl  Suppose $x\in\G$ and $y\in\G^*$ correspond, that is, $y=F(x)$.}
Then
$$
\gg(x)\gg^*(y) \ =\ 1
%\eqno{(\DU.3)}
$$ 
$$
f(x) + f^*(y) \ =\ 0
%\eqno{(\DU.4)}
$$
$$
\bra xy \ \equiv  \ N.
%\eqno{(\DU.5)}
$$

When $(\gg, \G)$ equals the basic example $(\det, \cp)$, the gradient map $F(A)\equiv D_A \log\,\det = A^{-1}$,
and $\gg^*=\gg$ is self-dual.
 Thus the first identity  above can be considered as a generalization of the fact that 
$\det (A^{-1} )= (\det A )^{-1}$.

We have the following result in Section 10 (where the items  are interchanged for the proof).

  \Theorem{\DU.1}  {\sl
 The dual  function $\gg^*$ on $\G^*$ has the following properties which mimick the properties of $\gg, \G$.
 
 (1)  \ \  $\gg^*$ extends to a continuous function on the closure $\overline{\G}^*$ where 
 $$
 \gg\bigr|_{\partial \G^*} \ \equiv\ 0 \and \gg\bigr|_{\G^* \  > \ 0.}
 $$
 
 (2)  \ \ $\gg^*$ is also homogeneous of degree $N$.
  
 (3) \ \ $\gg^*$ is real analytic on $\G^*$.
 
 (4)  \ \ $\log\,\gg^*$ is also strictly concave on $\G^*$.
 
 (5)  \ \ For $x\in\G$, the function $\psi(\yy) \equiv \bra \xx\yy - \log\, \gg^* (\yy)$ 
 is a strictly convex exhaustion function on $\G^*$.
 
 (6) \ \  {\bf (Monotonicity).} \  For $y, z \in \G^*$, $\gg^*(y+z)> \gg^*(y)$.  
  }
 
 As a result of this theorem, we shall refer to $\G^*$ as the {\bf G\aa rding-like cone associated to 
 the \Ga dual function $g^*$}.
 However, as discussed below, $g^*$ is not always a polynomial.

\Theorem{\DU.5} {\sl The dual gradient map
$$
F^*(y) \ \equiv \ D_y f^*\ =\ D_y(\log\, \gg^*) \quad \forall \, y\in \G^*
$$
is the inverse of the gradient map}
$$ 
F(x) \ =\ D_x f \ =\ D_x \log\,\gg  \hskip.2in \forall \, x\in\G.
$$

\Theorem {\DU.8} {\sl The {\bf dual potential}  $f^*(\yy) = \log\,\gg^*(y)$ of the potential function
 $f(\xx) \equiv \log \, \gg(\xx)$ is given by}
 $$
f^*(\yy) \ \equiv \ \inf_{\xx\in\G}  (\bra\yy\xx  -f(\xx) - N)  \quad {\rm for}\ \ y\in\G^*.
%\eqno{(\DU.15)}
$$

In \S \II \ we prove the following very useful result, referred to above, which for \Ga-Dirichlet 
operators is stated as follows.  (We continue to assume  that $E=\{0\}$ for clarity.)

\noindent
{\bf Exhaustion Theorem  \II.1}    {\sl
Suppose $\gg, \G$ is a \GD operator.  For each $B\in \G^*$, the open polar, the function 
$$
\psi(A) \ \equiv\ \bra BA - \log\,\gg(A), \qquad A\in \G
\eqno{(\II.1)}
$$
is a  exhaustion function for the G\aa rding cone $\G$, that is, the prelevel sets
$$
K_c \ \equiv\ \{A\in \G\cap S : \psi(A)\leq c\}, \qquad c\in\bbr,
$$
are compact, and exhaust $\G$, i.e., $$\bigcup_c K_c = \G.$$

Moreover $\psi$ is $C^\infty$ and strictly convex on $\G$.
}

Now by Theorem \DU.1 (5) above, this result also holds in the dual case.

\noindent
{\bf Interesting Note.} At this point the reader may be wondering:

{\bf Is  the homogeneous, degree $N$ dual function $\gg^*$    also
a polynomial?}

The answer is {\bf NO.}  

An explicit example $h(x) = (x_1+x_2)x_1x_2$, corresponding to the product of the Laplacian and the Monge-Amp\`ere operators in $\bbr^2$, whose dual $h^*$ is not even a rational function, 
is given in Proposition 16.2
Of course, $\gg^*$ is always an algebraic function, since it is the reciprocal  of $\gg$ composed with the inverse of the polynomial map  $D \log \gg$.  In Section 16 the dual function-cone pair 
$\gg^*, \G^*$ is computed explicitly for all cubics $\gg(x_1, x_2)$ in two variables.
Up to a linear coordinate change the example $h$ above is generic. This indicates that $\gg^*$ is 
not rational generically for higher dimensions and degrees.

However, $\gg^*$ has most of the remaining properties of $\gg$, and they are enough to apply
 the same process we used for  $\gg, \G$ to the pair  $\gg^*, \G^*$,  yielding Theorem \DU.11 which states that
$$
(\gg^*)^* \ = \ \gg.
$$
In this sense one can say that 
$$
\text{$\gg, \G$ and $\gg^*, \G^*$, are "dual"}.
$$

Section 13 examines duality after pulling back a \Ga pair $\gg, \G$ by a linear map $L$.
The pull back function $\gg_L$ and the pull back cone $\G_L$ also form a  \Ga pair (cf. Lemma 13.1), 
and the various quantities associated with $\gg_L, \G_L$, such as:  the polar cone $\G_L^*$,
the various gradient maps, the dual function $\gg_L^*$, and elliptcity, 
are expressed in terms of similar quantities for the given pair $\gg,\G$.
In the more difficult case where $L$ is not onto, a key point in understanding the gradient
diffeomorphism $D\,\log\,\gg_L$ (and its dual function) is that it factors through a submanifold
$\Sigma\ss \G_L^*$ as a composition of two, easier to understand, diffeomorphisms.
See the main Theorem 13.2 for the details.

 % \vskip.3in
  
  The reader may wonder {\bf  whether this duality takes differential operators to differential operators.}
  That is, if $\gg, \G$ is \Ga-Dirichlet operator, which means that $\cp \ss \overline\G$, could it also happen that 
  $\cp \ss \overline \G^*$ so that $\gg^*$ would also be a differential operator?   The answer is "never",  with one exception. 
  If both inclusions are true, then taking polars of $\cp\ss \overline\G$ reverses the inclusion, yielding 
  $\cp\ss  {\overline\G}^* \ss \cp^0 = \cp$ which proves that $\overline\G = \cp$.

   On the other hand {\bf  there is an interesting affirmative answer if one slightly generalizes the  question}.
  There exists a family of symmetric linear transformations $L_{t}: \Sn \to \Sn$ which,  as $t \downarrow 0$, widens the
 cone $\overline \G$ until at a certain first  point $t_0$, one has $\cp\ss L_{t_0} \overline \G$, and so at this point $\gg\circ L_{t_0}\iv$ becomes a \GD polynomial with \Ga cone $L_{t_0} \overline\G$.
  Furthermore, as in Krylov [K],  for $0< t <t_0$,  the \Ga-Dirichlet polynomial 
  operator-subequation pair $\gg \circ L_t\iv,  L_t \overline\G$
becomes uniformly elliptic. 
 
Now $L_t$  also transforms the dual space to $\Sn$, which is naturally identified with $\Sn$ by the inner
product $\bra AB = \tr(AB)$; and what is said above also  holds for the dual operator $\gg^*$.  That is,
in a similar  way $\gg^*, { \overline \G}^*$ can also be turned into %an elliptic and then 
a uniformly elliptic differential operator.
This gives {\bf  a natural duality between these differential operators}, which is associated to the basic duality 
of \S \DU.
  
  This family $L_t$ is  given by 
  $$
  L_t(A)  \ \equiv\ A +( t-1)  (\tr\, A) {1\over n} I, \quad {\rm for} \ t\in \bbr.
  \eqno{(\Ab.1)}
  $$
It has the following properties: 
 $
 \tr \, L_t(A) =  t \, \tr\, A,  L_1 \ =\ {\rm Id}, \ \ {\rm and} \ \  L_s \circ L_t \ =\ L_{st}
  $
for all $s, t \in \bbr$.    Since $L_1 = {\rm Id}$, this implies that if $t\neq 0$, $L_t$ is invertible with inverse
   $
   L_t\iv \ =\ L_{1\over t}.
  $

There is a natural transformation {$\phi:\Sn \overset {\cong} \arr  \ \{\tr = 0 \} \oplus (\bbr \cdot  {1\over  n} I)$}
given by 
   $$
   \phi(A) \ \equiv\ \left (A-(\tr\, A) {1\over n}I, \ ( \tr\, A )  {1\over n}I\right )  \ \equiv\ (\s, \tau).
      \eqno{(\Ab.5)}
  $$
  With respect to these coordinates we see that 
  $$
  \boxit
 {$ L_t(\s, \tau) \ =\ (\s, t \tau)$}
     \eqno{(\Ab.6)}
  $$
  which means that $L_t = {\rm Id}_{\{\tr=0\}} \oplus \L_t$ where $\L_t=$ scalar multiplication by $t$.
     
All of this is discussed in Appendix \Ab\  and in section \DU, Theorem \DU.15.

  \vskip.3in

  One of our principal results is in  \S \LO:
  
\noindent
{\bf A New Basic Inequality for the Operator $\gg(A)^{1\over N}$ via \Ga\ Duality.}
Let us review the Arithmetic-Mean-Geometric-Mean Inequality which can be stated in three ways.
For an $N\times N$ positive definite matrix $A$ one has that
$$
(\det  A)\ON \ = \ \underset {\det(B)=1}{ \inf_{{B>0 }}} {1\over N} \bra BA
\eqno{(\LO.1a)}
$$ 
where $\bra AB \equiv \tr\,AB$ is the natural inner product on ${\cal S}_N$, the space of real symmetric matrices.
More precisely, for $A>0$
$$
\begin{aligned}
&(\det A)\ON \ \leq \ {1\over N} \bra BA \quad \forall\, B>0, \det B = 1,   \\
 &\text{\sl and equality holds $\iff$ $B = (\det \, A)\ON A^{-1}$.}
\end{aligned}
\eqno{(\LO.1b)}
$$
The  equivalent symmetric form states that for all $A,B>0$
$$
\begin{aligned}
& (\det A)\ON (\det B)\ON  \ \leq \ {1\over N} \bra BA   \\
\text{\sl with}\ &\text{\sl  equality $\iff \ A=\l B^{-1}$ for some $\l>0$}.
\end{aligned}
\eqno{(\LO.1c)}
$$

See Section \LO\ for a discussion of the case $A\geq 0$  where $A$ is singular.

One purpose of this paper is to generalize this
 classical inequality, with determinant and the open cone $\{ P>0\}$ replaced by 
any G\aa rding polynomial $\gg$ and its open G\aa rding cone $\G$.
(The cone $\G$ for $\gg = \det$ is  $\{ P>0\}$. See Example \HH.2.)
The key new ingredient is the introduction  of our dual G\aa rding function $\gg^*$ discussed above.

Suppose that $\gg$ is a G\aa rding-Dirichlet operator of degree $N$  with G\aa rding cone $\G$
on a real vector space $V$.  Let  $\gg^*$ and $\G^*$  be the dual objects as above. 
We state our result in three different forms, in tandem with the three versions of the AM/GM inequality (\LO.1).

\Theorem {\LO.1}  \ \ {\sl  For $A\in \G$, one has
$$
\gg(A)\ON \ \ \leq\ \  \underset {\gg^*(B)=1}{ \inf_{{B\in \G^* }}} {1\over N} \bra BA 
\eqno{(\LO.2a)}
$$
and equality holds in (\LO.2a) $\iff$ $B \equiv D_A \gg\ON$ modulo the edge $E$ of $\G$.

Equivalently, 
$$
\gg(A)\ON \ \leq\ {1\over N}  \bra PA \quad \forall\, A\in \G \ \ {\rm and}\ \ \forall\, P\in \G^* \ {\rm with}\ \gg^*(P)=1,
\eqno{(\LO.2b)}
$$
with equality if and only if $P\equiv D_A \gg\ON$ modulo the edge $E$.

This inequality can be put in the following symmetric form:
$$
\gg(A)\ON \gg^*(B)\ON \ \leq \ {1\over N} \bra AB \quad \forall\, A\in \G \ \ {\rm and }\ \ B\in \G^*
\eqno{(\LO.2c)}$$
with equality $\iff$ up to a positive scale, $A$ and $B$ correspond under the gradient diffeomorphism
$A \mapsto D_A \log\, \gg$,  in which 
case $\bra AB = N$.
}

\medskip

The singular case $A\in\partial \G$ is also considered in Section \LO.

\bigskip
\noindent
{\bf  Concerning Lemma 4 in the paper [GPT] of Guo-Phong-Tong.}

For operators that are (in our terminology) invariant \Ga-Dirichlet operators, they proved that majorization of the determinant implies a lower bound for the determinant  of the gradient
of the operator.  Here we extend this important lemma by not requiring invariance and by proving   the reverse implication.

  \Prop{\LO.4} {\sl  Let $\gg, \G$ be a \Ga-Dirichlet operator of degree $N$  with \Ga cone $\G\ss \cs_n$
  and $\g>0$ a constant.  The following
are equivalent:

(1) \ \ (Majorization of the Determinant)
$$
n\g^{1\over n} (\det A)^{1 \over n} \ \leq \ \gg(A)^{1\over N}\qquad\forall\, A>0
$$

(2)\ \ (Lower Bound for the Determinant of the Gradient)
$$
              D_B \gg^{1\over N} >0 
             \ \ {\rm and}\ \ 
              \g \ \leq \ \det \, D_B \gg^{1\over N} \qquad \forall\, B \in \G.
$$
}

  Moreover, in [HL$_9$]  (cf. [HL$_8$]) we proved that under the assumption of
  
  \noindent
  {\bf The Central Ray Hypothesis:}  \qquad $D_I\gg \ =\ kI$ \  \ for some $k>0$,
  
  \noindent
and a certain coefficient condition, %which is valid for all invariant \GD operators $\gg, \G$,)
  \noindent
the conclusion (1) above holds (with the constant $n\g^{1\over n} = \gg(I)^{1\over N}$).
  
  This hypothesis, and the coefficient condition, are valid if the \GD operator $\gg, \G$ is O$(n)$-invariant, or invariant under any group
  $G\ss {\rm O}(n)$ for which the ray through $I$ is the only  $G$-invariant ray in $\Sn$ modulo the edge.
  These include U$(n)$ (or Sp$(n)$) invariant operators on the complex (or quaternionic) hermitian symmetric matrices.
  
  We also note that the Central Ray Hypothesis holds if and only if the \Ga\ Laplacian equals $k$
  times the standard Laplacian.  This follows from the first derivative formula (\EE.3).

 Finally in Section \DD\ we prove a  explicit uniform ellipticity formula  which has applications to \GD operators of compact manifolds.
 The result says that if $\gg, \G$ is a complete \Ga-Dirichlet operator and $u$ and $v$ are $C^2$ solutions to
 of the inhomogeneous \Ga-Dirichelt equation
$$
\gg(D^2 u) \ =\ e^f, \quad D^2 u \in \G,
%\eqno{(\DD.1)}
$$
then  the difference $w\equiv u-v$ satisfies a linear equation 
$$
Lw \ \equiv \ \bra {\mathbf M}{w} \ =\ 0
$$
where ${\mathbf M}$ is a continuous, positive-definite matrix-valued function,
that is to say that
$$
\gg(D^2u) - \gg(D^2 v) \ =\ \bra{\bbM(D^2u, D^2 v)}{D^2u-D^2v}
%\eqno{(\DD.2)}
$$
where the coefficient matrix $\bbM(D^2 u, D^2 v)$ is positive definite for $C^2$ functions $u$ and $v$ which are
strictly $\G$-subharmonic.
 
 The explicit formula for the term $\bbM(D^2u, D^2 v)$  utilizes the polarization
$\overline \gg$ of $\gg$ (see Section \BB).

\noindent
{\bf Definition \DD.3.}  Let ${\mathbf M}(G_1, ... , G_{N-1}) \in \Sn$ be the element defined,
using the canonical inner product $\bra B A \equiv \tr (BA)$, on ${\cal S}(n)$ by
$$
  \bra  {{\mathbf M}(G_1, ... , G_{N-1})} A \ \equiv \ \ggb(G_1, ... , G_{N-1}, A)
%\eqno{(\DD.4)}
$$
Note that ${\mathbf M}$ is an ${\cal S}(n)$-valued multi-linear form on $\Sn$.

We now set 
$$
\text{${\mathbf M}_k(A,B) \equiv {\mathbf M}(A, ... , A,B, ... ,B)$ with $k$ $A$'s}.
$$
for $k=0, ... , N-1$.  Then the basic theorem is the following.

\noindent
{\bf Theorem \DD.4.} {\sl  For all $A,B \in \Sn$
$$
\ggg(A) - \ggg(B) \ =\ \bra    {\bbM}{A-B}  
%\eqno{(\DD.6)}
$$
with 
$$
\bbM\ = \ \bbM(A,B)\ \equiv \ \sum_{k=0}^{N-1} {\mathbf M}_k(A,B).
%\eqno{(\DD.7)}
$$
Moreover, if $\gg, \G$ is complete, then
$$
\bbM (A,B) \ >\ 0 \quad{\rm for\  all}\ \ A, B \in \G.
%\eqno{(\DD.8)}
$$
In fact,
$$
\Mbf(G_1, ... , G_{N-1}) \ >\  0 \quad{\rm for\ all}\ \ G_1,... , G_{N-1}    \in \G.
%\eqno{(\DD.9)}
$$
}

This result has been applied [HL$_{10}$] to show uniqueness of solutions to a generalized Calabi-Yau
equation
$$
\ggg(I+\Hess(u)) \ =\ e^f
$$
where $\Hess(u)$ is the riemannian Hessian, and $I+\Hess(u) \in \G$. This is for compact manifolds and also the Euclidean case (Corollaries 4.13, 4.14 and Theorem 4.15).
 
 It also gives the uniform ellipticity of the operator $\ggg(I+\Hess(u))$
 which suggests that the continuity method might be used to establish  \note{We must  review this after working on the Calabi paper.}
 existence on compact manifolds.  This is discussed with some conjectures
 in [HL$_{10}$].

 %%%%%%%%%%%%%%%%%%%%%%%%%%%%%%%%%%%%%%%%%%%%%%%%%%%%%%%%%%%%%%%%%%%%%%%%%%%%%%%%%%%%%%%%%%%%%%%%%%%%%%%%%%%%%%%%%%%%%%%%%%%%%%%%%%%%%%%%%%%%%%%%%%%%%%%%%%%%%%%%%%%%%%%%%%%%%%%%%%%%%%%%%%%%%%%%%%%%%%%%%%%%%%%%%%%%%%%%%%%%%%%%%%%%%%%%%%%%%%%%%%%%%%%%%%%%%%%%%%%%%%%%%%%%%%%%%%%%%%%%%%%%%%%%%%%%%%%%%%%%%%%%%%%%%%%%%%%%%%%%%%%%%%%%%%%%%%%%%%%%%%%%%%%%%%%%%%%%%%%%%%%%%%%%%%%%%%%%%%%%%%%%%%%%%%%%%%%%%%%%%%%%%%%%%%%%%%%%%%%%%%%%%%%%%%%%%%%%%%%%%

\vskip .5in

\centerline{\bf \headfont  \BB.\  G\aa rding-Dirichlet Polynomial Operators.}

 \medskip
 
This important family of pure second order constant coefficient operators on $\rn$ is defined as follows
Let $\Sn \equiv \Sym(\bbr^n)$ denote the space of  real symmetric $n\times n$-matrices, so that
for a $C^2$-function $u$ in $\rn$, the second derivative or hessian of $u$ at $x$ is an element of $\Sn$:
$$
 D^2_x u\  \equiv \ \left( \frac {\partial^2 u}{\partial x_i \partial x_j} (x) \right) \in \Sn
 $$
 To keep notation simple we shall denote elements of $\Sn$ by $A$ rather than  $D^2_x u$.
 
Consider now a homogeneous polynomial $\ggg : \Sn \to \bbr$, of  degree $m$,  with the property
\note{Sometimes this is $m$ and sometimes it is $N$.}
 that $t\to \ggg(tI+A)$ has all real roots for each $A\in \Sn$. We assume $\ggg(I)>0$ and write
 $$
 \ggg(tI+A) \ = \ \ggg(I) \prod_{k=1}^m (t+\l_k(A))
 $$
where the $\l_k(A)$ are the $I$-{\sl eigenvalues} of $\ggg$ (the negatives of the roots).
These are sometimes called the {\sl G\aa rding eigenvalues}.

The set $\G$ where $\l_k(A) >0$ for all $k$ is called the open {\sl G\aa rding cone}. 

We know from G\aa rding's seminal paper  [G\aa r] (cf. [HL$_3$]]) that:

(1)\ \ $\G$ is a convex cone with 0 as vertex.

(2)\ \   $\G$ is the connected component  of $\{A : \ggg(A)\neq 0\}$ containing $I$.

(3)\ \ $\ggg(A)^{{1\over m}}$ is a concave function on the closed G\aa rding cone $\overline \G$,
and 

(4)\ \ (G\aa rding's Inequality)   \ \ \  
$$
\ggg(A_1)^{{1\over m}} \cdots \ggg(A_m)^{{1\over m}} \  \leq \ 
\overline{\ggg} (A_1, ... , A_m), \qquad \forall \, A_1, ... , A_m\in \G
\eqno{(\BB.1)}
$$
where $\overline\ggg(A_1, A_2, ... , A_m)$ is the polarization\footnote{See [G\aa r, \S 4] for example.}
 of $\ggg$.
 The fact that $\overline\ggg$ is a symmetric  multi-linear form with $\overline{\ggg} (A, ... , A) = \ggg(A)$
 characterizes $\overline \gg$.
 
(5)\ \ The {\bf nullity} $N\equiv \{A\in \Sn : \l_1(A) = \cdots = \l_m(A)=0\}$ is equal to the {\bf edge}
$E \equiv \overline \G \cap (-\overline\G)$, and is also  equal to the {\bf linearity} $L\equiv \{A : \ggg(A+B)=\ggg(B) \ \forall\, B\in \Sn\}$.

(6) \ \  If $G\in\G$, then for each $A\in \Sn$, $\ggg(tG+A)$ has all real roots, and $\ggg(tG+A) = \ggg(G)\prod_k (t+\l^G_k(A))$,
where the $\l^G_k(A)$ are called the $G$-eigenvalues\footnote{For $G\in \G$, $\l_j^G(I)=1/ \l_j^I(G)$. See also (2.5) below.}   
of $A$.  One has that  $\G = \{A : \l^G_k(A) >0\ {\rm for}\ k=1,...,m\}$, independent of $G\in\G$.

(7)\ \ If $G\in\G$, then $\gg(G+tA) = \gg(G)
\prod_{j=1}^N (1+t\l_j^G(A))$ for all $A\in \Sn$

(8) \ \  {\bf (Monotonicity)}  \ $\gg(A+B) > \gg(A)$ for all $A, B \in \G$.

For a proof of a stronger version of (8) which is different from the one in [G\aa r], see Lemma \EE.2.

\noindent
{\bf Definition \BB.1} If in addition  $\overline \G$ contains
the cone  $\cp \equiv \{A : A\geq 0\}$, we say that  $\ggg$  is a {\bf   G\aa rding-Dirichlet 
polynomial}, or equivalently a {\bf   G\aa rding-Dirichlet operator} (where we think of $\ggg(D^2 u)$ for smooth functions $u$.)
Note that
 $$
 \cp \ss\overline\G \ \iff\ \overline\G + \cp\ =\ \overline\G \ \iff\  \Int \cp \ss \G 
\ \iff\ \overline\G + \Int \cp=\G
\eqno{(\BB.2a)}
$$

The purpose of this hypothesis is to ensure $C^2$-coherence and operator ellipticity.
 {\bf $C^2$-coherence} means  that  on an open set $X\ss\rn$,
if $u\in C^2(X)$ satisfies the constraint condition
$D^2_x u \in \overline \G, \ \forall\, x\in X$ classically, then $D^2_x u  \in \overline \G$
in the viscosity sense on $X$.\footnote{This means that for all $x\in X$
 and  each $C^2$ test function $\vf$ for $u$ at $x$, $D^2_x \vf \in \cG$}  This should be thought of as the weakest possible form of ellipticity. 
The Dirichlet condition (\BB.2a) on $\overline \G$ will also be referred to as the {\bf positivity condition}.
However, to consider the inhomogeneous equation $\gg(D^2_x u ) = \psi(x)$ in the viscosity sense, we also need operator monotonicity.  This is provided by (8) above, that is,
 $$
\gg(A+P)\ > \ \gg(A) \quad {\rm for\ all\ } A\in \G, P\ > 0
\eqno{(\BB.2b)}
$$
since $\Int \cp \ss \G$. (See [HL$_6$] or [CHLP] for complete details.)

%%%%%%%%%%%%%%%%%%%%%%%%%%%%%%%%%%%%%%%%%%%%%%%%%%%%%%
%%%%%%%%%%%%%%%%%%%%%%%%%%%%%%%%%%%%%%%%%%%%%%%%%%%%%%
%%%%%%%%%%%%%%%%%%%%%%%%%%%%%%%%%%%%%%%%%%%%%%%%%%%%%%
%%%%%%%%%%%%%%%%%%%%%%%%%%%%%%%%%%%%%%%%%%%%%%%%%%%%%%
%%%%%%%%%%%%%%%%%%%%%%%%%%%%%%%%%%%%%%%%%%%%%%%%%%%%%%

%\vskip.3in

\vfill\eject

\centerline{\bf \headfont  \CC.  Useful Characterizations of Completeness.}
\medskip

Frequently an additional condition on a \GD operator, called {\sl completeness}, is required to eliminate
the extreme degeneracy found in examples such as $\gg(D^2 u) ={ \partial^2 u \over \partial x_1^2}$, 
or more generally, $\gg(D^2 u) ={ \partial^2 u \over \partial x_1^2} + \cdots + { \partial^2 u \over \partial x_k^2}$
in $\rn$ with $k<n$.  Here any continuous function $u$ of the variables $x_{k+1}, ... , x_n$ is a viscosity solution to 
this equation $\gg(D^2 u) =0$. This is usually considered an extreme violation of interior regularity.  

Roughly speaking, the \Ga  pair $\gg, \cG$ is {\bf complete}
if it cannot be pulled back from a similar pair $\gg_0, \cG_0$ on a proper vector subspace $V\ss \rn$ by a 
linear map $\rn\to V$.  Equivalently,  $\gg, \cG$ is {\bf incomplete} if there exists a choice of the $x_1$-axis
in $V\equiv \rn$ such that all continuous functions $u(x_1)$ satisfy $\gg(D^2 u)=0$.

\Def{\CC.1} Let $\gg$ be a \GD operator with \Ga cone $\G$.  
Then  {\bf completeness  for}  $\gg$ means that it is not possible that $\gg(D^2 u)$ can be written as 
 ${\mathfrak h}(D^2 u \bigl|_V)$  where ${\mathfrak h}$ is a \GD polynomial on a proper linear subspace $V$ of $\rn$.

Note that this is not the same concept  as requiring that  $\gg$ depend on all the variables in $\Sn$.  For example the diagonal operator $\D(A) = a_{11} + \cdots + a_{nn}$ is complete by Definition \CC.1   even though it does not depend on the variables $a_{ij}$ for $i\neq j$.

Completeness has some subtle consequences (actually, equivalences) which we now explore. (Some readers may want to skip this at first.)

Recall the  {\bf edge} $E\equiv \overline \G \cap (- \overline \G)$ is the largest linear subspace  contained in $\overline\G$
and the {\bf span} $S$ is defined to be its orthogonal complement $S=E^\perp$.
The following result is true for any closed convex cone $\overline\G\ss\Sn$ containing $\cp$.  The proof 
    can be found in Proposition 3.5 in  [HL$_6$] where it is shown that completeness ((1) $ \cong$ (2) in Prop. 3.5) 
    is equivalent to condition (1) below (condition (3) in Prop. 3.5). The equivalence of (1)  with (2) and (3)  below is straightforward.
The equivalences (4), (5) and (6) are proven in Proposition 3.10 in [HL$_7$].  

\Prop {\CC.2} {\sl  
Let $\overline\G\ss \Sn$ be a convex cone subequation, so $\cp\ss \overline\G$, which by definition means
$\G$ is G\aa rding-Dirichlet.
  Then $\overline\G$ is complete iff one of the following
equivalent conditions holds.

(1) \ \ $\cp\cap E = \{0\}$ (The Edge Condition).

(2)  \ \ $(\cp-\{0\})\cap E  = \emptyset$.

(3)  \ \ $\cp-\{0\} \ss \overline\G - E$.

(4) \ \ The projection $P_e \notin E$ (equivalently, $-P_e \notin \cG$) for all $|e|=1$.

(5) \ \ $S\cap (\Int\, \cp) \neq \emptyset$.

(6) \ \ $\G^* \ss\Int\,\cp$. (See Section 6 for a discussion of the open polar $\G^*$.)
}

Positive definiteness of elements of $\G^*$ will prove to be  a powerful way of stating completeness.

\Remark{\CC.3. (The support of $\cG$)} If $\gg, \G$ is a \GD operator which is  not complete,
it can always be replaced by another \GD operator  which is complete. This involves a subspace $W$ of $\rn$  (not of $\Sn$ like the span) called the {\bf support of $\cG$}.   This concept was introduced in Section 4 of [HL$_7$],
and the reader is referred there.

%%%%%%%%%%%%%%%%%%%%%%%%%%%%%%%%%%%%%%%%%%%%%%%%%%%%%%%%%%
%%%%%%%%%%%%%%%%%%%%%%%%%%%%%%%%%%%%%%%%%%%%%%%%%%%%%%%%%%
%%%%%%%%%%%%%%%%%%%%%%%%%%%%%%%%%%%%%%%%%%%%%%%%%%%%%%%%%%
%%%%%%%%%%%%%%%%%%%%%%%%%%%%%%%%%%%%%%%%%%%%%%%%%%%%%%%%%%
%%%%%%%%%%%%%%%%%%%%%%%%%%%%%%%%%%%%%%%%%%%%%%%%%%%%%%%%%%
%%%%%%%%%%%%%%%%%%%%%%%%%%%%%%%%%%%%%%%%%%%%%%%%%%%%%%%%%%

\vskip .5in

\centerline{\bf \headfont  \DD.  A Uniform Ellipticity Formula.}

\medskip

See Appendix A for  some standard background.

In this section we prove that the difference $w\equiv u-v$ of two solutions of the inhomogeneous \Ga-Dirichlet equation
$$
\gg(D^2 u) \ =\ e^f, \quad D^2 u \in \G
\eqno{(\DD.1)}
$$
satisfies a linear uniformly elliptic homogeneous equation by providing a formula of the form
$$
\gg(D^2u) - \gg(D^2 v) \ =\ \bra{\bbM(D^2u, D^2 v)}{D^2u-D^2v}
\eqno{(\DD.2)}
$$
and showing the coefficient matrix $\bbM(D^2 u, D^2 v)$ is positive definite for $C^2$ functions $u$ and $v$ which are
strictly $\G$-subharmonic, under the necessary assumption that $\gg, \G$ is complete (Theorem \DD.4).

As an immediate consequence we have the following.

\Theorem{ \DD.1}  {\sl
If $\gg, \G$ is a complete \Ga-Dirichlet operator and $u$ and $v$ are $C^2$ solutions to (\DD.1), 
then  the difference $w\equiv u-v$ satisfies a linear equation 
$$
Lw \ \equiv \ \bra {\mathbf M}{w} \ =\ 0
$$
where ${\mathbf M}$ is a continuous, positive-definite matrix-valued function.
}

We shall give explicit formulas for the term $\bbM(D^2u, D^2 v)$ in (\DD.2). The first version utilizes the polarization
$\overline \gg$ of $\gg$ (see Section \BB).

\noindent
{\bf Lemma \DD.2.} {\sl 
$$
\ggg(A) - \ggg(B) \ =\ \sum_{k=0}^{m-1} \overline\ggg (A, ... , A , B, ... , B, A-B) \qquad
\text{with $k$ $A$'s and $(N-k-1)$ $B$'s.}
$$
 }

\noindent
{\bf Proof.} 
$$ 
\eqno{(\DD.3)}
\begin{aligned}    
\ggg( A) &= \ggb(A-B,A, ... ,A) + \ggb(B, A, ... , A) \\
&= \ggb(A-B,A, ... ,A) + \ggb(B, A-B,A ... , A) + \ggb(B,B,A, ... ,A) = \cdots   \\
&=\ \sum_{k=0}^{m-1}\ggb (A-B, A, ... , A , B, ... , B)  +  \gg(B).  \qquad\mathqed
\end{aligned}
$$

We can put the RHS of Lemma \DD.2  in the more tractable form $\bra{\bbM}{A-B}$ by using the next definition.

 Fix $G_1, ... , G_{N-1} \in \Sn$ and consider the linear functional on $\Sn$ given by 
 $$
 A\  \to \  \overline\ggg(G_1, ... , G_{N-1}, A)
 $$ 
where $\overline\ggg$ is the polarization of $\ggg$.

\noindent
{\bf Definition \DD.3.}  Let ${\mathbf M}(G_1, ... , G_{N-1}) \in \Sn$ be the element defined,
using the canonical inner product $\bra B A \equiv \tr (BA)$, on ${\cal S}(n)$ by
$$
  \bra  {{\mathbf M}(G_1, ... , G_{N-1})} A \ \equiv \ \ggb(G_1, ... , G_{N-1}, A)
\eqno{(\DD.4)}
$$
Note that ${\mathbf M}$ is an ${\cal S}(n)$-valued multi-linear form on $\Sn$.

In other words,  if we set ${\mathbf M}_k(A,B) \equiv {\mathbf M}(A, ... , A,B, ... ,B)$ with $k$ $A$'s,
then the formula in Lemma \DD.2 can be written as
$$
\ggg(A) - \ggg(B) \ =\ \sum_{k=0}^{N-1}  \bra  { {\mathbf M}_k(A,B)} {A-B}.
\eqno{(\DD.5)}
$$
Summing the ${\mathbf M}_k$'s, puts (\DD.2) in our final form.

\noindent
{\bf Theorem \DD.4.} {\sl  For all $A,B \in \Sn$
$$
\ggg(A) - \ggg(B) \ =\ \bra    {\bbM}{A-B}  
\eqno{(\DD.6)}
$$
with 
$$
\bbM\ = \ \bbM(A,B)\ \equiv \ \sum_{k=0}^{N-1} {\mathbf M}_k(A,B).
\eqno{(\DD.7)}
$$
Moreover, if $\gg, \G$ is complete, then
$$
\bbM (A,B) \ >\ 0 \quad{\rm for\  all}\ \ A, B \in \G.
\eqno{(\DD.8)}
$$
In fact,
$$
\Mbf(G_1, ... , G_{N-1}) \ >\  0 \quad{\rm for\ all}\ \ G_1,... , G_{N-1}    \in \G.
\eqno{(\DD.9)}
$$
}

\pf
It remains to prove (\DD.9).  There are several ingredients in the proof.
One is a criterion for positive definiteness.

\Lemma{\DD.5} {\sl Given $\Mbf \in \Sn$,
$$
    \bra{\Mbf}{P} >0 \ \ \forall \, P\in \cp-\{0\}  \quad \iff  \quad  \Mbf >0 
\eqno{(\DD.10)}
$$
In particular, by Definition \DD.3
$$
\bbM (G_1, ... , G_{N-1}) >0 \ \iff\  \ggb(G_1, ... , G_{N-1}, P)>0\quad\forall\, P\in \cp-\{0\}.
\eqno{(\DD.11)}
$$
}
 Another criterion is number (3) for completeness  in Proposition \CC.4.
$$
\text {$\cG$ is complete \quad $\iff$\quad $\cG - E \supset \cp - \{0\}$. }
\eqno{(\DD.12)}
$$

The last ingredient is the following.

\Prop{\DD.6}
$$
\ggb (G_1, ... , G_{N_1}, A)>0 \quad\forall\, G_1, ... , G_{N_1} \in \G \ \ {\rm and}\ \ \forall\, A\in \cG-E.
$$

The positive definiteness   (\DD.9) in Theorem \DD.4 is an easy consequence of these results.  By the completeness criterion
(\DD.12), Proposition \DD.6 implies that $\gg(G_1, ... , G_{N-1}, P)>0$ for all $G_1, ... , G_{N_1}\in \G$ and $P\in \cp-\{0\}$.
By (\DD.11) this implies that $\bbM (G_1, ... , G_{N-1}) >0 $ completing the proof of Theorem \DD.4 once we prove Lemma \DD.5 and Proposition \DD.6.

\noindent
{\bf Proof of Lemma \DD.5.}
If $M>0$ is false, choose an eigenvector $e$ with eigenvalue $\l\leq 0$.  Take $P\equiv P_e \in \cp-\{0\}$.
Then $\bra M P = \bra M {P_e} = \l \leq 0$, proving that $\bra{\Mbf}{P} >0$ is false. 
 
The converse is  proven as follows.
Note that    $M>0 \Rightarrow M-\e I >0$ for some $\e>0$ 
$\Rightarrow  \bra {M-\e I} P\geq 0 \iff \bra M {P} \geq \e \, \tr P$ 
for all $P\in\cp$.  Finally, for $P\in\cp$, we have  $\tr P>0$ unless $P=0$. \qed

\noindent
{\bf Proof of Proposition \DD.6.}  To begin, assume all the $G_1, ... , G_{N-1}$ are equal to  $I \in\Int \cp\ss \G$.
Then   we have 
$$
\begin{aligned}
N\, \overline\ggg(I, ... , I, A) \ &=\ {d\over dt}\biggr|_{t=0} \overline \ggg(I+tA,I+tA, ... ,I+tA)
 \ =\ {d\over dt}\biggr|_{t=0}\ggg(I+tA)  \\
 & =\ {d\over dt}\biggr|_{t=0}  \prod_{k=1}^N (1+t\l_k(A))
 \ =\ \l_1(A) + \cdots + \l_N(A)   \ =\ \tr\, A,
 \end{aligned}
 $$
the {\sl G\aa rding trace} of $A$.
 Since $A\in \overline \G$, we have $\l_k(A)\geq 0$ for all $k$. 
 Hence, $\tr\, A>0$ unless $\l_k(A)=0$ for all $k$, that is, unless $A\in \Nul$, the nullity.  
 By (5) in Section 2, we have  $\Nul=E$, the edge,  and therefore $A\in E$.

Now for the general case, assume $G_1, ... , G_{N-1} \in\G$.  
We know that $\overline\ggg(G_1, ... , G_{N-1}, A) \geq 0$ from G\aa rding's Inequality
(\BB.1), so we only need to rule out the case where  $\overline\ggg(G_1, ... , G_{N-1}, A) = 0$.
In the case that this is zero, choose $\e>0$ so that $G_1-\e I \in \G$.  Then again by (\BB.1)
$$
\begin{aligned}
&0 \ \leq \ \overline\ggg(G_1-\e I, G_2, ... , G_{N-1}, A) \\
&= \ \overline\ggg(G_1, ... , G_{N-1}, A) - \e\overline\ggg(I, G_2, ... , G_{N-1}, A) \\
&=\ -\e \overline\ggg(I, G_2, ... , G_{N-1}, A).
\end{aligned}
$$
Since $\overline\ggg(I, G_2, ... , G_{N-1}, A) \geq 0$, this proves that 
 $\overline\ggg(G_1, G_2, ... , G_{N-1}, A) = 0$.
 
 Repeating this argument we have  $\overline\ggg(I, ... , I, A)= 0$.   Now our first argument applies
 to yield $A\in E$, but since in Proposition \DD.6 we assume that $A\notin E$, this  
 proves that $\gg(G_1, ... , G_{N-1}, A)>0$.    \qed
 
 This completes the proof of Theorem \DD.4.\qed

\noindent
{\bf Note.}
Since $$D_B\gg = {d\over dt}\ggb(B+tA, ... , B+tA) \bigr|_{t=0} = N \ggb(B, ... , B, A),$$  by definition \DD.3 this proves that
$$
D_B \gg \ =\ N \bbM (B, ... , B) \quad N-1 \ {\rm times}.
 \eqno{(\DD.13)}
$$
Therefore, the positive definiteness conclusion (\DD.9) in Theorem \DD.4 can be considered a
 significant strengthening of the more standard ellipticity statement that the linearizations of 
 $\gg$ at $B$ have coefficient matrix $D_B\gg >0$.
 %Theorem \Aa.1 in Appendix A.

It is interesting at this point to look at some of the consequences of Theorem \DD.4.
Suppose $u_1, u_2 \in C^2(X)$ for an open set $X\ss\rn$.  Consider the following cases.

\noindent
{\bf The Classical Case.}  \ \  Suppose that $u_1, u_2$ are strictly $\overline \G$-subharmonic on $X$, that is,
$D^2_x u_1, D^2_x u_2 \in \G$,   $\forall\,x\in X$.  Set $\bbM \equiv \bbM(D^2_x u_1, D^2_x u_2)$ as in 
Theorem \DD.6.

\noindent
{\bf The Quasi Case.}\ \ Suppose that $u_1, u_2$ are strictly quasi-$\overline \G$-subharmonic on $X$, that is,
$D^2_x u_1, D^2_x u_2 \in \G-I$,  $\forall\, x\in X$.  Set $\bbM \equiv \bbM(D^2_x u_1 + I, D^2_x u_2+I)$ as in 
Theorem \DD.6.  

In both cases, if $\G$ is complete, then $\bbM$ is positive definite.

It is convenient to refer to the to the linear operator $\bra \bbM A$ as the {\bf $\bbM$-Laplacian} and denote it by $\D_\bbM$.
Then the $\bbM$-Laplacian of the difference $u_1-u_2$ is given by:

\noindent
{\bf The Classical  Case.}\ \ $\D_\bbM(u_1-u_2) \ =\ \ggg(D^2_x u_1) - \ggg(D^2_x u_2)$,

\noindent
{\bf The Quasi Case.} \ \ $\D_\bbM(u_1-u_2) \ =\ \ggg(D^2_x u_1+I) - \ggg(D^2_x u_2+I)$.

\noindent
{\bf Corollary \DD.9.}
{\sl  If $\G$ is complete in the discussion above, then in either case
the difference $u_1-u_2$ satisfies the Strong Maximum Principle.}

\noindent
{\bf Corollary \DD.10.}
{\sl
Suppose $\ggg:\Sn \to \bbr$ is a G\aa rding-Dirichlet operator of degree $N$
with  G\aa rding  cone $\G$.  Given $f\in C(X)$, where $X$ is an open subset of $\rn$, 
consider either 

\noindent
{\rm The Classical Equation:}  \ \  $\ggg(D^2 u)\ =\ e^f$, \ \ with $D^2 u\in \overline \G$,  \qquad or

\noindent
{\rm The Quasi Equation:}  \ \  $\ggg(D^2 u+I)\ =\ e^f$, \ \ with $D^2 u + I \in \overline \G$.

\noindent
Suppose that $u, H\in C^2(X)$ with $u$ a subsolution and $H$  a solution.
Then the difference $u-H$ is $\D_\bbM$-subharmonic on $X$, that is,
$$
\D_\bbM(u-H) \ \geq \ 0.
$$
}

\noindent
{\bf Proof.} \rm
For Case 1,   $\D_\bbM(u-H) = \ggg(D^2 u)-\ggg(D^2H) = \ggg(D^2u)-e^f \geq 0$.
For Case 2,   $\D_\bbM(u-H) = \ggg(D^2 u+I)-\ggg(D^2H+I) = \ggg(D^2u+I)-e^f \geq 0$.\qed

The assumption that $\G$ is complete is always necessary for Theorem \DD.4 to hold.
We  show this in the next  Example,  first for  the particular case where $\G$ is pulled back from a 
two-dimensional subspace of $\rn$ ($n\geq 3$), and then that it holds if  $\G$ is pulled
back from any proper subspace.

\noindent
{\bf  Example \DD.11.}  Consider the polynomial $\ggg(A) = a_{11}a_{22} - a_{12}^2$     on $\rn$.
This  is just the polynomial $\ggg(A) = \det(A\bigr|_V)$ where $V$ is the $(x_1,x_2)$-plane.
This polynomial is G\aa rding-Dirichlet, and the G\aa rding cone is just 
$$
\G=\{A : a_{11}a_{22}-a_{12}^2 >0 \ \ {\rm and}\ \ a_{11}+a_{22}>0\}.
$$
One computes that the polarization is
$$
\overline\ggg(G,A) \ =\ \half(g_{11}a_{22} + g_{22}a_{11} - 2 g_{12}a_{12}).,
$$
which, ignoring the $\half$, can be written as 
$$
\overline\ggg(G,A) \ =\ \left \langle
{\left(\begin{matrix} g_{22} & - g_{12} \\ -g_{12} & g_{11}
\end{matrix}\right), } 
{\left(\begin{matrix} a_{11} &  a_{12} \\ a_{12} & a_{22}
\end{matrix}\right) } 
\right\rangle
$$
Thus
$$
\overline\ggg(G,A) \ =\ \bra {{\mathbf M}(G)} {A}
$$
where 
$$
\text{${\mathbf M}(G)$ is the cofactor matrix of $G\bigr|_V$}.
$$
If $n\geq 3$ (and so $V$ is a proper subspace and $\gg, \cG$ is not complete), then  ${\mathbf M}(G)$ as an $n\times n$-matrix which has zero eigenvalues.  Hence ${\mathbf M}(G) \notin \Int \cp$  and Theorem \DD.4 fails. 

More generally,  consider the polynomial $\ggg(A) = \det(A\bigr|_V)$ where $V\ss \rn$
is a proper subspace of dimension $m<n$.  In this case 
$$
\overline\ggg(G, ... , G, A) \ =\ \bra {{\mathbf M}(G, ... , G)} {A}
$$
and
$$
\text{${\mathbf M}(G, ... , G)$ is the cofactor matrix of $G\bigr|_V$}.
$$
Note that ${\mathbf M}(G_1, ... , G_{N-1})$ is just the polarization of ${\mathbf M}(G, ... , G)$.
Again since ${\mathbf M}(G_1, ... , G_{N-1})$ is supported in $V$, it has zero eigenvalues
as an $n\times n$-matrix and hence does not belong to $\Int\cp$.  As before, this shows that 
Theorem \DD.4 fails whenever $\G$ is not complete. So Theorem \DD.4 cannot be strengthened.
(See Example \HH.3 for more about $\gg=\det$.)\ Finally note that Corollaries \DD.9 and \DD.10
 can all be combined into a uniqueness result.

\noindent
{\bf Theorem \DD.12.\ (Uniqueness).}  
{\sl
Suppose that $\G$ is complete and $X$ is connected.  Let  $u, H\in C^2(X)$ be a subsolution
and a solution respectively  to either of the two equations in Corollary \DD.14.
Then if $u-H$ has a maximum point in $X$,
$$
u\ =\ H+c \qquad\text{for some constant $c\in\bbr$.}
$$
Hence the subsolution $u$ must, in fact, be a solution.  

In particular, if $H_1$ and $H_2$ are two
solutions on $X$ with $H_1-H_2$ having either a maximum or a  minimum on $X$, 
then they differ by a constant.}

%%%%%%%%%%%%%%%%%%%%%%%%%%%%%%%%%%%%%%%%%%%%%%%%%%%%%%%%%%
%%%%%%%%%%%%%%%%%%%%%%%%%%%%%%%%%%%%%%%%%%%%%%%%%%%%%%%%%%
%%%%%%%%%%%%%%%%%%%%%%%%%%%%%%%%%%%%%%%%%%%%%%%%%%%%%%%%%%
%%%%%%%%%%%%%%%%%%%%%%%%%%%%%%%%%%%%%%%%%%%%%%%%%%%%%%%%%%
%%%%%%%%%%%%%%%%%%%%%%%%%%%%%%%%%%%%%%%%%%%%%%%%%%%%%%%%%%
%%%%%%%%%%%%%%%%%%%%%%%%%%%%%%%%%%%%%%%%%%%%%%%%%%%%%%%%%%
%%%%%%%%%%%%%%%%%%%%%%%%%%%%%%%%%%%%%%%%%%%%%%%%%%%%%%%%%%

%\vskip.3in

\vskip .5in

\centerline{\bf \headfont  \EE.  Derivatives.}
\medskip

  At a point $G\in\G$, it is very easy to compute the derivatives of $\log\,\gg$
(as opposed to $\gg$ or $\gg\ON$), because of  formula (6) for the $G$-eigenvalues of $A\in \Sn$ in Section 2:
$$
\gg(G+tA) \ =\ \gg(G) \prod_{j=1}^N (1+t \l^G_j(A)), \ \ \ {\rm or}
\eqno{(\EE.1)}
$$
$$
\log\, \gg(G+tA) \ =\ \log\, \gg (G) + \sum_{j=1}^N \log\,(1+t\l^G_j(A)).
\eqno{(\EE.2)}
$$
Computing the first and second derivatives at $t=0$ yields $D_G\log\,\gg$ and $D^2_G\log\,\gg$ as linear and quadratic 
forms in $A$ respectively.

\Prop{ \EE.1} 
$$
{\bf (1st Derivative)}\qquad
 \bra{D_G \log\,\gg}{A} \ =\ \sum_{j=1}^N  \l^G_j(A)
\eqno{(\EE.3)}
$$
{\sl which we also denote as either $\tr^G(A)$, the $G$-trace of $A$, or as $\bra{M(G,...,G)}{A} = \gg(G,...,G, A) = \D_M A$,
the $M=M(G,...,G)$-Laplacian of $A$ as in Section \DD.
$$
{\bf (2nd Derivative)}\ \  
(D_G^2 \log\,\gg)(A,A) \ =\ - \sum_{j=1}^N (\l^G_j(A))^2 \ \equiv\ -|\l^G(A)|^2.
\eqno{(\EE.4)}
$$
}

Continuing, the $k$th derivative of $\log \gg(G+tA)$ is easily computed, and setting $t=0$ yields $\pm$
the $k$th power function
$$
\left(D_G^{(k)}  (\log\, \gg)   \right) (A, ... , A) \ =\ (-1)^{k-1} \sum_{j=1}^N \left (\l_j^G(A)\right)^k.
$$

The formulas (\EE.3) and (\EE.4) can be used to prove 
$$
(D_G^2\gg)^{1\over N}(A) \ =\ - {\gg(G)^{1\over N} \over N} {\rm Discr} \l^G(A),
\eqno{(\EE.5)}
$$
which can be used to prove our new basic inequality for $\gg^{1\over N}$ in Section \LO.

For a rounded discussion of derivative formulas, including a proof of (\EE.5),  the reader is referred to [HL$_8$, App. D].

Note that by using the polarization ${\overline \gg} (A_1, ... , A_N)$ of $\gg$, we can compute the
polynomial 
$$\begin{aligned}
\vf(t) \equiv \gg(G+tA) &= {\overline\gg}(G+tA, ... , G+tA) \\
&= \sum_{k=0}^N \binom N k
{\overline\gg}(G,...,G,A,...,A)t^k \ \ \text{$A$ $k$-tiimes}.
\end{aligned}
\eqno{(\EE.6)}
$$
Expanding out (\EE.1) yields 
$$
\vf(t) \ =\ \sum_{k=0}^N \gg(G) \s_k(\l^G(A))t^k.
\eqno{(\EE.7)}
$$
This proves that
$$
\begin{aligned}
{1\over k!} (D_G^{(k)} \gg)(A, ... ,A) \ &= \ \gg(G)  \s_k(\l^G(A)) \ \\
&=\ \binom Nk 
{\overline\gg}(G,...,G,A,...,A)t^k \ \ \text{$A$ $k$-tiimes}.
\end{aligned}
\eqno{(\EE.8)}
$$

  The open polar $\G^*$ is the focus of the next Section \FF.  The key to Section \HH\  is that the gradient $D_A \log \, \gg \in \G^*$ if $A\in \G$ (Prop. \HH.1).  
  
  Finally, in this Section \EE, as an application of the first derivative formula (\EE.3), we prove the lemma needed for this Proposition \HH.1.
  This lemma  describes the maximum amount of "positivity" in $D_A \log \, \gg$.
  The example $\gg(D^2_x u) \equiv {\partial^2 u \over \partial x_1^2}$ in $\rn$, $n\geq 2$,
  shows that $D_A \log \, \gg = 2P_{[a_{11}]}$ need not be positive definite for a \GD operator.
  In fact, the additional hypothesis needed for positive definiteness is {\sl completeness}
  (see Prop. \CC.2 (6) and Cor. \HH.3).

\Lemma{\EE.2} {\sl
Fix $A\in \G$.  Then}
$$
\bra {D_A \log\, \gg}{B} \ >\ 0\quad \forall\, B\in \cG - E.
\eqno{(\EE.9)}
$$

\pf 
Since $A\in\G$ and $B\in \cG$ here, each $A$-eigenvalue of $B$ satisfies $\l^A_j(B)\geq 0$, $j=1, ... , N$.
Therefore, by the first derivative formula 
(\EE.3), we have $\bra {D_A \log\, \gg}{B} >0$ unless all the $A$-eigenvalues of $B$ vanish,
that is, unless $B$ is in the nullity of $\gg, \G$.  As mentioned in the \Ga Property (5) of Section 2, the nullity equals the edge, so
this proves (\EE.9). \qed

\Remark{\EE.3}  Interpreting (\EE.9) geometrically we have proved that:
$$
\begin{aligned}
&\ \ \ \text{$\log\,\gg$, or equivalently $\gg$, is strictly increasing} \\
&\text{along any curve in $\G$ whose tangents lie in $\cG-E$}
\end{aligned}
\eqno{(\EE.10)}
$$

 \Cor {\EE.4} {\sl
 If $A\in\G$, then 
 
 (1) \ \ $D_A \gg\geq 0$ and $D_A \gg \neq 0$.
 
 \noindent
 This implies that 
 
 (2) \ \ The level sets of $\gg$ in $\G$, namely $\{B\in\G : \gg(B)=c>0\}$, are smooth codimension one submanifolds of $\G$.
 }

\pf
$D_A\gg = \gg(A) D_A\log \gg$ and $\gg(A)>0$ for $A\in \G$.  If $B>0$, then $B\in \G \ss \cG-E$.
Thus $\bra {D_A \gg}{B} >0$ for all $B>0$, which implies  $D_A \gg\geq 0$ and $D_A\gg\neq 0$.\qed

%%%%%%%%%%%%%%%%%%%%%%%%%%%%%%%%%%%%%%%%%%%%%%%%%%%%%%%%%%%%%%%
%%%%%%%%%%%%%%%%%%%%%%%%%%%%%%%%%%%%%%%%%%%%%%%%%%%%%%%%%%%%%%%
%%%%%%%%%%%%%%%%%%%%%%%%%%%%%%%%%%%%%%%%%%%%%%%%%%%%%%%%%%%%%%%
%%%%%%%%%%%%%%%%%%%%%%%%%%%%%%%%%%%%%%%%%%%%%%%%%%%%%%%%%%%%%%%
%%%%%%%%%%%%%%%%%%%%%%%%%%%%%%%%%%%%%%%%%%%%%%%%%%%%%%%%%%%%%%%
%%%%%%%%%%%%%%%%%%%%%%%%%%%%%%%%%%%%%%%%%%%%%%%%%%%%%%%%%%%%%%%
%%%%%%%%%%%%%%%%%%%%%%%%%%%%%%%%%%%%%%%%%%%%%%%%%%%%%%%%%%%%%%%

%\vskip.3in

\vskip .5in

\centerline{\bf \headfont \FF.  The Relative Interior of the Polar Cone..}

\medskip

Here we describe some useful criteria for determining when a vector is in the relative interior of the polar of
a convex cone.  We start with a non-empty open convex cone  $\G\ss V$ in a finite dimensional 
real vector space $V$.  The polar
$$
\G^0\ \equiv \  \{y\in V^* : \bra yx \geq 0 \ \forall \, x\in \G\}
$$
is a non-empty closed convex cone.  As such it has interior when considered as a cone in its vector space $S \equiv {\rm span}(\G^0)\ss V^*$.  

\Def{\FF.1} We denote this relative interior by $\G^*$ and shall refer to it as the {\bf open polar}, keeping in mind
that it is an open convex cone  in its span $S$, which may be have  lower dimension than that of $V^*$.
It is easy to see that  the orthogonal complement $S^\perp$ is  the {\bf edge} $E$ of  $\G$, defined
to be $\overline\G \cap (-\overline\G)$ (see (6) in Section 2).
It is characterized as the linear subspace of $\G$ which contains all the lines (through the origin) in $\G$.

Note that  $\overline \G$ is {\bf self-polar} (i.e., $\G^0 = \overline \G^0 = \overline \G)$, if and only if $\G^*=\G$.
The classical examples of  self-polar cones are  $\cp$ and $\rn_{>0}$.
 \def\ccc{{\mathcal C}}
 
\noindent
{\bf Remark.}  Some caution is required as the set
$$
{\mathcal C} \ \equiv \ \{y\in S : \bra yx>0\ \ \forall \, x\in \G\}
$$
might seem a good  candidate for the open polar $\G^*$.  It is true that $\G^*\ss {\mathcal C} \ss \G^0$.
However,  first note that  for $\G \equiv \{x\in \bbr^2 : x_1>0, x_2>0\}$,   the polar
$\G^0 = \{x\in \bbr^2 : x_1\geq 0, x_2\geq 0\}$, that $\ccc$ equals $\G^0-\{0\}$,
and that $\G^*=\G$, so that $\G^* \neq {\mathcal C}$ and ${\mathcal C}\neq \G^0$.
Second, note that   for $\G\equiv \{x\in \bbr^2 : x_1>0\}$ where 
$\G^0= \{y\in \bbr^2 : y_1\geq 0, y_2 = 0\}$ is a closed ray with span $\bbr\times \{0\}$ and 
$\G^* =\{y\in \bbr\times \{0\} : y_1>0\}$ and here $\ccc = \G^*$.  So $\ccc$ does work as the open polar here.

Nevertheless,  the following gives  several criteria for $y$ to belong to the open polar $\G^*$, 
with (3) a particularly useful modification of the condition defining $\ccc$.

%\vfill\eject
\noindent
{\bf Lemma \FF.2. (The Open Polar Criteria).}    {\sl
Suppose $y\in S$.  The following are equivalent.

(1) \ \ $y \in \G^*$.

(2) \ \ $\exists\,  \e>0 \ \ \text{such that } \ \ \bra yx \geq \e|x|\ \ \forall\, x \in \overline\G \cap S$.

(3) \ \ $\bra yx >0 \quad \forall \, x\in \overline \G - E$.

(4) \ \ $\bra yx >0\quad \forall\, x\in (\partial \G-\{0\})\cap S$.

(5) \ \  $\bra yx >0\quad \forall\, x\in \overline \G \cap S, \ x\neq 0$.
}

 \noindent
 {\bf Proof.} (1) $\iff$ (2):  Note that:  
 $$
 \begin{aligned}
 y\in \G^* \ &\iff \ \text{there exist an $\e$-ball $B_\e(y) \ss \G^0$ about $y$} \\
 &\iff \exists \e>0 \ {\rm such\  that}\ \bra{y+\e z}x \geq 0 \ \forall\,|z|\leq 1, x\in \G \\
&\iff  \text{(a): $\exists \, \e>0 \  {\rm such\  that}\ \bra yx \geq \e \bra zx \ \ \forall\, |z| \leq 1, x\in\G$.}
  \end{aligned}
$$ 
Taking $z=x/|x|$ in (a) yields (2),  while (2) yields (a) since $\bra zx \leq |x|$ if $ |z|\leq 1$.

To see  (2) $\Rightarrow$ (3), suppose $x\in \overline\G - E$ and decompose $x$ into $x=a+b$ where
$a\in E, b\in S$.  Then $b\in \overline \G\cap S$ and $b\neq 0$, so $\bra yx = \bra bx >0$, which proves (3).

That (3) $\Rightarrow$ (4) is obvious.

Since $\overline \G$ is convex, it is easy to see that (4) $\Rightarrow$ (5).

To see  (5) $\Rightarrow$ (2): We can  assume that $E=\{0\}$ and $S=V$.  Obviously, (5) $\Rightarrow$ $y\in \G^0$.  Hence,
$$
\e \ \equiv \ \inf_{x\in \overline \G, |x|=1} \bra yx \ \geq\ 0.
$$ 
Since $\overline \G \cap \{|x|=1\}$ is compact, there exists $x_0 \in \overline \G \cap \{|x|=1\}$ with $\e= \bra y{x_0}$.
Now (5) $\Rightarrow$ $\e>0$.   Hence, $\bra y {{x\over |x|}} \geq \e \ \ \forall\, x\in \overline\G-\{0\}$, or 
$\bra yx \geq \e|x|\ \ \forall\, x\in \overline\G$.  \qed

This Lemma \FF.2 is needed in the next section.

%%%%%%%%%%%%%%%%%%%%%%%%%%%%%%%%%%%%%%%%%%%%%%%%%%%%%%%%%%%%%%%
%%%%%%%%%%%%%%%%%%%%%%%%%%%%%%%%%%%%%%%%%%%%%%%%%%%%%%%%%%%%%%%
%%%%%%%%%%%%%%%%%%%%%%%%%%%%%%%%%%%%%%%%%%%%%%%%%%%%%%%%%%%%%%%
%%%%%%%%%%%%%%%%%%%%%%%%%%%%%%%%%%%%%%%%%%%%%%%%%%%%%%%%%%%%%%%
%%%%%%%%%%%%%%%%%%%%%%%%%%%%%%%%%%%%%%%%%%%%%%%%%%%%%%%%%%%%%%%
%%%%%%%%%%%%%%%%%%%%%%%%%%%%%%%%%%%%%%%%%%%%%%%%%%%%%%%%%%%%%%%
%%%%%%%%%%%%%%%%%%%%%%%%%%%%%%%%%%%%%%%%%%%%%%%%%%%%%%%%%%%%%%%

\vfill\eject

%\vskip .5in

\centerline{\bf \headfont \HH.  The Gradient Diffeomorphism $D \log \, \gg$.}
\medskip

The map
$$
F(A) \ \equiv \ D_A(\log\, \gg) \ =\ {1\over \gg(A)} D_A \gg, \quad {\rm for} \ A\in\G
\eqno{(\HH.1)}
$$
given by  a \GD operator $\gg,\G$ is the focus of this section.  
A second important gradient map $D_A \gg\ON$ will be analyzed  in Section  \GM \  utilizing the \Ga dual function  $\gg^*$.

Now (\HH.1)  defines a smooth map $F:\G \to \Sn$.
In (\EE.3) of Proposition \EE.1, we computed that the element $F(A) = D_A \log\, \gg \in \Sn$,
considered as a linear functional on $\Sn$, with its natural inner product, is given by 
$$
\bra {D_A \log\, \gg}{A} \ =\ \sum_{j=1}^N \l^A_j(A), \qquad A\in \G \ \ {\rm and}\ \ A\in \Sn.
\eqno{(\HH.2)}
$$
Next we show that $F(A) \in \G^*$, the open polar of $\G$, if $A\in \G$.

\Prop{\HH.1} {\sl
The gradient map $D \log \, \gg$ takes its values in $\G^*$.
}

\pf
Fix $A\in \G$. One criterion for $F(A) \in\G^*$, given by (3) of Lemma \FF.2, is that 
$$
\bra {F(A)}{B} \ >\ 0\quad \forall\, B\in \cG - E.
\eqno{(\HH.3)}
$$
This is exactly what was proved in Lemma \EE.2. \qed

\Cor {\HH.2} {\sl    Suppose $A\in \G$.
If \ $\cG$ is complete, then $D_A \log \, \gg$ is positive definite.
}

\def\dlog{D_A \log \, \gg}
\noindent
{\bf Note.} It is easy to find examples that  show that the extra completeness hypothesis is required for the gradient to be positive definite.

\pf
Part (6) of the completeness criterion (Prop. \CC.2)
states that completeness is actually equivalent to the assertion $\G^* \ss \Int \cp$  and
hence $\dlog >0$ follows from $\dlog \in \Int \cp$ which was proved in Theorem \HH.1. \qed

We now continue with Example \DD.11.

\Ex{\HH.3. (The Determinant)}  The most basic example is $\gg = \det$, $\cG = \cp$.  Its gradient is its cofactor matrix, i.e., 
$$
D_A\det \ =\  A^{\rm cof}\ =\  \Mbf(A,...,A)  \quad \text{ (see (\DD.13)).}
\eqno{(\HH.3)}
$$
Since $A A^{\rm cof} = A^{\rm cof} A = (\det\, A)I$, this proves by (\HH.1) that 
$$
F(A) \ =\ D_A \log\, \det \ =\ {D_A \det \over \det\, A}\ =\ { A^{\rm cof}\over \det\, A}\ =              \ A\iv.
\eqno{(\HH.4)}
$$
Now $\G^*=\G\equiv \{P: P>0\}$ is self-polar, and hence the gradient diffeomorphism is just the matrix inverse map
$F(A) = A\iv$, which is a diffeomorphism of $ \{P: P>0\}$.
The complex and quaternionic cases are the same with $F(A) = A_\bbc\iv$ and $F(A)= A_\bbh\iv$ respectively.

  In general, $F(A) \equiv D_A \log\, \gg$ has some
remnants of an inverse.  For instance, the identity $\gg(A)\gg^*(F(A)) =1$, for all $A\in \G$, will be discussed in Section 10.

Proposition  \HH.1 frames the main theorem.  Note that since $\Sn$ is the orthogonal sum of the edge $E$ of $\G$
and the span $S$ of $\G^*$, we have the choice of considering either of the two isomorphic cones
$\G/E \cong \G\cap S$.
For our second gradient map $D\gg\ON$ considered in Section \GM, we choose $\G\cap S$.

\Theorem{\HH.4} {\sl
The gradient map
$$
F = D \log\,\gg :\G/E \ \arr\ \G^*
\eqno{(\HH.5)}
$$
is a diffeomorphism.  Moreover, the Jacobian  $F'(A) = D^2_A \log\,\gg$, at a point $A\in \G$, 
expressed as a quadratic form on $\Sn/E$, is the negative definite form}
$$
F'(A)(B,B) \ \equiv \ - \sum_{j=1}^N(\l^A_j(B))^2, \quad \forall\, B\in \Sn/E = S.
\eqno{(\HH.6)}
$$

\pf
Recall from (5) in Section 2 the the linearity 
$$
L\ \equiv\ \{B\in \Sn : \gg(A+B) = \gg(A) \ \ \forall\, A\in\Sn\}
$$
equals the edge $E$.  Hence, 
$$
(\log\,\gg)(A+B) = (\log\,\gg)(A) \ \text{is constant in $B$ for all $B\in E= L$.}
$$
This proves that
$$
F(A+B) \ \equiv\ D_{A+B}\log\,\gg \ =\ D_A\log\,\gg \ \equiv \ F(A).
$$
Thus the mapping $F:\G\cap S \to \G^*$ is well defined in (\HH.5).

At each point $A \in \G\cap S$, the Jacobian $F'(A) = D^2_A \log\, \gg$, which by formula (\HH.6) 
is a  negative definite quadratic form on $\G\cap S = \G/E$.  This proves that $F:\G/E \to \G^*$ is a local differmorphism.

To prove that $F$ is one-to-one, 
suppose $F(A_1) \equiv D_{A_1} \log\,\gg = D_{A_2} \log\,\gg \equiv F(A_2)$ for $A_1, A_2\in \G\cap S$.
If $A_1$ and $A_2$ are distinct, consider the affine line through them.  The one variable function $\vf$ equal 
to $\log\,\gg$ restricted to this line, has the same derivative at $A_1$ as at $A_2$. Since $\vf$ is strictly concave,
$A_1$ must equal $A_2$.

To prove that $F:\G/E \ \arr\ \G^*$ is onto we use  the Exhaustion Theorem in the next section. It says that for each $B\in\G^*$,
$
\bra BA - \log\, \gg(A)$ is a proper exhaustion function for $\G$.  In particular, it must have a minimum point $A\in \G$.
At the minimum point the derivative $B-\ D_A \log\,\gg$ must vanish, i.e., $B=F(A)$. Thus we have completed the proof of Theorem \HH.4. 
 \qed

%%%%%%%%%%%%%%%%%%%%%%%%%%%%%%%%%%%%%%%%%%%%%%%%%%%%%%%%%%%%%%%
%%%%%%%%%%%%%%%%%%%%%%%%%%%%%%%%%%%%%%%%%%%%%%%%%%%%%%%%%%%%%%%
%%%%%%%%%%%%%%%%%%%%%%%%%%%%%%%%%%%%%%%%%%%%%%%%%%%%%%%%%%%%%%%
%%%%%%%%%%%%%%%%%%%%%%%%%%%%%%%%%%%%%%%%%%%%%%%%%%%%%%%%%%%%%%%
%%%%%%%%%%%%%%%%%%%%%%%%%%%%%%%%%%%%%%%%%%%%%%%%%%%%%%%%%%%%%%%
%%%%%%%%%%%%%%%%%%%%%%%%%%%%%%%%%%%%%%%%%%%%%%%%%%%%%%%%%%%%%%%
%%%%%%%%%%%%%%%%%%%%%%%%%%%%%%%%%%%%%%%%%%%%%%%%%%%%%%%%%%%%%%%

%\vfill\eject

\vskip .5in

\centerline{\bf \headfont  \II. The Exhaustion Theorem.}

\medskip

The following theorem is essential for the main results in the coming sections.

\noindent
{\bf EXHAUSTION THEOREM \II.1}    {\sl
Suppose $\gg, \G$ is a \GD operator.  For each $B\in \G^*$, the open polar, the function 
$$
\psi(A) \ \equiv\ \bra BA - \log\,\gg(A), \qquad A\in \G\cap S
\eqno{(\II.1)}
$$
is a  exhaustion function for the G\aa rding cone $\G\cap S$, that is, the prelevel sets
$$
K_c \ \equiv\ \{A\in \G\cap S : \psi(A)\leq c\}, \qquad c\in\bbr,
$$
are compact, and exhaust $\G$, i.e., $$\bigcup_c K_c = \G.$$

Moreover $\psi$ is $C^\infty$ and strictly convex on $\G\cap S$.
}
\note{At some point we should think about changing $x$ to $\x$ or $A$.}

\pf
It suffices to prove the theorem when $\G$ is regular, since otherwise $\gg\bigr|_{S}$ has regular G\aa rding cone $\G_S = 
\G\cap S$ with the same open polar $\G^*_S = \G^*$.

Recall from (\EE.4) that $-\log\, \gg(x)$ has second derivative at a point $x\in \G$ given by
$$
\left\{  D_x^2(-\log \, \gg)   \right\} (\xi, \xi) \ =\ \sum_{j=1}^N\left ( \l_j^{\gg, x}(\x)   \right)^2
$$
for all $\x \in V$, where the G\aa rding eigenvalues are taken with respect to the direction $x\in \G$.
By G\aa rding [G\aa r] (or see  [HL$_4$]), the nullity set
$$
\{ \x \in V : \l_1^{\gg, x}(\x) = \cdots + \l_N^{\gg, x}(\x) = 0\}
$$
 equals the edge $E$.
Hence,  the function $-\log\, \gg(x)$ is strictly convex on $S=E^\perp$, and so $\psi(x)  \equiv \bra y x -\log\, \gg(x)$ 
has the same property  since 
$\bra y x$ is affine.  

Notice that since $\gg \equiv 0$ on $\partial \G$  and $\bra yx$ is finite,
the function $\psi \equiv +\infty$ on $\partial \G$.
We conclude that $K_c$ is a closed  subset of $\G$.

It remains to  show that $K_c$ is bounded.  For this we use the full hypothesis that $y$ belongs to the open polar  $\G^*$
of $\G$, which equals the interior of $\overline\G^0$.
By Lemma A.1 (2)  this is equivalent to the statement

\noindent (2) \ \ $\exists\ \e>0$  such that $\bra yx  \geq \e|x|\ \forall\,x\in \overline \G$.

This implies that:
$$
\begin{aligned}
K_c \ &\ss\ \{x\in  \overline \G : e^{\e|x|} \leq e^c \gg(x)\} \ss \left \{x\in  \overline \G : {(\e|x|)^{N+1} \over (N+1)!} \leq e^c\gg(x) \right \}  \\
&= \left \{x\in  \overline  \G : |x| \leq {(N+1)! e^c \over \e^{N+1}} \gg\left({x\over |x|}  \right)\right\} \ \ss\ B_R(0)
\end{aligned}
$$
with
$$\qquad
R \equiv {(N+1)! e^c\over \e^{N+1}} \sup_{\x\in \overline\G, |\x|=1} \gg(\x)\ <\ \infty.    \hskip.5in \mathqed
$$

%%%%%%%%%%%%%%%%%%%%%%%%%%%%%%%%%%%%%%%%%%%%%%%%%%%%%%%%%%%
%%%%%%%%%%%%%%%%%%%%%%%%%%%%%%%%%%%%%%%%%%%%%%%%%%%%%%%%%%%
%%%%%%%%%%%%%%%%%%%%%%%%%%%%%%%%%%%%%%%%%%%%%%%%%%%%%%%%%%%
%%%%%%%%%%%%%%%%%%%%%%%%%%%%%%%%%%%%%%%%%%%%%%%%%%%%%%%%%%%
%%%%%%%%%%%%%%%%%%%%%%%%%%%%%%%%%%%%%%%%%%%%%%%%%%%%%%%%%%%
%%%%%%%%%%%%%%%%%%%%%%%%%%%%%%%%%%%%%%%%%%%%%%%%%%%%%%%%%%%

%\vskip.3in

\vskip .5in

\centerline{\bf \headfont \JJ. Uniform Ellipticity and Auxiliary Operators.}
\medskip

 As an immediate consequence of Theorem \II.1  we have the following.
 
 \Prop {\JJ.1.}  {\sl
 For each $B\in\G$, the \GD operator  $\gg$ restricted to 
 $$
 K(\e, C) \  \equiv \  \{A\in \G : \gg(A)\geq \e\ \ {\rm and}\ \ \sum_{j=1}^N\l^B_j(A)\leq C\}
 $$
 with $\e>0$ and  $C>0$,  is uniformly elliptic.
 }

\pf
The operator $\gg$ restricted to $K(\e,C)$ being uniformly elliptic means that there exists $\eta>0$
such that the linearization $(D_A\gg)(B)$ at $A\in K(\e, C)$ has minimum eigenvalue $\l_{\rm min}(D_A\gg) \geq  \eta$
independent  of the point $A\in K(\e, C)$.  It suffices to show that the set of coefficient matrices $\{D_A\gg : A\in K(\e, C)\}$,
i.e., the image of $K(\e, C)$ under the gradient map $A\mapsto  D_A\gg$ is compact.  Finally, $K(\e, C) \ss K_{C-\e}$
(from the previous section) 
is immediate and easily implies that $K(\e, C)$ is compact.  \qed

\note{Maybe mention Phong's extension?}

%%%%%%%%%%%%%%%%%%%%%%%%%%%%%%%%%%%%%%%%%%%%%%%%%%%%%%%%%%%%%%%%%%%%%%%%%%%%%%%%%%%%%%%%%%%%%%%%%%%%%%%%%%%%%%%%%%%%%%%%%%%%%%%%%%%%%%%%%%%%%%%%%%%%%%%%%%%%%%%%%%%%%%%%%%%%%%%%%%%%%%%%%%%%%%%%%%%%%%%%%%%%%%%%%%%%%%%%%%%%%%%%%%%%%%%%%%%%%%%%%%%%%%%%%%%%%%%%%%%%%%%%%%%%%%%%%%%%%%%%%%%%%%%%%%%%%%%%%%%%%%%%%%%%%%%%%%%%%%%%%%%%%%%%%%%%%%%%%%%%%%%%%%%%%%%%%%%%%%%%%%%%%%%%%%%%%%%%%%%%%%%%%%%%%%%%%%%%%%%%%%%%%%%%%%%%%%%%%%%%%%%%%%%%%%%%%%%%%%%%%%%%%%%%%%%%%%%%%%%%%%%%%%%%%%%%%%%%%%%%%%%%%%%%%%%%%%%%%%%%%%%%%%%%%%%%%%%%%%%%%%%%%%%%%%%%%%%%%%%%%%%%%%%%%%%%%%%%%%%%%%%%%%%%%%%%%%%%%%%%%%%%%%%%%%%%%%%%%%%%%%%%%%%%%%%%%%%%%%%%%%%%%%%%%%%%%%%%%%%%

\def\E{E}
\def\fpsi{{F_f(\psi)}}
\def\bL{{\bf \Lambda}}
\def\bdf{{\bf f}}
\def\UU{U}
\def\bbm{{\bf M}}
\def\gg{{\mathfrak g}}

\vskip .5in

\centerline{\bf \headfont  \DU.   A  G\aa rding Duality Theory.}

\medskip

In this section we present one of the main topics of  this paper, 
  the duality theory for G\aa rding polynomials and their G\aa rding cones.
More precisely, for every G\aa rding polynomial $\gg$ on a vector space $V$ with G\aa rding cone $\G$,
there is a dual G\aa rding function $\gg^*$ defined on $\G^* \ss V^*$, where $\G^*$ is the
 relatively open polar of $\G$ in the dual space of $V$ (see  Section \FF). 
 We shall establish the following  G\aa rding-like results for $\gg^*$, although $\gg^*$ is not always a polynomial
 (see subsection \EX.6).
  As a result of this theorem, we shall refer to $\G^*$ as the {\bf G\aa rding-like cone associated to $\gg^*$}, and to $\gg^*,\G^*$ as a {\bf \\G\aa rdng-like pair}.
 
 \Theorem{\DU.1}  {\sl
 The function $\gg^*$ on $\G^*$ has the following properties.

 (1) \ \ $\gg^*$ is real analytic on $\G^*$.
 
 (2)  \ \ $\log\,\gg^*$ is strictly concave on $\G^*$.

 (3)  \ \ $\gg^*$ is homogeneous of  the same degree $N$ as $\gg$.

 (4)  \ \  {\bf (Monotonicity).} \  For $y, z \in \G^*$, $\gg^*(y+z)> \gg^*(y)$.  
 
  (5)  \ \  $\gg^*$ extends to a continuous function on the closure $\overline{\G}^*$ where 
 $$
 \gg\bigr|_{\partial \G^*} \ \equiv\ 0 \and \gg\bigr|_{\G^*} \  > \ 0.
 $$

 (6)  \ \ For $x\in\G$, the function $\psi(\yy) \equiv \bra \xx\yy - \log\, \gg^* (\yy)$ 
 is a strictly convex exhaustion function on $\G^*$.
  }

 %Note that Theorem \DU.1 says that the dual cone $\G^* $ ($\equdef$ the interior of the polar $\G^0$ of $\G$)
 %has all the properties of the \Ga\ cone of $\g^*$ were in fact a polynomial.

 Here the ordering (1) --- (6) reflects the order of the proof; it is different in the introduction.
 
As we have seen (Thm. \HH.3), the gradient map            
$$
F \equiv D \log\,\gg : \G/E \ \arr\ \G^*
$$
is a diffeomorphism, where $E$ is the edge of $\G$. Recall that the  {\bf edge} is  the  vector subspace 
 $E = \gg\cap (-\gg)$ of $V$, and $\G/E\ss V/E$.  As noted in Section \BB (5),  $E$ is also the nullity of $\gg$,
 which is important for this discussion.

   For clarity of exposition we shall  assume the edge $E = \{0\}$.  The reader can easily carry the exposition over to the general case above. One simply replaces  $V$ with $V/E$.

\Def{\DU.2} The {\bf dual \Ga function} $\gg^*$ is defined by 
$$
\gg^*(y) \ \equiv \ {1\over  \gg(F\iv (y))} \quad\text{for all $y\in \G^*$.}
\eqno{(\DU.1a)}
$$
or equivalently,
$$
\gg^*(F(x)) \ =\ {1\over \gg(x)} \qquad \text{for all $x\in\G$.}
\eqno{(\DU.1b)}
$$
The {\bf dual potential} $f^*$ of the potential $f= \log\, \gg$ is defined by
$$
f^*(y) \ \equiv \log\, \gg^*(y).
\eqno{(\DU.2)}
$$

We have the following three identities. 

\Prop{\DU.3. (The Duality Identities)}  {\sl  Suppose $x\in\G$ and $y\in\G^*$ correspond, that is, $y=F(x)$.}
Then
$$
\gg(x)\gg^*(y) \ =\ 1
\eqno{(\DU.3)}
$$ 
$$
f(x) + f^*(y) \ =\ 0
\eqno{(\DU.4)}
$$
$$
\bra xy \ \equiv  \ N.
\eqno{(\DU.5)}
$$

\pf
The first two are just the definitions.
The third is just Euler's formula for a degree $N$ homogeneous function $\gg$ on $\G$, $\bra {D_x\gg}{x} = N\gg(x)$, 
since $D_x(\log\, \gg) = D_x\gg / \gg(x)$.\qed

\Remark{\DU.4} We point out that 
$$
\bra xy \ \equiv \ N \ \ \text{\sl is actually equivalent to $\gg$ being homogeneous of degree $N$}
$$
as a straightforward computation shows.

Now we can begin the proof of Theorem \DU.1.

\noindent
{\bf Proof of (1):}
Note that $F(x) = D_x \gg/\gg(x)$ is a rational map, and hence its inverse $F\iv$ is real analytic,
which implies that  $\gg^* = 1/(\gg\circ F\iv$) is real analytic. \qed

\Theorem{\DU.5} {\sl The dual gradient map
$$
F^*(y) \ \equiv \ D_y f^*\ =\ D_y\, \log\, \gg^* \quad \forall \, y\in \G^*
$$
is the inverse of the gradient map}
$$ 
F(x) \ =\ D_x f \ =\ D_x \log\,\gg  \hskip.2in \forall \, x\in\G.
$$

\pf
Combining the identities (\DU.4) and (\DU.5) we have 
$$
f(x) + f^*(F(x)) \ =\ \bra x {F(x)} - N  \quad \forall\, x\in\G.
\eqno{(\DU.6)}
$$
Differentiating both sides yields
$$
F(x) + \left(D_{F(x)}f^*\right) F'(x) \ =\ F(x) + xF'(x)  \quad \forall\, x\in\G.
%\eqno{(\DU.7)}
$$
or
$$
\left(   D_{F(x)} f^* - x  \right) F'(x) \ =\ 0  \quad \forall\, x\in\G.
\eqno{(\DU.7)}
$$
where $F'(x)$ is the Jacobian of the gradient map $F$.  Since $F'(x)$ is negative definite,  it is invertible.  Hence
$$
D_{F(x)}f^*\ =\ x\quad \forall\, x\in\G,
\eqno{(\DU.8)}
$$
i.e., $F^*$ is the inverse of $F$,  \qed

\noindent
{\bf Proof of (2):}
Since $F:\G\to\G^*$ and $F^*:\G^*\to \G$ are inverses (Theorem \DU.5), the Jacobians 
$F'(x) \equiv D_x^2 \log\, \gg$ and $(F^*)'(y) = D_y^2 \log\,\gg^*$ are inverses (where $y=F(x)$ 
correspond).   Since $F'(x)$ is negative definite (Theorem \HH.3), its inverse $(F^*)'(y)$ is negative definite, which proves that $f^*=\log\,\gg^*$ is strictly concave.   \qed
\vskip.3in 

Next we address various homogeneity properties,  including the {\bf Proof of (3).} First,
$$
\text{$F(\xx) \equiv {1\over \gg(\xx)} D_\xx \gg$ is homogeneous of degree -1}.
\eqno{(\DU.9)}
$$
(By homogeneous of degree $k\in \bbr$ we mean homogeneous for positive scaling.)
This is obvious, as is the fact that $f(\xx) \equiv \log\, \gg(\xx)$ is {\bf  log-homogeneous with parameter $N$}, meaning that:
$$
f(t\xx) \ =\ f(\xx) + N\log\, t\qquad \text{for all $\xx\in\G$ and $t>0$.}
\eqno{(\DU.10)}
$$

If $G$ is a invertible map which is homogeneous of degree $k\neq 0$, then, since $x \mapsto \yy = G(\xx)$
implies $t\xx \mapsto \yy = t^k G(\xx)$, we see that, with $s=t^{1\over k}$, 
we have $s^{1\over k} \xx \mapsto s\yy$, so that $G^{-1}(sy) = s^{1\over k} G^{-1}(y)$.  
We conclude that
$$
\begin{aligned}
&\text{$G$ is homogeneous of degree $k\neq 0$} \\ 
\implies\    &\text{$G^{-1}$  is homogeneous of degree $1/k$.}
\end{aligned}
\eqno{(\DU.11)}
$$

This proves the following.

\Prop {\DU.6}  {\sl
Note that the degrees are the same in (\DU.11) when $k=1$ or $k=-1$.  In particular,
$$
F^{-1} \ =\ F^* \ =\ Df^*\ \ \ \text{is homogeneous of degree -1}.
\eqno{(\DU.12)}
$$
Therefore, by the identity $\gg(\xx)\gg^*(\yy) =1$
$$
\gg^*(\yy) \ =\ {1\over \gg(F^{-1}(\yy))} \ \ \text{is homogeneous of the same degree $N$ as $\gg$,}
\eqno{(\DU.13)}
$$
and, just as for $f$ its dual $f^*$ is log-homogeneous with parameter $N$.
$$
f^*(t\yy) \ =\ f^*(\yy) + N \log\, t \qquad\text{for all $y \in \G^*$ and $t>0$.}
\eqno{(\DU.14)}
$$
}
\vskip.3in

Before proving the rest of Theorem \DU.1, we give, as an application of the diffeomorphism
Theorem \HH.4 and parts (1), (2) and (3) of Theorem \DU.1, a "level set" diffeomorphism 
result, which is a companion to Theorem \HH.4.  Here we work in $\G\cap S$ instead of the 
isomorphic space $\G/E$.  We shall denote the level sets of $\gg$ and $\gg^*$ by:
$$
\G(c) \ \equiv \ \{x\in \G\cap S : \gg(x) = c\}, \quad \G^*(k) \ \equiv\ \{y\in \G^* : \gg^*(y) = k\}
\eqno{(\DU.15)}
$$ 
where  $c, k >0$.

By Corollary \EE.4 and Theorem \DU.5 we have:

\Remark{\. } {\sl
The gradients satisfy $D\gg \neq 0$ on $\G$ and $D\gg^*\neq 0$ on $\G^*$.  Therefore, the level sets $\G(c)$ and $\G^*(k)$ are smooth hypersurfaces.
}

\Theorem{\DU.7}  {\sl The diffeomorphism $D\log\, \gg : \G\cap S\cong \G^*$ of Theorem \HH.4,
restricts to a diffeomorphism 
$$
D \log\, \gg : \G(c) \ \overset  \cong \arr \ \G^*\left(  {1\over c} \right) \qquad \text{for each $c>0$},
\eqno{(\DU.16)}
$$
with inverse
$$
D \log\, \gg^* : \G^*(k) \ \overset  \cong \arr \ \G\left(  {1\over k} \right) \qquad \text{for each $k>0$}.
\eqno{(\DU.17)}
$$
 }

\pf 
By Definition \DU.1 we have $\gg(x) \gg^*(D_x \log\, \gg) = 1$ for $x\in \G$.
Therefore, $x\in \G(c) \iff y \equiv  D_x \log \, \gg \in \G^*(1/c)$.
Applying the diffeomorphism in Theorem \HH.4 completes the proof.\qed

Two more diffeomorphism theorems using the gradient maps $D\gg$ and $D\gg\ON$   will be given in the next Section \GM.  

\vskip.3in

\noindent
{\bf Proof of (4).}   Since $\G^*$ is convex,  we can apply the  Mean Value Theorem to the interval $[y, y+z]$ reducing the proof 
to showing  that the directional derivative  of $\gg^*$  at $\bar y\in [y, y+z] \ss \G^*$ in the direction $z\in \G^*$ is positive, namely 
$$
\bra {D_{\bar y} \gg^*}{z} \ >\ 0   \quad\forall \, {\bar y},z\in \G^*.
\eqno{(\DU.18)}
$$
Equivalently, 
$$
\bra {D_{\bar y} \log\, \gg^*}{z} \ >\ 0 \quad\forall \, \bar y,z\in \G^*.
\eqno{(\DU.19)}
$$
since 
$$
D_{\bar y} \log\, \gg^*  \ =\ {D_{\bar y} \gg^*      \over  \gg^*({\bar y})    }  \and \gg^*({\bar y})>0 \qquad \forall\, {\bar y}\in \G^*.
$$
Now by Theorem \DU.5, $D_{\bar y}(\log\, \gg^*) = F\iv({\bar y})$ the inverse of the gradient map $F(x) \equiv D_x \log\, \gg$ which maps 
$\G$ to $\G^*$.  Since $D_{\bar y}\log\, \gg^* \in \G$ and $z\in \G^*$, the open polar, we have that (\DU.19) holds.
\qed

Before addressing the boundary behavior of  $\gg^*$, which proves (5), we need the infimum characterization
of the dual potential.

\Theorem {\DU.8} {\sl The {\bf dual potential}  $f^*(\yy) = \log\,\gg^*(y)$ of the potential function
 $f(\xx) \equiv \log \, \gg(\xx)$ is given by}
 $$
f^*(\yy) \ \equiv \ \inf_{\xx\in\G}  (\bra\yy\xx  -f(\xx) - N)  \quad {\rm for}\ \ y\in\G^*
\eqno{(\DU.20)}
$$

\pf
By the Exhaustion Theorem \II.1,   $\bra  \yy\xx -f(\xx) - N$ is a strictly convex exhaustion function for $\G$, so that it has 
a unique minimum point $\xx\in \G$. This minimum point $\xx$ is critical for the exhaustion function, i.e., 
$$
\yy  -  D_\xx f\  = \ 0.
\eqno{(\DU.21)}
$$
Evaluating the exhaustion function at the minimum $\xx$ yields the minimum value $f^*(\yy)$:
$$
f^*(\yy) \ =\ \bra \yy\xx - f(\xx) -N \ \equiv  \ -f(\xx) 
\eqno{(\DU.22)}
$$
where the last equality is by the  identity $\bra yx=N$. \qed

\vskip.3in

\noindent
{\bf Proof of  Property (5) in Theorem \DU.1} which, for clarity,  we restate here.

%In order to apply  the Exhaustion Theorem \II.1 to $\gg^*$ instead 
%of $\gg$, and conclude that the double dual of $f$ (or $\gg$) equals 
%$f$ (or $\gg$), we need some further properties of the dual G\aa rding function $\gg^*$.

\Prop {\DU.9} {\sl
Extend $\gg^* \in C^\infty(\G^*)$ to a function $\gg^* : \overline{\G^*}  \to \bbr$ by defining $\gg^*$ to be $\equiv 0$
on $\partial\G^*$.  Then $\gg^* \in C(\overline{\G}^*)$, and of course, $\gg^*>0$ on $\G^*$ and $\gg^*\bigr|_{\partial \G^*} =0$.
}

\pf  For all $\yy\in \G^*$ we have $\gg^*(\yy) = 1/ \gg(F^{-1}(\yy))$.
By (1), $\gg^*(\yy)$ is real analytic on $\G^*$ Also   $\gg^*(\yy) >0$ for all $\yy\in \G^*$, and $\gg^*\bigr|_{\partial \G^*} =0$ by definition.

It remains to prove continuity at boundary points. 
For this we show that if $\{y_j\}_j$ is a sequence in $\G^*$ converging to a boundary point $y \in \G^*$,
then $\lim_{j\to\infty} \gg^*(y_j) =0$,
by proving that   $\lim_{j\to\infty} \log\, \gg^*(y_j) = \lim_{j\to\infty} f^*(y_j) = -\infty$.

Since $y\in \overline{\G}^*= \G^0$ but $\yy\notin\G^*$, condition (1) in the Open Polar Criterion Lemma is false.
Hence Condition (4) is false, i.e., 
$$
\exists \, \xx_0 \in \partial\G - \{0\} \ \ \ \text{such that}  \ \bra\yy{\xx_0} = 0.
$$

\def\eb{\zz}

Recall that $\G+ \overline \G=\G$, since it is an open subset of $\overline \G$ containing $\G$.  In particular, 
for all $\zz\in \G$ the ray $\zz +t\xx_0, t>0$ belongs to $\G$.
Setting $\xx \equiv \zz + t \xx_0$   in (\DU.20) in Theorem \DU.8  yields
$$
f^*(\yy_j) \ \leq \ \bra {\yy_j}{\zz + t\xx_0} -f(\zz+t\xx_0) - N.
$$
Now 
$$\lim_{j\to\infty} \bra{y_j}{\zz + t\xx_0} = \bra \yy \zz - t\bra \yy {\xx_0} = \bra \yy \zz
$$
since our assumption is that $\bra \yy {\xx_0} =0$.  However,
$$
\lim_{t\to\infty}  f(\zz+t\xx_0) = f(\zz) + \lim_{t\to \infty} \sum_j \log\, (1+t\l_j^\zz(\xx_0))
$$
by the basic formula (7) in Section 2 that  we used to compute the derivatives of $\log\,\gg$ at $z$ in terms
 of the $z$-eigenvalues. 
 %This gives
 %$$
 %\gg(\eb+t\xx_0) \ =\ \gg(\eb) \prod_j (1+t\l^j_{\eb}(\xx_0)).
 %\eqno{(\DU.21)}
%$$
Since  $\l_j^zz(\xx_0)\geq 0$ for all $j$, we see that unless all 
$\l^j_{\eb}(\xx_0)= 0$, we have $\lim_{t\to\infty} f(\zz+t\xx_0)= \infty$.  \qed

\vskip.3in

\noindent

Now  we prove {\bf The Dual Exhaustion Lemma}, which is the remaining  assertion (6) of Theorem \DU.1,  namely:

\Lemma{\DU.10.  (The Exhaustion Lemma for $\gg^*$)}  {\sl
Suppose $\gg$ is a G\aa rding polynomial of degree $N$ on $V$ with open G\aa rding cone $\G$.
Fix $x\in \G$.  Then the function
$$
\psi(y)\ \equiv \ \bra \xx\yy - \log\,\gg^*(\yy)
$$
is a strictly convex exhaustion function  for $\G^*$, so that its prelevel sets 
$$
K ^*_c \ \equiv\ \{\yy\in \G^*  : \bra \xx\yy - \log\, \gg^*(\yy) \leq c\}  \qquad\text{for all $c\in\bbr$}
$$
are compact.
}

\pf
The function $\psi$ is strictly convex since $f^*(\yy) = \log\, \gg^*(\yy)$ is strictly concave on $\G^*$.
The proof of the ExhaustionTheorem  \II.1  carries over with $\gg^*$ replacing $\gg$ since
$$
\gg^* \in C(\overline{\G}^*), \ \  \gg^*\bigr|_{\partial \G^*}  = 0, \ \ \gg^*\bigr|_{\G^*} >0, \ \ {\rm and} 
\ \ \gg^*\ \text{is homogeneous of degree $N$}.
$$
\qed

\Cor{\DU.11} {\sl
For constants $\e, c>0$ and $x\in\G$, the sets 
$$
K_{\e,c} \ \equiv \ \{\yy\in \G^* : \bra \xx\yy \leq c\ \ {\rm and} \ \ \gg^*(\yy)\geq \e\}
$$
are compact.
}

\pf
The proof of Corollary C.3 in [HL$_8$] carries over. \qed

 The properties of the dual \Ga function $\gg^*$  established in this section are sufficiently G\aa rding-like to apply our dual construction to obtain 
 a double dual function , even though the dual function $\gg^*$ may not be a polynomial.

 Suppose that $\gg, \G$ are a pair  as  the pair $\gg^*, \G^*$ in Theorem \DU.1, with all the properties stated there.
 We will call such a pair a {\bf G\aa rding-like pair}.  
For such a pair the  arguments given in proving  Theorem \DU.1  can be applied  to construct a new operator and domain
  $\gg^{**}, \G^{**}$.

 \Theorem{\DU.12} {\sl  If $\gg,\G$ is a G\aa rding-like pair on a vector space $V$, then $$\gg^{**} = \gg.$$  
 In particular,
 if $\gg$ is a \Ga polynomial with \Ga cone $\G$ on $\Sn$, then
 $$
\gg \ =\ \gg^{**} \and \G \ =\ \G^{**}.
 $$
 }
 
 \pf  By the Bipolar Theorem $\G^{**}=\G$. Since $D\, \log\,\gg^*$ is the inverse of $D\,\log\, \gg$ by Theorem 10.5, if we apply Definition 10.2 with $\gg$ and $\gg^*$ interchanged, it follows that $\gg^{**}= \gg$.\qed

 %%%%%%%%%%%%%%%%%%%%%%%%%%%%%%%%%%%%%%%%%%%%
%%%%%%%%%%%%%%%%%%%%%%%%%%%%%%%%%%%%%%%%%%%%
%%%%%%%%%%%%%%%%%%%%%%%%%%%%%%%%%%%%%%%%%%%%
%%%%%%%%%%%%%%%%%%%%%%%%%%%%%%%%%%%%%%%%%%%%
%%%%%%%%%%%%%%%%%%%%%%%%%%%%%%%%%%%%%%%%%%%%
%%%%%%%%%%%%%%%%%%%%%%%%%%%%%%%%%%%%%%%%%%%%
%%%%%%%%%%%%%%%%%%%%%%%%%%%%%%%%%%%%%%%%%%%%
%%%%%%%%%%%%%%%%%%%%%%%%%%%%%%%%%%%%%%%%%%%%
%%%%%%%%%%%%%%%%%%%%%%%%%%%%%%%%%%%%%%%%%%%%
%%%%%%%%%%%%%%%%%%%%%%%%%%%%%%%%%%%%%%%%%%%%
%%%%%%%%%%%%%%%%%%%%%%%%%%%%%%%%%%%%%%%%%%%%

\vskip .5in

\centerline{\bf \headfont  \GM.  The Gradient Maps $D\gg$ and $D\gg\ON$.}

\medskip

 Both are rescalings of the gradient map which sends $A \in  S(n)$ to $D_A \log\, \gg$ in Section \HH.
 They are given by
 $$
 D_A \gg \ =\ \gg(A) D_A \log\, \gg, \qquad A\in S(n)
 \eqno{(\GM.1)}
 $$
 and
  $$
  D_A \gg\ON \ =\ \gg(A)\ON D_A \log\, \gg, \qquad A\in S(n).
 \eqno{(\GM.2)}
 $$
As such, most of their properties are easy to derive from Section \HH, but we shall also use the duality
of the last section \DU.  Note that they are positively homogeneous with:
 $$
\text{$D_A\gg$ degree $N-1$,  $D_A \gg\ON$  degree 0, and $D_A\log \,\gg$  degree -1.}
 \eqno{(\GM.3)}
 $$
 Here we elect to consider $D_A\gg : \G/E \to \G^*$ as a map defined on $\G/E$, while
 for $D \gg\ON$ we will consider it as a map  defined on $(\G\cap S)/\bbr_{>0}$.
 Note that each level set $\G(c) \equiv \{A\in \G : \gg(A)=c\}$ ($c>0$),  is diffeomorphic to 
 $(\G\cap S)/\bbr_{>0}$, so that each $\G(c)\cap S$ can also be taken as the domain of the gradient map $D\gg\ON$.
 
 \Theorem{\GM.1} {\sl
 Suppose $\gg, \overline \G$ is a \GD operator subequation pair on $\rn$ of degree $N$.
 
 \noindent
 (a) The gradient map
 $$
 D\gg : \G/E \ \to \ \G^* \ \ \text{is a diffeomorphism}.
 \eqno{(\GM.4)}
  $$
 (b) The gradient map
  $$
  D\gg\ON : \G(c)\cap S \cong (\G\cap S)/\bbr_{>0} \ \to\ \G^*(1)  \ \ \text{is a diffeomorphism}.
 \eqno{(\GM.5)}
 $$
 }
 
 \noindent
 {\bf Proof of (a):} Let $G(A) \equiv D_A\gg$ denote the gradient map $D_A\gg$, and let $F(A)$ denote the gradient map $D_A\log\,\gg$, so that $G(A) = \gg(A)F(A)$.
 First we  prove $G$ is onto $\G^*$.  Given $B\in\G^*$, pick $A\in \G/E$ with $F(A)=B$ by Theorem \HH.4. Then $G(\gg(A)^{-1 \over N-1} A) = G(A)/\gg(A) = F(A) = B$,
 since $G$ is positively homogeneous of degree $N-1$.
 
 To prove that $G:\G/E\to \G^*$ is one-to-one, assume $G(A_1)=G(A_2)$ with $A_1, A_2 \in \G$.
 We use the dual function $\gg^*$ on $\G^*$ as follows.  Since $\gg^*$ is positively homogeneous
 of degree $N$, we have that $\gg^*(G(A_i)) = \gg^*(\gg(A_i) F(A_i)) = \gg(A_i)^N \gg^*(F(A_i))
 = \gg(A_i)^N / \gg(A_i) = \gg(A_i)^{N-1}$,
since $\gg(A_i)\gg^*(F(A_i))=1$.
 Hence $G(A_1)=G(A_2) \Rightarrow \gg(A_1)^{N-1} = \gg(A_2)^{N-1} \Rightarrow\gg(A_1)=\gg(A_2)$, and therefore $F(A_1)=F(A_2)$.
 By Theorem \HH.4 this implies that $A_1-A_2 \in E$. \qed

 \noindent
 {\bf Proof of (b):} By Theorem \DU.7, $D \log\, \gg : \G(1) \to \G^*(1)$ is a diffeomorphism.
 At points $A\in \G(1)$ we have, since $\gg(A)=1$, that $D_A \gg\ON = \gg(A)\ON D_A \log\, \gg = D_A \log\, \gg$.  Thus the diffeomorphism $D\log \,\gg : \G(1) \to \G^*(1)$ equals the desired 
  diffeomorphism $D\gg\ON : \G(1) \to \G^*(1)$.\qed
  
  \Remark {\GM.2. (Jacobians of the Gradient Map)} 
 For the convenience of the reader we list the Jacobians of the three gradient maps
 $D\log\, \gg, D \gg\ON$ and $D\gg$ at a  point $B\in\G$  as quadratic forms utilizing the 
 $B$-eigenvalues $\l_j^B(\x)$ of $\x \in S(n)$, $j=1, ... , N$.
 See Appendix D of [DetMajCalcVar] for more details.
 
 \Prop{\GM.3. ([HL$_8$])}
 
 {\sl
 (1)  The gradient map $F(B) \equiv D_B \log \,\gg$ has Jacobian $F'(B) = D^2_B \log \, \gg$ at $B\in \G$, given by
 $$
 F'(B)(\x,\x) \ =\ -\left |  \l^B(\x)   \right|^2.
 \eqno{(\GM.6)}
 $$
     
     (2)  The gradient map  $G(B) = D_B\gg$ has Jacobian $G'(B)= D^2_B \gg$, at $B\in \G$, given by
  $$
 G'(B)(\x,\x) \ =\  2\gg(B) \s_2(\l^B(\x)).
  \eqno{(\GM.7)}
 $$   '
     
     (3) The gradient map $H(B) =D_B\gg\ON$ has Jacobian  $H'(C) = D^2_BH$ at $B\in\G$, given by
      $$
 H'(B)(\x,\x) \ =\  - {1\over N^2} \gg(B)\ON {\rm Discr}(\l^B(\x)).
  \eqno{(\GM.8)}
 $$    
     
     \noindent
     Here $\s_2(\l)  \equiv \sum_{i<j} \l_i\l_j$ and {\rm Discr}$(\l) \equiv \sum_{i<j} (\l_i-\l_j)^2$.
     }
     
     \noindent
     The null spaces of $D^2_B \log\,\gg$ and $D_B\gg$ are the same, 
     namely the nullity of $\gg$ $\{\x\in S(n) : \l^B_j(\x) = 0, \, j=1, ..., N\}$,  which equals the 
     edge $E$ of $\G$.  Hence, restricting these quadratic forms to $S=E^\perp = S(n)/E$,
     the quadratic form $F'(B) = D^2_B\log\,\gg$ is negative definite, and
     the quadratic form $G'(B)= D^2_B\gg$ is of Lorentzian type, positive in the radial direction and negative definite in the perpendicular hyperplane, (i.e., the future light cone is $\G$.)
     Finally $H'(B)= D^2_B \gg\ON\geq 0$ with null space $E+\bbr\cdot  B$ of dimension one more than the edge.

%%%%%%%%%%%%%%%%%%%%%%%%%%%%%%%%%%%%%%%%%%%%%%%%%%%%%%%%%%%%%%%%%%%%%%%%%%%%%%%%%%%%%%%%%%%%%%%%%%%%%%%%%%%%%%%%%%%%%%%%%%%%%%%%%%%%%%%%%%%%%%%%%%%%%%%%%%%%%%%%%%%%%%%%%%%%%%%%%%%%%%%%%%%%%%%%%%%%%%%%%%%%%%%%%%%%%%%%%%%%%%%%%%%%%%%%%%%%%%%%%%%%%%%%%%%%%%%%%%%%%%%%%%%%%%%%%%%%%%%%%%%%%%%%%%%%%%%%%%%%%%%%%%%%%%%%%%%%%%%%%%%%%%%%%%%%%%%%%%%%%%%%%%%%%%%%%%%%%%%%%%%%%%%%%%%%%%%%%%%%%%%%%%%%%%%%%%%%%%%%%%%%%%%%%%%%%%%%%%%%%%%%%%%%%%%%%%%%%%%%%%%%

%\vskip.3in

\def\ON{^{1\over N}}
\def\On{^{1\over n}}

\vskip .5in

\centerline{\bf \headfont  \LO. The Duality Version of  G\aa rding's }
\centerline{\bf \headfont   Basic Inequality for the Operator $\gg(A)^{1\over N}$.}

\medskip

We begin by discussing the fundamental classical case of the operator $\det(A)$ on $\rn$.
Let  $\bra \cdot \cdot$ denote  the natural inner product on ${\cal S}(n)$, the space of real symmetric $n\times n$ matrices, and denote these matrices by $A, B, C, ...$  etc.
 
 We start with the symmetric form for two positive definite matrices $A,B$.
 
 \noindent
 (1) For $A,B >0$
 $$
\begin{aligned}
& (\det A)\On (\det B)\On  \ \leq \ {1\over n} \bra BA   \\
\text{\sl with}\ &\text{\sl  equality $\iff \ A=\l B^{-1}$ for some $\l>0$}.
\end{aligned}
\eqno{(\LO.1a)}
$$

Restricting to $\det B =1$, $B>0$ yields an equivalent form.

\noindent
(2) For $A>0$
$$
\begin{aligned}
(\det A)\On \ &\leq \ {1\over n} \bra BA \quad \forall\, B>0,\ \  \det B = 1, \ {\rm and}   \\
 \text{\sl Equality holds} & \iff \text{\sl  for $A>0$ one has $B = (\det \, A)\On A^{-1}$}.
 \end{aligned}
 \eqno{(\LO.1b)}
 $$
We will also find the degenerate case useful (cf. Proof of Proposition\ \Ph.1).

\noindent
(3) For $A\geq 0$, possibly singular,
$$ \begin{aligned}
   \text{$\det(A)^{1\over n}=   \underset {\det(B)=1}{ \inf_{{B>0 }}} {1\over n} \bra BA $.   }
\end{aligned}
 \eqno{(\LO.1c)}
$$

\noindent
{\bf Proof of (1):}
One can reduce to the standard Geometric Mean -- Arithmetic Mean (GM-AM) inequality as follows.
Set $P\equiv \sqrt A B \sqrt A$. Then $P$ is obviously symmetric, and $P \geq 0$ since
 $\bra{\sqrt A B \sqrt Ax} x = \bra { \sqrt B \sqrt Ax} { \sqrt B \sqrt Ax} \geq 0$ for all $x$. 

Diagonalizing   $P \equiv \sqrt A B \sqrt A$ with eigenvalues $\l_1, ... , \l_n$,  one sees that 
if  $P>0$ the inequality 
$$
\det(P)\On \leq {1\over n} \bra I {P}
 \eqno{(\LO.2a)}
$$
becomes   the classical GM-AM inequality 
$$
\text{(GM-AM)} \qquad   \qquad
{(\l_1 \cdots \l_n)}\On \leq {1\over n} (\l_1+ \cdots + \l_n).  \qquad\qquad
 \eqno{(\LO.2b)}
$$ 
Equality holds if and only if all the eigenvalues of $P$ are equal, i.e., $P=\l I, \l>0$.
(Obviously the same inequality  holds for $P\geq 0$, but if $P \not >0$, the left hand side is zero, so
equality holds only if $P=0$.)

Now   $\det\, P = \det(\sqrt A B \sqrt A) = \det(AB) = ( \det\, A)( \det\, B)$ and 
$\bra BA  = \bra  B {\sqrt A \sqrt A} = \bra  {\sqrt A B \sqrt A} I = \tr \, P$.
So the inequality (\LO.2a)    for $P$ and  the (GM-AM)  inequality in (\LO.2b) are the same.
For the  equality in (\LO.1a),
note that  the eigenvalues of $P$ are all equal to $\l$ means that $P=\l I = \sqrt A B \sqrt A$, or equivalently, $B=\l A^{-1}$.
 \qed
 
 \noindent
 {\bf Proof that (1)  $\Rightarrow $ (2):}
 The inequality in (2) is just the special case of (1) where the determinant of $B$ equals 1.
 For the equality in (2) note that for $\l>0$,  we have  $A = \l B^{-1} \iff B= \l A^{-1}$.
 In order for $1=\det B = \l^n(\det A)^{-1}$, $\l$ must equal $(\det A)\On$. \qed

 \noindent
 {\bf Proof of (3):}  Using  (2) it suffices to find $B_\e>0$ with $\det B_\e =1$ and $\bra {B_\e}A$ arbitrarily small.  Let $N$ denote the null-space of $A\geq 0$.  The restriction of $A$
 to $N^\perp$ is positive definite.  Given $\e>0$, set $B_\e =\e P_{N^\perp}+\l(\e) P_N$
 with $\l(\e)$ chosen (to be the right power of $1/\e$) so that $\det B_\e=1$.  
 Now  $\bra {B_\e}A = \e \bra {P_{N^\perp} }A$ is independent of $\l(\e)$ 
 so it has  limit zero as $\e\searrow0$.\qed

One purpose of this paper is to generalize this
 classical inequality, with determinant and the open cone $\{ P>0\}$ replaced by 
any G\aa rding polynomial $\gg$ and its open G\aa rding cone $\G$.
(The cone $\G$ for $\gg = \det$ is  $\{ P>0\}$.)
The key new ingredient is the introduction  of a dual G\aa rding function $\gg^*$
defined on the open dual cone $\G^*$ of $\G$. 
(Recall that $\gg = \det$ is self-dual, that is, $\gg^*=\gg$ and  $\G^*=\G$.)

We suppose that $\gg$ is a G\aa rding operator of degree $N$  with G\aa rding cone $\G$
on a real vector space $V$.  Let $\G^*$ denote the open polar discussed in Section 6, and let
$\gg^*$ denote the dual G\aa rding function defined on  $\G^*$ by (\DU.1).
Just as with the three  versions (1), (2) and (3) above, we state our result in three forms.

\Theorem {\LO.1}  \ \ {\sl 

\noindent
(1) For  $A\in  \G$ and $B\in \G^*$, 
$$
\gg(A)\ON \gg^*(B)\ON \ \leq \ {1\over N}\bra AB
\eqno{(\LO.3a)}
$$
with equality  $\iff$ $B \equiv \l\, D_A \log \,\gg$  (modulo the edge $E$ of $\G$), for some $\l>0$,
in which case $\bra AB = \l N$.

\noindent
(2) Restricting to the level set $\gg^*(B)=1$ yields
$$
\gg(A)\ON \ \leq\ {1\over N} \bra BA \quad \forall\, A \in \G \ \ {\rm and}\ \ \forall \ B\in \G^* \ \ {\rm with}\ \ \gg^*(B)=1
 \eqno{(\LO.3b)}
$$
with equality if and only if $B \equiv D_A \gg\ON$ module the edge $E$.

The singular case $A\in\partial \G$ is included in the final statement.

\noindent
(3) For $A\in \overline \G$, one has
$$
\gg(A)\ON \ \ = \ \  \underset {\gg^*(B)=1}{ \inf_{{B\in \G^* }}} {1\over N} \bra BA.
\eqno{(\LO.3c)}
$$
}

\noindent
{\bf Proof of (1):}  Given $B\in \G^*$, the gradient diffeormorphism Theorem \HH.4 says that there exists a unique $C\in \G$ with $D_C\log\, \gg = B$.  Taking $C$-eigenvalues, $\l^C_j(A)$
the inequality (1) reduces to the classical GM-AM for $\l^C(A)$ as follows.  
First, the right hand side of  (1) is the arithmetic mean of $\l^C(A)$
$$
{1\over N}\bra AB \ =\ {1 \over N}\bra {D_C \log\, \gg}{A} \ =\ {1\over N} \sum_{j=1}^N \l^C_j(A)
$$
by the first derivative formula (Prop.\ \EE.1 (\EE.3)). Second, regarding the left hand side of 
(1), we have $\gg(A)\equiv \gg(C) \prod_{j=1}^N \l^C_j(A)$, while $\gg^*(B) = \gg^*(D_C\log \,\gg)$<
which by Definition \DU.2 of $\gg^*$, equals $1/\gg(C)$.  
Hence, $\gg(A)\gg^*(B) = \prod_{j=1}^N \l^C_j(A)$, 
so the LHS is the geometric mean of $\l^C(A)$.

For the equality part in (1), note first that by the classical GM-AM inequality, all the $C$-eigenvalues of $A$ must be equal, i.e., $\l^C_j(A)=\l>0$ for $j=1, ... ,N$. 
Now $\l^C_j(C)=1$, so equivalently $\l^C_j(A-\l C)=0 \ \forall\, j$.
Hence, $A-\l C$ belongs to the nullity which, as noted in (5) of Section \BB, equals the edge.
Thus $D_A\log\,\gg= D_{\l C}\log\, \gg = \l^{-1} D_C\log\, \gg = \l^{-1} B$
by homogeneity -1 of $A\mapsto D_A\log\,\gg$.  This proves that $B=\l D_A\log\,\gg$.
Finally, $\bra AB = \l \bra {A}{D_A\log\,\gg} = \l \bra{A}{D_A\gg} / \gg(A) = \l N$ by Euler's Formula. \qed

\noindent 
{\bf Proof that (1) $\Rightarrow$ (2):}  Just restrict (1) to $B\in \G^*$ with $\gg^*(B)=1$.  Now in the case of equality, we have $\l= \gg(A)\ON$ because $1=\gg^*(B) = \gg^*(\l D_A \log\, \gg)
=    \l^N\gg^*(D_A \log\,\gg) = \l^N /\gg(A)$.  Therefore, $B=\gg(A)\ON D_A\log\, \gg = D_A\gg\ON$.\qed

\noindent
{\bf Proof of (3):}  If $A\in\G$, then(2) (\GM.3b) implies that 
$$
\gg(A)\ON \leq   \underset {\gg^*(B)=1}{ \inf_{{B\in \G^* }}}     {1\over n}\bra BA.
\eqno{(\LO.4)}
$$
Equality holds at the unique point $B \equiv D_A\gg\ON \in \G^*$ by (2).
If $A\in \partial \G$, then approximate $A\equiv\lim_j A_j$ with $A_j \in \G$ and take the 
limit as $j\to \infty$ to get (\LO.4).  Taking $B_j \equiv D_{A_j} \gg\ON$ yields 
$\gg(A)\ON \equiv \lim_{j\to \infty} \gg(A_j)\ON = \lim_{j\to \infty} {1\over N}\bra{B_j}{A}$ which proves (3).\qed

%%%%%%%%%%%%%%%%%%%%%%%%%%%%%%%%%%%%%%%%%%%%%
%%%%%%%%%%%%%%%%%%%%%%%%%%%%%%%%%%%%%%%%%%%%%
%%%%%%%%%%%%%%%%%%%%%%%%%%%%%%%%%%%%%%%%%%%%%
%%%%%%%%%%%%%%%%%%%%%%%%%%%%%%%%%%%%%%%%%%%%%
%%%%%%%%%%%%%%%%%%%%%%%%%%%%%%%%%%%%%%%%%%%%%
%%%%%%%%%%%%%%%%%%%%%%%%%%%%%%%%%%%%%%%%%%%%%
%%%%%%%%%%%%%%%%%%%%%%%%%%%%%%%%%%%%%%%%%%%%%
%%%%%%%%%%%%%%%%%%%%%%%%%%%%%%%%%%%%%%%%%%%%%
%%%%%%%%%%%%%%%%%%%%%%%%%%%%%%%%%%%%%%%%%%%%%
%%%%%%%%%%%%%%%%%%%%%%%%%%%%%%%%%%%%%%%%%%%%%
%%%%%%%%%%%%%%%%%%%%%%%%%%%%%%%%%%%%%%%%%%%%%
%%%%%%%%%%%%%%%%%%%%%%%%%%%%%%%%%%%%%%%%%%%%%

\font\tpt=cmr10 at 12 pt
\font\fpt=cmr10 at 14 pt

\font \fr = eufm10

%\font\AAA=Times.dfont  at 12pt
 %\font\BBB=Times.dfont at 8pt

%\font\AAA=cmr10 at 12pt
%\font\BBB=cmr10 at 8pt

%\def\AAA{\bf}
%\def\BBB{\bf}

\overfullrule=0in

\def\boxit#1{\hbox{\vrule
 \vtop{%
  \vbox{\hrule\kern 2pt %
     \hbox{\kern 2pt #1\kern 2pt}}%
   \kern 2pt \hrule }%
  \vrule}}

  \def\harr#1#2{\ \smash{\mathop{\hbox to .3in{\rightarrowfill}}\limits^{\scriptstyle#1}_{\scriptstyle#2}}\ }

\def\AA{1}
\def\BB{2}
\def\CC{3}
\def\DD{4}
\def\EE{5}
\def\FF{6}
\def\GGG{7}
\def\HH{8}
\def\II{9}
\def\JJ{10}
\def\KK{11}
\def\LL{12}
\def\MM{13}

\def\ALL{1}
\def\BTA{2}
\def\BL{3}
\def\BRE{4}
\def\CNS{5}
\def\CIL{6}
\def\CRA{7}
\def\DDD{8}
\def\DDR{9}
\def\GEO{10}
\def\HYP{11}
\def\BEL{12}
\def\AC{13}
\def\SURVEY{14}
\def\NOTES{15}
\def\AET{16}
\def\LAG{17}
\def\KRY{18}
\def\PLI{19}
\def\RT{20}
\def\SLO{21}
\def\TRU{22}
\def\TWC{23}
\def\WAL{24}

 \def\GG{{{\bf G} \!\!\!\! {\rm l}}\ }

\def\GL{{\rm GL}}

\def\bll{I \!\! L}

\def\IFF{\qquad\iff\qquad}
\def\bra#1#2{\langle #1, #2\rangle}
\def\bbf{{\bf F}}
\def\bbj{{\bf J}}
\def\Jtn{{\bbj}^2_n}  \def\JtN{{\bbj}^2_N}  \def\JoN{{\bbj}^1_N}
\def\jt{j^2}
\def\jtx{\jt_x}
\def\Jt{J^2}
\def\Jtx{\Jt_x}
\def\bpp{{\bf P}^+}
\def\bpt{{\wt{\bf P}}}
\def\fsh{$F$-subharmonic }
\def\mo{monotonicity }
\def\jet{(r,p,A)}
\def\ss{\subset}
\def\sse{\subseteq}
\def\half{\hbox{${1\over 2}$}}
\def\smfrac#1#2{\hbox{${#1\over #2}$}}
\def\oa#1{\overrightarrow #1}
\def\dim{{\rm dim}}
\def\dist{{\rm dist}}
\def\codim{{\rm codim}}
\def\deg{{\rm deg}}
\def\rank{{\rm rank}}
\def\log{{\rm log}}
\def\Hess{{\rm Hess}}
\def\Hessyp{{\rm Hess}_{\rm SYP}}
\def\trace{{\rm trace}}
\def\tr{{\rm tr}}
\def\max{{\rm max}}
\def\min{{\rm min}}
\def\span{{\rm span\,}}
\def\Hom{{\rm Hom\,}}
\def\det{{\rm det}}
\def\End{{\rm End}}
\def\Sym{{\rm Sym}^2}
\def\diag{{\rm diag}}
\def\pt{{\rm pt}}
\def\Spec{{\rm Spec}}
\def\pr{{\rm pr}}
\def\Id{{\rm Id}}
\def\Grass{{\rm Grass}}
\def\Herm#1{{\rm Herm}_{#1}(V)}
\def\arr{\longrightarrow}
\def\supp{{\rm supp}}
\def\Link{{\rm Link}}
\def\Wind{{\rm Wind}}
\def\Div{{\rm Div}}
\def\vol{{\rm vol}}
\def\foral{\qquad {\rm for\ all\ \ }}
\def\fpsh{{\cal PSH}(X,\f)}
\def\Core{{\rm Core}}
\def\dis{f_M}
\def\Re{{\rm Re}}
\def\rn{\bbr^n}
\def\pp{\cp^+}
\def\plp{\cp_+}
\def\Int{{\rm Int}}
\def\cix{C^{\infty}(X)}
\def\Gr#1{G(#1,\rn)}
\def\Symn{{\Sym(\rn)}}
\def\SymN{{\Sym(\bbr^N)}}
\def\Gpn{G(p,\rn)}
\def\fd{{\rm free-dim}}
\def\SA{{\rm SA}}
 \def\cd{{\cal C}}
 \def\cdt{{\widetilde \cd}}
 \def\cm{{\cal M}}
 \def\cmt{{\widetilde \cm}}

\def\Theorem#1{\medskip\noindent {\bf THEOREM \bf #1.}}
\def\Prop#1{\medskip\noindent {\bf Proposition #1.}}
\def\Cor#1{\medskip\noindent {\bf Corollary #1.}}
\def\Lemma#1{\medskip\noindent {\bf Lemma #1.}}
\def\Remark#1{\medskip\noindent {\bf Remark #1.}}
\def\Note#1{\medskip\noindent {\bf Note #1.}}
\def\Def#1{\medskip\noindent {\bf Definition #1.}}
\def\Claim#1{\medskip\noindent {\bf Claim #1.}}
\def\Conj#1{\medskip\noindent {\bf Conjecture \bf    #1.}}
\def\Ex#1{\medskip\noindent {\bf Example \bf    #1.}}
\def\Qu#1{\medskip\noindent {\bf Question \bf    #1.}}
\def\Exercise#1{\medskip\noindent {\bf Exercise \bf    #1.}}

\def\HoQu#1{ {\AAA T\BBB HE\ \AAA H\BBB ODGE\ \AAA Q\BBB UESTION \bf    #1.}}

\def\pf{\medskip\noindent {\bf Proof.}\ }
\def\qed{\hfill  $\vrule width5pt height5pt depth0pt$}
\def\equdef{\buildrel {\rm def} \over  =}
\def\qedqed{\hfill  $\vrule width5pt height5pt depth0pt$ $\vrule width5pt height5pt depth0pt$}
\def\mathqed{  \vrule width5pt height5pt depth0pt}

\def\V{W}

\def\df{d^{\phi}}
\def\hk{\_{\rm l}\,}
\def\n{\nabla}
\def\w{\wedge}

\def\cu{{\cal U}}   \def\cc{{\cal C}}   \def\cb{{\cal B}}  \def\cz{{\cal Z}}
\def\cv{{\cal V}}   \def\cp{{\cal P}}   \def\ca{{\cal A}}
\def\cw{{\cal W}}   \def\co{{\cal O}}
\def\ce{{\cal E}}   \def\ck{{\cal K}}
\def\ch{{\cal H}}   \def\cm{{\cal M}}
\def\cs{{\cal S}}   \def\cn{{\cal N}}
\def\cd{{\cal D}}
\def\cl{{\cal L}}
\def\cp{{\cal P}}
\def\cf{{\cal F}}
\def\ccr{{\cal  R}}

\def\gerG{{\fr{\hbox{g}}}}
\def\gerB{{\fr{\hbox{B}}}}
\def\gerR{{\fr{\hbox{R}}}}
\def\p#1{{\bf P}^{#1}}
\def\vf{\varphi}

\def\wt{\widetilde}
\def\wh{\widehat}

\def\and{\qquad {\rm and} \qquad}
\def\arr{\longrightarrow}
\def\ol{\overline}
\def\bbr{{\mathbb R}}\def\bbh{{\mathbb H}}\def\bbo{{\mathbb O}}
\def\bbc{{\mathbb C}}
\def\bbq{{\mathbb Q}}
\def\bbz{{\mathbb Z}}
\def\bbp{{\mathbb P}}
\def\bbd{{\mathbb D}}

\def\a{\alpha}
\def\b{\beta}
\def\d{\delta}
\def\e{\epsilon}
\def\f{\phi}
\def\g{\gamma}
\def\k{\kappa}
\def\l{\lambda}
\def\o{\omega}

\def\s{\sigma}
\def\x{\xi}
\def\z{\zeta}

\def\D{\Delta}
\def\L{\Lambda}
\def\G{\Gamma}
\def\O{\Omega}

\def\bd{\partial}
\def\bdf{\partial_{\f}}
\def\lag{Lagrangian}
\def\psh{plurisubharmonic }
\def\ph{pluriharmonic }
\def\pph{partially pluriharmonic }
\def\omp{$\omega$-plurisubharmonic \ }
\def\ffl{$\f$-flat}
\def\PH#1{\widehat {#1}}
\def\lloc{L^1_{\rm loc}}
\def\dbar{\ol{\partial}}
\def\lp{\Lambda_+(\f)}
\def\lpp{\Lambda^+(\f)}
\def\bo{\partial \Omega}
\def\Ob{\overline{\O}}
\def\fc{$\phi$-convex }
\def\PSH{{ \rm PSH}}
\def\SH{{\rm SH}}
\def\totr{ $\phi$-free }
\def\BM{\lambda}
\def\Der{D}
\def\CH{{\cal H}}
\def\RH{\overline{\ch}^\f }
\def\pconv{$p$-convex}
\def\MA{MA}
\def\lagpsh{Lagrangian plurisubharmonic}
\def\hermsk{{\rm Herm}_{\rm skew}}
\def\PSHl{\PSH_{\rm Lag}}
 \def\ppsh{$\pp$-plurisubharmonic}
\def\fp{$\pp$-plurisubharmonic }
\def\fh{$\pp$-pluriharmonic }
\def\Symn{\Sym(\rn)}
 \def\ci{C^{\infty}}
\def\USC{{\rm USC}}
\def\LSC{{\rm LSC}}
\def\fa{{\rm\ \  for\ all\ }}
\def\ppc{$\pp$-convex}
\def\cpt{\wt{\cp}}
\def\ft{\wt F}
\def\ob{\overline{\O}}
\def\Be{B_\e}
\def\K{{\rm K}}

\def\M{{\bf M}}
\def\N#1{C_{#1}}
\def\ds{Dirichlet set }
\def\dir{Dirichlet }
\def\Fa{{\oa F}}
\def\TR{{\cal T}}
 \def\ISO{{\rm ISO_p}}
 \def\Span{{\rm Span}}

\def\ALL{1}
\def\AV{2}
\def\BTA{3}
\def\BL{4}
\def\BRE{5}
\def\CNS{6}
\def\CP{7}
\def\CPW{8}
\def\CIL{9}
\def\CRA{10}
\def\DTT{11}
\def\DON{12}
\def\CG{13}
\def\DDD{14}
\def\DDR{15}
\def\GEO{16}
\def\HYP{17}
\def\BEL{18}
\def\SURVEY{19}
\def\AC{20}
\def\NOTES{21}
\def\AET{22}
\def\TANG{23}
\def\TANGG{24}
\def\LAG{25}
\def\SLE{26}
\def\JTY{27}
\def\KRY{28}
\def\PLI{29}
\def\RT{30}
\def\SLO{31}
\def\SPR{32}
\def\TRUU{33}
\def\TRU{34}
\def\TWC{35}
\def\TWCC{36}
\def\TWCCC{37}
\def\WAL{38}

\def\II{1}
\def\AA{2}
\def\AB{3}
\def\EE{4}
\def\BB{5}
\def\CC{6}
\def\CCC{7}
\def\DD{8}

\vskip .4in

\def\ip{\bra}
\def\Gc{\Gamma}
\def\E{E}
\def\fpsi{{F_f(\psi)}}
\def\bL{{\bf \Lambda}}
\def\bdf{{\bf f}}
\def\UU{U}
\def\bbm{{\bf M}}
\def\gg{{\mathfrak g}}
\def\ggL{\gg_L}
\def\GL{\G_L}
\def\gL{\ggL}
\def\Imm{{\rm Im}}
\def\EL{E_L}
\def\eqdef{\overset {\rm def} =}
\def\Gcs{\Gc^*}
\def\GLs{\GL^*}
\def\Sig{\Sigma}
\def\Dlog{D\, \log\,}
\def\Lt{L^t}
\def\gLs{\gL^*}
\def\gs{\gg^*}
\def\qed{\hfill $\square$}
\def\Ga{G\aa rding\ }
\def\ed{\end{document}}
\def\IV{^{-1}}
\def\nm#1{\eqno{(#1)}}

\font\headfont=cmr10 at 14 pt
\font\aufont=cmr10 at 11 pt

\bigskip

\centerline{\bf \headfont 13. Duality Under Pull-Back of G\aa rding Polynomials}
%\label{sec:13}

In this section we describe the general case. Some of the many important special cases will
be covered in Section 14.

Suppose $g,\Gc$ is a degree $N$ G\aa rding polynomial with G\aa rding cone $\Gc$
on an $n$-dimensional inner product space $W,\bra \cdot\cdot$. The
{\bf pull back} $\gg,\GL$ of $g,\Gc$ to another inner product space
$X,\ip\cdot,\cdot$, with $X$ of dimension $m$, by a linear map
$$
L : X \longrightarrow W
$$
is defined as follows.
$$
\gL(x) \equiv g(Lx) \quad \forall x \in X
\qquad\text{and}\qquad
\GL \equiv L^{-1}(\Gc).
\eqno{(13.1)}
$$
We will always assume (the \emph{Pull Back} hypothesis)
$$
\GL \neq \emptyset,
\qquad\text{or equivalently}\qquad
(\Imm L)\cap\Gc \neq \emptyset.
\eqno{(13.2)}
$$

The objective of this section (Lemma 13.1 and Theorem 13.2)  is to show that the various quantities associated
with $\gL,\GL$ can be computed in terms of the similar quantities for $\gg,\Gc$.
We start by computing $\gL$-eigenvalues of $x\in X$ in terms of $\gg$-eigenvalues
of $Lx$.

\vfill\eject

\Lemma{13.1. (The Pull Back)} {\sl

\noindent
  (a) {\rm (G\aa rding Polynomial.)} $\gL$ is a degree $N$ G\aa rding
        polynomial on $X$ with G\aa rding cone $\GL$.
        
        \noindent
 (b) {\rm (Eigenvalues)} For each $e\in\GL$,
       $$
          \lambda_{j}^{\gL,e}(x) \;=\; \lambda_{j}^{g,Le}(Lx)
          \qquad \forall\, j=1,\dots,N \text{ and } \forall\, x\in X.
          \eqno{(13.3)}
             $$
        
        \noindent
(c) {\rm (Edge)} $\EL = L^{-1}(\E)$.

}

\noindent
{\sl Proof  (b) $\Rightarrow$ (a)}.
In particular, b) says that the $\gL,e$-eigenvalues of $x\in X$ are real since the $\gg, Le$ eigenvalues of $Lx$ are real; and
that they are $>0$ iff the $g,Le$-eigenvalues of $Lx$ are $>0$. Thus $\gL$ has
G\aa rding cone
\[
\{ x\in X : \lambda_{j}^{\gL,e}(x) > 0,\ j=1,\dots,N \}
\;=\;\{ x\in X : Lx \in \Gc \}.   \qquad \square
\]

\noindent
{\sl Proof of (b)}
\[
\gL(te+x) \equiv g(tLe+Lx) \;=\; g(Le)\prod_{j=1}^{N}\bigl(t + \lambda_{j}^{g,Le}(Lx)\bigr),
\]
as well as
\[
\gL(te+x) \;=\; \gL(e)\prod_{j=1}^{N}\bigl(t + \lambda_{j}^{\gL,e}(x)\bigr). \qquad\square
\]

\noindent
{\sl Proof of (c).}
Since $\overline{\Gamma_L} = L^{-1}(\overline\Gamma)$, we have $E_L = L^{-1}(E)$. \qed

Now we turn to the dual pair $(g_L)^{*},\Gamma_L^{*}$. The computations can be organized around  the useful diagram of maps.

$$
\begin{array}{ccccccc}
\Gamma_L/E_L & & \Gamma_L \subset X & \xrightarrow{\ L\ } & W \supset \Gamma & & \Gamma/E \\[4pt]
\Big\downarrow{\scriptstyle \cong} & & & & \Big\downarrow{\scriptstyle F \,\equiv\, D\log g} & & \Big\downarrow{\scriptstyle \cong} \\[4pt]
\Gamma_L^{*} & & \Sigma \subset \Gamma^{*} & \xleftarrow{\ L^t\ } & \Gamma^{*} & & \Gamma^{*} \\[4pt]
& & {\scriptstyle \cap} & & {\scriptstyle \cap} & & \\[4pt]
& & X & & W & &
\end{array}
\eqno{(13.4)}
$$

\noindent
where the left vertical map is $F_L \equiv D\log g_L$.
Now this left vertical map can be rewritten as the composition of $L^t\circ F\circ L$ in the middle diagram, divided by the edge $E_L$.  The vertical map on the right is $F$ pushed down to the quotient of $\Gamma$ by the edge.  In the following, 
we will be  studying these maps and the isomorphisms.

By the gradient diffeomorphism Theorem 7.4,
$$
F \eqdef D\,\log \, g : \Gc/\E \longrightarrow \Gcs
\eqno{(3.5)}
$$
is a diffeomorphism. By Lemma 13.1,
$$
F_{L} \eqdef D\, \log\, \ggL : \GL/\EL \longrightarrow \GLs
\eqno{(13.6)}
$$
is also a diffeomorphism because of Thm.~7.4.

We define $\Sig$ to be
$$
\Sig \eqdef F\bigl(\Gc \cap L(X)\bigr) \;=\; \{ F(Lx) : x\in \GL \},
\eqno{(13.7)}
$$
the Image of the slice of $\Gc$ by the subspace $L(X)$ ($\equiv \Imm L$) under the
gradient diffeomorphism $F\equiv \Dlog g$.

\medskip
\noindent\textbf{Note.} Our pull back hypothesis (13.2) that $\Gc\cap L(X)\neq\emptyset$
is true iff $\Sig$ is a non-empty closed submanifold of $\Gcs$ of dimension equal to  the
dimension of $L(X)$, or codimension equal to the dimension of the $\ker L^t$.

\medskip
Now we can state the general duality formula under pull backs.

\Theorem {13.2. (The Pull Back Dual)} {\sl
The dual pair $\gLs,\GLs$ of the pull back pair $\gL,\GL$ is determined by the
dual pair $\gs,\Gcs$ of $\gg,\Gc$ as follows.

\noindent
{\bf (1).  (Polar Cones)}
      $$
          (\GL)^{*} = \Lt\Gcs,
          \qquad\text{or equivalently}\qquad
          (\GL)^{0} = \Lt\Gc^{0}.
         \eqno{(13.8)}
       $$
 
\noindent {\bf  (2).  (Gradient Maps)}

 (a)
        $$
          D_{x}\gL \;=\; \Lt D_{Lx} \gg \qquad \forall\, x\in\GL.
              \eqno{(13.9)}
                    $$
                    
   (b)
        $$
          F_{L}(x) \eqdef D_{x}\log \gL \;=\; \Lt D_{Lx}\log \gg
          \;\equiv\; \bigl(\Lt\circ F\circ L\bigr)(x)
              \eqno{(13.10)}
                     $$
                     
                     \noindent for all $x\in\GL$. This factors the gradient diffeomorphism
        $$
          F_{L}\equiv \Dlog \gL : \GL/\EL \longrightarrow \GLs
           \eqno{(13.11)}
       $$
        into a diffeomorphism

     (c)  $$ F\circ L : \GL/\EL \longrightarrow \Sig
        \eqdef \{F(Lx) : x\in\GL\}, 
        \eqno{(13.12)}
        $$

(d)    followed by a diffeomorphism 
$$
\Lt : \Sig \longrightarrow \GLs.
\eqno{(13.13)}
$$

\noindent {\bf (3).   (Dual Function)}
        $$
          \gLs(y) \;=\; \gs\!\bigl((\Lt)^{-1}(y)\bigr) \qquad \forall\, y\in\GLs,
        \eqno{(13.16)}
     $$
        where $(\Lt)^{-1} : \GLs \to \Sig$ denotes the inverse of the
        diffeomorphism in (d).
        
        \smallskip
\noindent {\bf (4).  (Ellipticity)}
  $\gL^*,\GL^*$ is elliptic if and only if  
  
  \noindent
 { \bf Case 1. $(X= \rn)$:}  $L(\rn_{>0}) \supset \G$, or equivalently, $\rn_{>0} \ss L^t\G^*$,

  \noindent
 { \bf Case 2. $(X= S(n))$:}  $L(\Int \cp) \supset \G$, or equivalently, $\Int \cp\ss L^t\G^*$.

}

\Remark {13.3} 
If the inverse $F^{-1}$ of $F\equiv \Dlog g : \Gc \to \Gcs$ is known explicitly,
then the inverse $F_{L}^{-1}$ of $F_{L}\equiv \Dlog \gL : \GL \to \GLs$ can be
computed in terms of $(L^t)^{-1}$ followed by $(f\circ L)^{-1}$ because of  (2).

\medskip

The meaning of part 4 above will be clearly explained in Remark 13.4 below where ellipticity is discussed 
for {\bf any G\aa rding - like pair} which includes    $\gL^*,\GL^*$.

\noindent
{\sl Proof of (1)}.
The fact that the polar $K^{0}$ of the preimage $K\equiv L^{-1}(C)$ of a closed
convex cone $C$ equals $\Lt(C^{0})$ for a linear map $L$ is standard in convex
analysis (see, e.g., [Ro, Theorem 14.1]). Hence
$(\GL)^{0} = \Lt(\Gc^{0})$. The argument goes briefly as follows. The proof
that $\Lt(C^{0}) \subset (L^{-1}(C))^{0}\equiv K^{0}$ is straightforward from
the definitions. The proof that $(L^{-1}(C))^{0} \subset \Lt C^{0}$ is a little
harder. It requires the use of the Hahn--Banach Theorem to extend the linear
functional on $\Imm L$ given by $\varphi_{y}(L(x))\equiv\bra xy$ for
$y\in(L^{-1}(C))^{0}$, after showing  that  
$\vf_y$ is well-defined on
${\rm Im} L$.  \qed

\noindent
{\sl Proof of (2a) and (2b).}
Note that for all $x,y\in X$,
$$
\begin{aligned}
\ip{D_{x}\gL,y} \; &=\; \left.\frac{d}{dt}\gL(x+ty)\right|_{t=0} 
\; =\; \left.\frac{d}{dt} \gg\bigl(L(x)+tL(y)\bigr)\right|_{t=0} \\
\;&=\; \ip{D_{Lx}\gg,\,Ly} \;=\; \ip{\Lt D_{Lx}\gg,\,y},
\end{aligned}
$$
which proves that $D_{x}\, \ggL = \Lt D_{Lx}\,\gg$. Dividing both sides by
$\gL(x) = \gg(L(x))$ proves that $D_{x}\log\,  \gL = \Lt D_{Lx}\log \,\gg$.

%%%%%%%%%%%%%%%%%%%%%%%%%%%%%%%%%%%%%%%%%%%%
%%%%%%%%%%%%%%%%%%%%%%%%%%%%%%%%%%%%%%%%%%%%
%%%%%%%%%%%%%%%%%%%%%%%%%%%%%%%%%%%%%%%%%%%%
%%%%%%%%%%%%%%%%%%%%%%%%%%%%%%%%%%%%%%%%%%%%
%%%%%%%%%%%%%%%%%%%%%%%%%%%%%%%%%%%%%%%%%%%%

\noindent
{\sl  Proof of (2c).}
We have $F\circ L : \Gamma_L \to \Sigma$ since $x\in \Gamma_L$ implies $Lx\in  \Gamma$ and $F(Lx)\in \Sigma$.
By definition of $\Sigma$ the map $F\circ L : \Gamma_L \to \Sigma$ is onto.
If $x\in E_L$, then $Lx\in E$, and hence by the gradient diffeomorphism (13.7), we have $F(Lx)=0$, which proves 
that $F\circ L:\Gamma_L/E_L\to \Sigma$ is well defined and onto.
To prove it is injective  it suffices to prove that $L^t\circ F\circ L :  \Gamma_L/E_L \to \Gamma_L^*$ is injective.
Now by (2b) we have $L^t\circ F\circ L = F_L$ is the gradient diffeomorphism  for $\log \gg_L$, so it is injective. \qed
\medskip

\noindent
{\sl  Proof of (2d).}  Since $F_L = L^t\circ (F\circ L)$ and $F\circ L$ are diffeomorphisms, 
$L^t= F_L\circ (F\circ L)^{-1} : \Sigma \to {\Gamma_L}^*$  
 is a diffeomorphism with inverse
$$
(L^t)^{-1} \ =\ F\circ L\circ F_L^{-1} : \Gamma_L^* \to \Sigma.  
\eqno{(13.15)}
$$
\qed

\medskip
\noindent
{\sl  Proof of (3).}
Suppose that $w\in\Sig$, i.e., $w = F(Lx)\equiv D_{Lx}\log g$ for some $x\in\GL$, and
that $\Lt w = y$. Then
\[
y \;=\; \Lt D_{Lx}\log g \;=\; D_{x}\log \gL \qquad\text{by (2)}.
\]

Therefore by the identities $\gg(w)\gg^*(F(w))=1$ and $\gg_L(x)\gg_L^*(F_L(x))=1$, we have

\[
\gLs(y) \;=\; \frac{1}{\gL(x)} \;=\; \frac{1}{g(Lx)}
\;=\; \gs(D_{Lx}\log g) \;=\; \gs(w). 
\]
\qed

\Remark{13.4. (Explaining  \lq\lq Ellipticity")}
Our use of the term {\sl elliptic} for a \Ga pair depends  on the context.
There are two cases.

{\bf Case 1. (The Eigenvalue or diagonal case).}
Here $\gg,\G$ is a \Ga pair on $X\equiv \rn$ and the actual elliptic  operator on $S(n)$ in obtained by either (assuming $\gg,\G$ is invariant) applying $\gg$ to the eigenvalues of $A\in S(n)$, or by using projection $L:S(n) \to \bbd$ onto the diagonal 
$\bbd \equiv \rn$ to pull back $\gg,\G$ to $S(n)$.  The pair $\gg,\G$ is {\bf elliptic} if
$$
\Gamma\supset {\mathbb R}^n_{>0} \qquad\text{or equivalently}\qquad  {\mathbb R}^n_{>0}\supset \G^*
\eqno{(13.16)}
$$
(recall that ${\mathbb R}^n_{>0}$ is self-polar).

%$X\equiv {\mathbb R}^n$ and $\gg, \Gamma$ is {\sl elliptic} if $\Gamma\supset {\mathbb R}^n_{>0}$ 

{\bf Case 2. (G\aa rding polynomial on $S(n)$).} Here $X\equiv S(n)$ and 
$\gg,\Gamma$ is {\sl elliptic} if
$$\
\Gamma \supset {\rm Int} {\mathcal P}  \qquad\text{or equivalently}\qquad  {\rm Int} {\mathcal P} \supset \G^*
\eqno{(13.17)}
$$
(siince $\cp$ is self-polar).
As a consequence of our definitions, $\gg,\G$ is elliptic if and only if $\gg,\G$ is a G\aa rding-Dirichlet operator.

The same discussion applies  to the dual pair $\gg^*,\G^*$ of a \Ga pair $\gg,\G$ so that (13.16) and (13.17) apply with
$\G$ and $\G^*$ interchanged.

\medskip
\noindent
{\sl Proof of (4).}  With these definitions, along with part (a) that $\G_L^* = L^t\G^*$, the statement (4) follows easily. \qed

%%%%%%%%%%%%%%%%%%%%%%%%%%%%%%%%%%%%%%%%%%%%%%%%%%%%%%%%%%%%%%%
%%%%%%%%%%%%%%%%%%%%%%%%%%%%%%%%%%%%%%%%%%%%%%%%%%%%%%%%%%%%%%%
%%%%%%%%%%%%%%%%%%%%%%%%%%%%%%%%%%%%%%%%%%%%%%%%%%%%%%%%%%%%%%%
%%%%%%%%%%%%%%%%%%%%%%%%%%%%%%%%%%%%%%%%%%%%%%%%%%%%%%%%%%%%%%%
%%%%%%%%%%%%%%%%%%%%%%%%%%%%%%%%%%%%%%%%%%%%%%%%%%%%%%%%%%%%%%%
%%%%%%%%%%%%%%%%%%%%%%%%%%%%%%%%%%%%%%%%%%%%%%%%%%%%%%%%%%%%%%%
%%%%%%%%%%%%%%%%%%%%%%%%%%%%%%%%%%%%%%%%%%%%%%%%%%%%%%%%%%%%%%%
%%%%%%%%%%%%%%%%%%%%%%%%%%%%%%%%%%%%%%%%%%%%%%%%%%%%%%%%%%%%%%%
%%%%%%%%%%%%%%%%%%%%%%%%%%%%%%%%%%%%%%%%%%%%%%%%%%%%%%%%%%%%%%%

\vskip .5in

\vfill\eject

\centerline{\bf \headfont 14. Important Special Cases of the Pull Back Theorem 13.2}
%\label{sec:13}
\medskip

\noindent
{\bf 14.1.   Linear Coordinate Changes - $ L \in GL_n(X)$.}  Suppose $L : X \arr X$ is a linear isomorphism, and $\gg,\G$  is a  \Ga pair on $X$. The hypothesis (13.2) states that $\G_L = L^{-1}(\G)$ is non-empty; this is automatic.

Since  $L(X) = X$, the subset $\Sigma = F(\G\cap L(X))$ defined in (13.7) equals the image of $\G$ under the gradient map 
$D\,\log\, \gg$, which is just $\Sigma = \G^*$ by the Gradient Diffeomorphism Theorem.

By Parts (1) (Polar Cones) and (2) (Gradient Maps) of Theorem 13.2, the same simplification remains. Parts (3) and (4) simplify as follows.

\Cor{14.1. $L \in  GL(X)$}  {\sl

(3a)
 $$
 L^t : \G^* \arr \G^*_L  \text{ is a linear isomorphism of cones.}
\eqno{(14.1)}
$$

}

 (3b)
 $$
 \gg_L^* = \gg^*((L^t)^{-1}(y))  \qquad\text{\sl for all $y\in \G_L^*$}
\eqno{(14.2)}
 $$
 
 \noindent
{\sl  where $(L^t)^{-1}$ denotes the inverse in (3a).
}

\noindent
In particular, if $\gg^*$ is a polynomial, then the dual $\gg_L^*$ is also a polynomial.

\noindent
{\bf 14.1 A (Matrices).} Let $L : \rn \arr \rn$  be a linear transformation given by a matrix with rows $a^1, ... ,a^n$  in $\rn$ and so $\gg_L(x) = \gg(a^1\cdot x, ... , a^n\cdot x)$.
The dual basis $b_1,..., b_n$ in $\rn$ (i.e.,  $a^i \cdot  b_j = \delta_{ij}$) is the set of column vectors for the matrix representing 
$(L^t)^{-1} : \rn \arr\rn$, so that:
$$
\gg_L^*(y)=\gg^*(b_1 \cdot y,...,b_n \cdot y),   \quad{\rm with} \quad  \g_L^* = \{y:(b_1 \cdot y, ...  \, ,b_n \cdot y)\in \G^*\}.
\eqno{(14.3)}
$$
This is a useful formula for the dual function $\gg_L^*$. 

For example, if  $\gg = \det$, then:
$$
 \gg_L^*(y) =  \prod_{i=1}^n (b_i  \cdot y) \quad {\rm for\ all } \quad y \in \G_L^* = \{y : b_1 \cdot y, ... \, ,b_n \cdot y > 0\}.
\eqno{(14.4)}
$$

\noindent
{\bf 14.1 B (Orthogonal Coordinate Changes). } If $ L \in {\rm  O}(n)$, then $L\IV = L^t$, so $\G_L = L^t\G$  and
$\gg_L(x) = \gg(Lx)$ for all $x \in  \G_L$. Therefore:
$$\gg_L^* (y) = \gg^*(Ly) \quad {\rm for \ all}\quad  y \in \G_L^* = L^t\G^* = L\IV\G^*.
\eqno{(14.5)}
$$
(The rows $a^1, . . . , a^n$ form an orthonormal basis.)

\noindent
{\bf 14.1 C (Symmetric Coordinate Changes). }  If $L \in  \Sym(\rn)$, i.e.,$ L^t = L$, then $\G_L  = L\G$
and $ \gg_L(x) = \gg(Lx)$ for all $x \in \G_L$. Moreover, 
$$
\gg_L^* (y) = \gg^*(L\IV y) \quad  {\rm for\  all}\quad  y \in \G_L^* = L\IV\G^*.
\eqno{(14.6)}
$$

\noindent
{\bf 14.1 D (Double Rescalings).}  A simple but very important special case is when the symmetric invertible linear transformation 
$L \in {\rm GL}(W)$ has exactly two eigenvalues, with eigenspaces $W = E\oplus E^\perp$. By rescaling $L$, 
we may assume that $L$ belongs to the family of symmetric transformations:
$$
 L_t =P_{E^\perp} +tP_E = I+(t - 1)P_E, \qquad  t\in \bbr^*   = R- \{0\}.
\eqno{(14.7)}
$$

\noindent
There is a simple composition law:
$$
 L_t \circ  L_s = L_{ts},\quad  t, s \in  \bbr^*.
\eqno{(14.8)}
$$
Since $ P_{E^\perp} \circ  P_E = P_E \circ  P_{E^\perp} = 0$, $L_t$ is invertible with inverse:
$$
L_t\IV = L_{1\over t}, \quad t\in \bbr^*.
\eqno{(14.9)}
$$

\noindent
{\bf 14.2. Generalized Elliptic Regularization.} The  pullback Theorem 13.2 applied to this last case 14.1D can be stated in a simpler form as follows.

\Theorem {14.2} {\sl Given a \Ga pair $\gg, \G$ on $W$ , and a subspace $ E \ss W$ , 
define $ L_{t,x } \equiv  P_{E^\perp} x + tP_Ex$ for all $x\in W$, $t\neq 0$. Then
$$
\G_t \equiv \{x\in W : L_{t,x} \equiv x+(t-1)P_Ex\in \G\},
\nm{(14.10)} 
$$
and
$$
\gg_t(x) \equiv \gg ( x + (t -1)P_E x).
\nm{(14.11)}
$$
The dual 
$$
\G^*_t = \left \{ y \in  W : y + \left ({1\over t } - 1 \right) P_E y \in  \G^*  \right \} ,
$$
and
$$
\gg^*_t(y) = \gg^* \left(y +\left( {1\over t} - 1\right)P_Ey \right)\quad {\rm  for\ all}\ \  y \in  \G^*_t .
$$
}

\noindent
{\bf 14.3. Elliptic Regularization. } The case where $\dim  E = 1$ is particularly important. Choose a
unit vector $e$ in $E = \bbr e$. Then $P_Ex =  \bra xe e$ and
$$
L_{t,x}  = P_{E^\perp}x+t \bra xe e=x+(t-1) \bra xe e.
$$
Consequently,

\Lemma{ 14.3} {\sl  $L_0 = P_{E^\perp}$ is orthogonal projection onto the hyperplane $e^\perp \equiv \{x : \bra xe = 0\}$ with kernel $\bbr \cdot e = E$.   Moreover, $ L_t$ fixes this hyperplane $e^\perp = \{x : \bra xe  = 0\}$, while:

$\bullet$ if  $t \geq 1$,  then $L_t$  expands the upper (and lower) halfspace $\{\bra xe > 0\}$ ($\{\bra xe < 0\}$) by $t$,

$\bullet$  and if $0<t\leq 1$,  then $ L_t$ contracts these half spaces by $t$.

\noindent
If  $t = -1$,  then $L_{-1}x = P_{E^\perp}x - P_Ex$ is reflection through the hyperplane $e^\perp$  along the perpendicular line $\bbr e$. Hence for $t \neq  -1$ and negative, $L_t$  flips the half spaces $\{ \bra xe > 0\}$ and $\{ \bra xe < 0\}$  followed by $L_{|t|}$ acting as above.
}

%%%%%%%%%%%%%%%%%%%%%%%%%%%%%%%%%%%%%%%%%%%%%%%%%%%%%%%%%%%%%
%%%%%%%%%%%%%%%%%%%%%%%%%%%%%%%%%%%%%%%%%%%%%%%%%%%%%%%%%%%%%
%%%%%%%%%%%%%%%%%%%%%%%%%%%%%%%%%%%%%%%%%%%%%%%%%%%%%%%%%%%%%
%%%%%%%%%%%%%%%%%%%%%%%%%%%%%%%%%%%%%%%%%%%%%%%%%%%%%%%%%%%%%
%%%%%%%%%%%%%%%%%%%%%%%%%%%%%%%%%%%%%%%%%%%%%%%%%%%%%%%%%%%%%
\def\Sn{S(n)}
\def\Ltt{L_t}
\def\Fcal{{\mathcal F}}
\def\Fcals{\Fcal^*}

\vskip .5in

\centerline{\bf \headfont 15. Elliptic Regularization and Duality}

\medskip

All of the previous results apply to the case $W \equiv \Sn$, which we are
primarily interested in. The identity $I \in \Sn$ corresponds to singling out
a vector $e \in W$, so that the previous family $\Ltt$ becomes the {\bf elliptic regularization family} of Krylov [K]:
$$
\begin{aligned}
\Ltt A \;&=\; A + (t-1)(\tr A)\,\tfrac{1}{n}\, I \; \\
&=\; P_{\{\tr = 0\}} A + t\,(\tr A)\,\tfrac{1}{n}\, I,
\qquad \forall\, A \in \Sn,\ t \neq 0.
\end{aligned}
\nm {15.1}
$$
We wish to examine the action of $\Ltt$ on a closed convex cone
$\Fcal \subset \Sn$ with vertex at the origin, under the following two
conditions.

\

\noindent
{\bf Assumption 1.}
$\Fcal \subset \{ A \in \Sn : \tr A \geq 0 \}$ (abbreviated $\{\tr \geq 0\}$).

\noindent Recall the edge of $\Fcal$ is defined to be
$$
E(\Fcal) \;=\; \Fcal \cap (-\Fcal).
\nm{15.2}
$$

\noindent
{\bf Assumption 2.}
$E(\Fcal) \;=\; \Fcal \cap \{\tr = 0\}$; i.e., $E(\Fcal) \subset \{\tr = 0\}$.

\noindent That is, $\Fcal$ is contained in the $\{\tr \geq 0\}$ half-space
with its edge $E(\Fcal)$ contained in the $\{\tr = 0\}$ hyperplane.

%\medskip
%\noindent[The following lemma \emph{was} Lemma~B.1 in Appendix~B and is retained
%here for reference. It records the elementary action of $\Ltt$ on the
%half-spaces about $\{\tr = 0\}$.]

\Lemma{15.0}  {\sl
$L_0$ is projection onto the hyperplane $\{\tr = 0\}$ with kernel $\R I$. If
$t > 0$, $\Ltt$ fixes $\{\tr = 0\}$, and
\begin{itemize}
  \item if $t \geq 1$ it expands the upper and lower half-spaces by $t$;
  \item if $t \leq 1$ it contracts the upper and lower half-spaces by $t$.
\end{itemize}
If $t < 0$, then $\Ltt(\sigma,\tau) = (\sigma,-|t|\tau)$, so the action is the
same after interchanging the upper and lower half-spaces.
}

\medskip
We now suppose that $\Fcal$ is a closed convex cone in $\{\tr \geq 0\} \subset \Sn$
with vertex at the origin, and which contains $I$. We recall that the edge of
$\Fcal$ is
\[
E(\Fcal) \;\equiv\; \Fcal \cap (-\Fcal) \;\subset\; \{\tr = 0\}.
\]
We shall make the additional assumption that
$$
E(\Fcal) \;=\; \Fcal \cap \{\tr = 0\}.
\nm{15.3}
$$
However, this condition is \textbf{always satisfied} if $\Fcal = \Gc$, the
G\aa rding cone for a G\aa rding polynomial $g$ on $\Sn$, because here the edge
equals the nullity of $g$ (see (5) in Section~2).

% Renumbering: formerly Lemma B.2 -> now Lemma 15.1
\Lemma{15.1} {\sl
Under these two assumptions the following conditions are equivalent.
\begin{enumerate}
  \item[\textup{(a)}] $E(\Fcal) = \{0\}$. \medskip
  
  \item[\textup{(b)}] $B \equiv \{ A \in \Fcal : \tr A = 1 \}$ is a compact set
        (and so  a compact base for $\Fcal$).\medskip
        
  \item[\textup{(c)}] $I \in \Fcals$ ($\equiv$ the open polar of $\Fcal$).
\end{enumerate}
}

\noindent
{\sl Proof.} $(a)\Rightarrow(b)$: Assume (a). Then $\Fcal \cap \{\tr = 1\}$ is closed since
$\Fcal$ is closed. It is bounded since, if not, there exists
$A_j \in \Fcal \cap \{\tr = 1\}$ with $|A_j| \to \infty$. We can assume that
$A_j/|A_j|$ converges to a non-zero element $A \in \{\tr = 0\}$. Then the ray
$\overline{0 A_j}$ converges to the ray $\overline{0 A} \subset \{\tr = 0\}
\subset \Fcal$ (since $\Fcal$ is closed), which with Assumption 2
contradicts (a).

$(b)\Rightarrow(c)$: The statement (b) says that the base $B$ for $\Fcal$ is a
compact set. We have assumed that $I \in \Fcal$, so $\tfrac{1}{n}I \in B$. Thus
$\ip{I, B} > 0$ for all $B \in \Fcal - \{0\}$, which means that (c) holds.

$(c)\Rightarrow(a)$: Note that $I^{\perp} = \{\tr = 0\}$, so if (c) holds, then
a neighborhood of $I$ lies in the polar of $\Fcal$. This means that for any
direction $B$ in $\{\tr = 0\}$, there exist $A \in \Fcals$ with $\bra AB < 0$.
Thus $B = 0$ and (a) holds.

\medskip
It is easy to see the following.

% Renumbering: formerly Lemma B.3 -> now Prop. 15.2
\Prop{15.2} {\sl
Suppose $\Fcal$ satisfies Assumptions 1 and 2,  and the
conditions in Lemma 15.1. Then:
\begin{itemize}
  \item As $t \uparrow \infty$, $\Ltt(\Fcal)$ contracts monotonically down to
        the ray $\{tI : t \geq 0\}$, and $L_{-t}(\Fcal)$ contracts monotonically
        down to the ray $\{tI : t \leq 0\}$.
              \medskip
  \item As $t \downarrow 0$, $\Ltt(\Fcal)$ expands monotonically to
        $\{\tr > 0\} \cup \{0\}$, and $L_{-t}(\Fcal)$ expands monotonically to
        $\{\tr < 0\} \cup \{0\}$.
\end{itemize}
}

%%%%%%%%%%%%%%%%%%%%%%%%%%%%%%%%%%%%%%%%%%%%%%%%%%%%%%%%%%%%%
%%%%%%%%%%%%%%%%%%%%%%%%%%%%%%%%%%%%%%%%%%%%%%%%%%%%%%%%%%%%%
%%%%%%%%%%%%%%%%%%%%%%%%%%%%%%%%%%%%%%%%%%%%%%%%%%%%%%%%%%%%%
%%%%%%%%%%%%%%%%%%%%%%%%%%%%%%%%%%%%%%%%%%%%%%%%%%%%%%%%%%%%%
%%%%%%%%%%%%%%%%%%%%%%%%%%%%%%%%%%%%%%%%%%%%%%%%%%%%%%%%%%%%%

\vskip .5in

\centerline{\bf \headfont 16. Duality for General Cubic}\centerline{\bf\headfont G\aa rding Polynomials in Two Variables}

\medskip

Let $\gg,\Gamma \subset \mathbb{R}^2$ be a G\aa rding polynomial with G\aa rding cone $\Gamma$. 
Note that since $\gg$ vanishes on the boundary $\partial\Gamma$, it must be a product of three linear functions, two of which define $\partial\Gamma$. 
These can be normalized to be the standard coordinates $x_1,x_2$ in $\mathbb{R}^2$, and so $\gg(x) = x_1x_2(\alpha_1x_1 + \alpha_2x_2)$ with $\alpha_k \ge 0$ for $k=1,2$. 
Thus $x_1 > 0$, $x_2 > 0$ defines $\Gamma$, and $x\in\G
\Rightarrow \a_1x_1+\a_2x_2>0$

 If both $\a_k>0$ for $k=1,2$, we can normalize by the linear map 
$$
\bar x_1 = {1\over \a_1}x_1 \and  \bar x_2 = {1\over \a_2}x_2 
$$

.

%If both $\alpha_k > 0$ for $k=1,2$, we can normalize by the linear map
%\\[\\bar x\_1 = \\frac{1}{\\alpha\_1} x\_1   \\qquad\\text{and}\\qquad  \\bar x\_2 = \\frac{1}{\\alpha\_2} x\_2\\]
%which reduces to considering $g(x) = x_1x_2(x_1 + x_2)$.

Hence the general case of a cubic comes down to understanding the cases:

\medskip
\noindent
\textbf{General Case (1).} $\gg(x) = x_1x_2(x_1 + x_2)$.

\noindent
\textbf{Degenerate Case (2).} $\gg(x) = x_1x_2^2$ (or $x_1^3$, which we omit).

\medskip
\noindent
\textbf{Case (2):} The gradient map is
$$F(x) = D_x \log \,\gg = \left(\frac{1}{x_1}, \frac{2}{x_2}\right) = (y_1,y_2)$$
and
$$
\gg^*(y) = \frac{1}{\gg(x)} = \frac{1}{x_1x_2^2} = {1\over 4}y_1 y_2^2
$$
is the dual function.

This reduces us to Case (1) (which in eigenvalue space is the determinant times the trace). 
We begin with the following general fact.

\medskip
\noindent
\textbf{Lemma 16.1.} 
If $p(x,y)$ is G\aa rding hyperbolic in the direction $e_2$ in $\mathbb{R}^2$, then
\begin{equation}
 p(x,y) = a \prod_{j=1}^N (y - r_jx) \tag{16.1}
\end{equation}
for $a>0$, and
$$
\Gamma = \{(x,y) \in \mathbb{R}^2 : y >  r_{\rm max} x \ {\rm and} \ 
y >  r_{\rm min} x \}.
$$

\

\noindent
{\sl Proof.}
So $p(x,y) \in \mathbb{R}[x,y]$ is homogeneous of degree $N \ge 1$ and hyperbolic in the direction $e=(0,1)$ with $p(0,1)=a>0$. 
Then $p(te + (x,y)) = p(x,t + y)$ has $N$ real roots in $t$. 
Take $x=1$, $y=0$, so that $p(1,t) = a\prod_{j=1}^N(t-r_j)$, where the $r_j$'s are real. 
Then
$$p(x,y) = a x^N p(1,y/x) = a x^N \prod_{j=1}^N \left( \frac{y}{x} - r\_j \right)       = a \prod_{j=1}^N (y - r\_jx)$$
is the product of linear factors. Now
$$
p(te+(x,y)) = p(x,t+y) = a\prod_{j=1}^N (t + (y-r\_jx))
$$
so the $e$-eigenvalues of $(x,y)$ are $y-r_jx$, $j=1,\dots,N$. 
Hence,
$$
\Gamma^* = \{(x,y) \in \mathbb{R}^2 : y > r_{\max}x \text{ and } y > r\_{\min}x\}.
$$
Note that if $y > r_{\max}x$ and $y > r_{\min}x$, then $y > r_jx$ for all other $j$, so (16.1) is proved. 
\hfill $\Box$

\

\bigskip
\noindent
{\bf Case $N=3$}
The projective linear group acts simply and transitively on the lines in $\mathbb{R}^2$. 
Hence we can assume that when $N=3$, $g(x) = x_1x_2(x_1 + x_2)$, and $\Gamma = \mathbb{R}^2_{>0}$ if the three lines are distinct. 

Now we return to the generic case.

\medskip
\noindent
\Prop{16.2} {\sl
Set $g(x) \equiv x_1x_2(x_1 + x_2)$ and $\Gamma = \mathbb{R}^2_{>0}$. 
Then:
\begin{enumerate}
 \item[(1)] \textsl{(Self Dual Cone).} 
 $\Gamma^{*}= \mathbb{R}^2_{\geq 0}, \quad\text{or equivalently,}\quad \Gamma^* = \Gamma^0 = \mathbb{R}^2_{\geq0}.$

 \item[(2)] \textsl{(Gradient Map).}
 $$y = D_x \log \, \gg = \left(\frac{1}{x\_1} + \frac{1}{x_1 + x_2}, \frac{1}{x_2} + \frac{1}{x_1 + x_2}\right).$$

 \item[(3)] \textsl{(Inverse Gradient Map).}
 For $y_1,y_2>0$, set $\Delta \equiv y_1^2 - y_1y_2 + y_2^2$, which is $>0$. 
 Then
 $$x_1 = \frac{3}{2y_1 - y_2 + \sqrt{\Delta}},    \qquad   x_2 = \frac{3}{2y_2 - y_1 + \sqrt{\Delta}}$$
 $$
 x_1 + x_2 = \frac{3}{y_1 + y_2 - \sqrt{\Delta}}.
 $$
 
 \item[(4)] \textsl{(Dual Function).}
 $$\gg^*(y) = \frac{1}{27}\bigl(2y_1 - y_2 + \sqrt{\Delta}\bigr)                  \bigl(2y_2 - y_1 + \sqrt{\Delta}\bigr)                  \bigl(y_1 + y_2 - \sqrt{\Delta}\bigr)
 $$
 for $y_1,y_2>0$, and $\gg^*(y)$ is not a polynomial.
\end{enumerate}
}

\medskip
\noindent
\emph{Proof of (1).}
That $x_1x_2 + x_2x_2 \ge 0$ for all $x_1,x_2 \ge 0$ implies $y_1,y_2 \ge 0$ is obvious.

\medskip
\noindent
\emph{Proof of (2).}
This is a straightforward calculation.

\medskip
\noindent
\emph{Proof of (3).}
To solve
$$
y_1 = \frac{1}{x_1} + \frac{1}{x_1 + x_2}  \qquad\text{and}\qquad  y_2 = \frac{1}{x_2} + \frac{1}{x_1 + x_2},
$$\
consider the reciprocal variables
$$
a = \frac{1}{x_1},\quad b = \frac{1}{x_2},\quad c = \frac{1}{x_1 + x_2}.
$$
Then
$$
y_1 - c = \frac{1}{x_1} = a,  \qquad  y_2 - c = \frac{1}{x_2} = b.
$$
Hence,
$$
\frac{1}{c} = x_1 + x_2 = \frac{1}{y_1 - c} + \frac{1}{y_2 - c}              = \frac{y_1 + y_2 - 2c}{(y_1 - c)(y_2 - c)},
$$
or
$$
y_1y_2 - (y_1 + y_2)c + c^2 + c\bigl(2c - (y_1 + y_2)\bigr) = 0.
$$
Thus $c$ satisfies the quadratic equation
$$
 3c^2 - 2(y_1 + y_2)c + y_1y_2 = 0 
 \nm{16.2}
$$
with coefficients $A=3$, $B = - (y_1 + y_2)$, $C = y_1y_2$, and discriminant
$$
\Delta \equiv B^2 - AC = y_1^2 + y_2^2 + 2y_1y_2 - 3y_1y_2 = y_1^2 + y_2^2 - y_1y_2.
$$
Therefore
$$
c = \frac{-B \pm \sqrt{B^2 - AC}}{A}    = \frac{y_1 + y_2 \pm \sqrt{\Delta}}{3}.
$$
Selecting $- \sqrt{\Delta}$ yields
$$
\begin{aligned}
&c = \frac{1}{3}(y_1 + y_2 - \sqrt{\Delta}) \qquad {\rm and} \\
a= y_1-c = {1\over 3}(2y_1-&y_2+\sqrt{\Delta}), \quad b= y_2-c = {1\over 3}(2y_2-y_1+\sqrt{\Delta}) 
\end{aligned}
\nm{16.3}
$$

\bigskip
Taking the reciprocals of $a,b,c$ produces the formulas in (3). 

It is obvious that if $x_1,x_2 > 0$, then $y_1,y_2 > 0$. 
The gradient diffeomorphism theorem says that the map $y(x) = D_x\log g$ from $\mathbb{R}^2_{>0}$ to $\mathbb{R}^2_{>0}$ is a diffeomorphism. 
To see that the correct sign in $\pm\sqrt{\Delta}$ has been selected, we complete the proof of (3) by showing the following.

\Lemma{16.3} {\sl
If $y_1,y_2>0$, then
$$
x_1 = \frac{3}{2y_1 - y_2 + \sqrt{\Delta}}  \qquad\text{and}\qquad  x_2 = \frac{3}{2y_2 - y_1 + \sqrt{\Delta}}
$$
are both $>0$.
}

\noindent
\emph{Proof.}
It suffices to show that:
$$
 y_1,y_2>0 \Rightarrow 2y_1 - y_2 + \sqrt{\Delta} > 0 
 \quad\text{and}\quad 
 2y_2 - y_1 + \sqrt{\Delta} > 0. 
 \nm{16.4}
$$
By symmetry it is enough to show that:
$$
y_1,y_2>0 \  \Rightarrow\  2y_1 - y_2 + \sqrt{\Delta} > 0,  \text{ i.e. } \sqrt{\Delta} > y_2 - 2y_1.
 \nm{16.5}
$$
If $y_2 - 2y_1 \le 0$, then $\sqrt{\Delta} > 0$ is obvious, so we can assume that $y_2 - 2y_1 > 0$. 
Now
$$
\sqrt{\Delta} > y_2 - 2y_1  \iff \Delta > (y_2 - 2y_1)^2
$$
$$
\iff \ \ \Delta>(y_2-2y_1)^2
$$
$$
\iff \ y_1^2 +y_2^2-y_1y_2 \ >\ (y_2-2y_1)^2 \ =\ y_2^2+4y_1^2-4y_1y_2
$$
$$
\iff \ 3y_1y_2-3y_1^2 = 3y_1(y_2-y_1)>0
$$
$$
\iff\ y_2>y_1
$$
which is true since we were assuming that $y_2 > 2y_1$.
(Note that $\Delta - y_1y_2 = (y_1-y_2)^2 \ge 0$.)  \qed

\noindent
{\sl Proof of (4).}
Using the formulas in part (3) for the inverse gradient map and the defining identity for the dual function
$\gg^*(y)$, we have that
$$
\gg^*(y) = {1\over \gg(x)} = {1\over x_1x_2(x_1+x_2)} =  \frac{1}{27}\bigl(2y_1 - y_2 + \sqrt{\Delta}\bigr)                  \bigl(2y_2 - y_1 + \sqrt{\Delta}\bigr)                  \bigl(y_1 + y_2 - \sqrt{\Delta}\bigr).
$$
\qed

Finally, to see that $\gg^*(y)$ is not a polynomial, first note that 
$\gg^*(y) =  p(y) + q(y) \sqrt{\Delta}$ with $p$ a cubic polynomial and $q$ a quadratic polynomial.
If $\gg^*(y) = P(y)/Q(y)$ with $P,Q \in \bbr[y]$, i.e., if  $\gg^*(y)$ is a rational function, then 
$\sqrt{\Delta} = (P-pQ)/qQ$.

But $\Delta$ is a prime in $\bbr[y]$ and therefore in $\bbr(y)$.  If $\Delta$ were not prime, it would have a linear factor
and hence would vanish on a line in $\bbr^2$.  However, 
$$
\Delta \overset {\rm def} = y_1^2+y_2^2-y_1y_2 = \left(y_1-{1\over 2}y_2 \right)^2 +{3\over 4} y_2^2  >0
$$
unless $y_1=y_2=0$, and $\Delta$ is of degree 2.

This completes the proof of Proposition 16.2.\qed

\medskip
\Remark{16.4}
This Proposition 16.2 is a special case of the Pull Back Theorem 13.2 which illustrates the factoring of the  gradient diffeomorphism for the pull back, into two diffeomorphisms through a submanifold 
$\Sigma$.   Start with the \Ga pair: $h(v) = v_1v_2v_3$ and $ \G(h) = \bbr^2_{>0}$ on $\bbr^3$,
and pull back to $\bbr^2$ via the linear map $L:\bbr^2\to \bbr^3$ defined by 
$L(x)= (x_1, x_2, x_1+x_2) = v$.  Then the pull back $h_L = x_1x_2(x_1+x_2)\equiv g(x)$.

The reciprocal variables $(a,b,c)\in \bbr^3$, defined in the proof of (3), have
$$
\Sigma \equiv \{H(L(x) : x\in \bbr^2_{>0}\} = \{ 1/c = 1/a +1/b\},
$$
which equals $\{(a,b,c) : ab= c(a+b)\}$.  The map $L^t:\Sigma \arr \G^*_L(h)$ is given by
$y=L^t(a,b,c)=(a+c, b+c)$  with the inverse of this  diffeomorphism $L^t$ given by the formulas (16.3) for $a,b$ and $c$ in terms of $y_1,y_2 \in \G^*_L(h)$ at the end of the proof of part (c) of Proposition 16.2..

%%%%%%%%%%%%%%%%%%%%%%%%%%%%%%%%%%%%%%%%%%%%%%%%%%%%%%%%%
%%%%%%%%%%%%%%%%%%%%%%%%%%%%%%%%%%%%%%%%%%%%%%%%%%%%%%%%%
%%%%%%%%%%%%%%%%%%%%%%%%%%%%%%%%%%%%%%%%%%%%%%%%%%%%%%%%%
%%%%%%%%%%%%%%%%%%%%%%%%%%%%%%%%%%%%%%%%%%%%%%%%%%%%%%%%%
%%%%%%%%%%%%%%%%%%%%%%%%%%%%%%%%%%%%%%%%%%%%%%%%%%%%%%%%%
%%%%%%%%%%%%%%%%%%%%%%%%%%%%%%%%%%%%%%%%%%%%%%%%%%%%%%%%%
%%%%%%%%%%%%%%%%%%%%%%%%%%%%%%%%%%%%%%%%%%%%%%%%%%%%%%%%%
%%%%%%%%%%%%%%%%%%%%%%%%%%%%%%%%%%%%%%%%%%%%%%%%%%%%%%%%%
%%%%%%%%%%%%%%%%%%%%%%%%%%%%%%%%%%%%%%%%%%%%%%%%%%%%%%%%%
%%%%%%%%%%%%%%%%%%%%%%%%%%%%%%%%%%%%%%%%%%%%%%%%%%%%%%%%%

\def\LO{12}

\def\AAA{1}
\def\BB{2}
\def\CC{3}
\def\DD{4}
\def\EE{5}
\def\FF{6}
\def\HH{7}
\def\II{8}
\def\JJ{9}
\def\DU{10}
\def\GM{11}
\def\LO{12}
\def\BL{13}
\def\EX{14}
\def\Ph{15}

\def\Ph{17}

\def\ON{^{1\over N}}
\def\On{^{1\over n}}
\def\GD{G\aa rding-Dirichlet\ }

\vskip.5in

\centerline{\bf \headfont  \Ph. Completing the  Guo-Phong-Tong Lemma.}

%{\bf Concerning Lemma 4 in the paper of Guo-Phong-Tong.}

For operators that are (in our terminology) invariant \Ga-Dirichlet operators,  Guo, Phong and Tong [GPT] proved that majorization of the determinant implies a lower bound for the determinant  of the gradient
of the $N^{\rm th}$-root of the operator.  Here we extend this important lemma by not requiring invariance and by proving   the reverse implication.  Finally, the determinant majorization result in 
[HL$_9$] is expanded to a two-sided bound in Theorem \Ph.5.

\Prop{\Ph.1} {\sl  Let $\gg$ be a \Ga-Dirichlet operator of degree $N$ on $\rn$ with \Ga cone $\G$.  For each constant $\g>0$ the following
are equivalent:

(1) \ \ (Majorization of the Determinant)
$$
n\g^{1\over n} (\det A)^{1 \over n} \ \leq \ \gg(A)^{1\over N}\qquad\forall\, A>0
$$

(2)\ \ (Lower Bound for the Determinant of the Gradient)
$$
     D_B \gg^{1\over N}>0 
             \ \ {\rm and}\ \ 
\g \ \leq \ \det \, D_B \gg^{1\over N} \qquad \forall\, B \in \G.
$$
}

\noindent
To prove that (1) $\Rightarrow$ (2) we need the following upper bound for $\gg(A)\ON$.

\Lemma {\Ph.2} {\sl
For all $A, B \in \G$
$$
\gg(A)\ON \ \leq \bra{A}{D_B\gg\ON}.
\eqno{(\Ph.1)}
$$
Equality holds if and only if  $A$ and $B$ are positively proportional 
module the edge $E$  of $\G$.
}

\pf
Two hypotheses on $f\equiv \gg\ON$ are used as follows.

First, homogeneity-one says that equality holds when $A=B$, that is, 
$$
\text{(Euler) \hskip .5in  $f(B) \ =\ \bra{D_B}B  \qquad \forall \, B\in\G$}.
\eqno{(\Ph.2)}
$$
(In particular, Euler and positive 1-homogeneity in $A$ implies equality in (\Ph.1) if 
$A=\l B$ for $\l\geq 0$.)

Second,  concavity can be stated as the slope of the secant line above the segment 
$[B,A]$ is less than or equal to the slope of the tangent line to the graph of $f$ at $B$, that is,
$$
{\rm (Concavity)}  \hskip .3in f(A) - f(B) \ \leq \   \bra{D_B f}{A-B} \qquad \forall\, A,B \in\G.
 \eqno{(\Ph.3)}
$$
Adding these two facts  yields (\Ph.1)
$$
f(A) \ \leq\ \bra{D_Bf}{A}
 \eqno{(\Ph.1')}
$$

\noindent 
{\bf Note.} This proves (\Ph.1)$'$ with $\gg\ON$ replaced by any $f$ which is concave and homogeneous of degree one on a convex cone.

Finally, we  suppose that equality holds
in  (\Ph.1) for $A, B \in \G$ (the {\sl open} \Ga cone).  Since $\gg\ON$ is invariant under translations by the edge $E$,
 we may assume the edge  is 0.
 
 By (\EE.5) we see that at any $B\in\G$, the second derivative of $f$ is $\leq 0$, and it $=0$ on any tangent
vector $\z$ if and only if all the $B$-eigenvalues of $\zeta$ are equal to some $\l\in \bbr$.
Since $\l^j_B(B)=1$ for $j=1, ... , N$  (reviewed below\footnote{\ Note that  $\gg(tB + B)$  equals both  $(t+1)^N\gg(B)$ and $\gg(B) \prod_{i=1}^N(t+\l^i_B(B)).$ }),  this implies that all the $B$-eigenvalues of
$\z-\l B$ vanish, that is, $\z-\l B \in {N} \equiv$ the nullity of  $\gg, \G$.
By (5) in Section 2, we have that the nullity equals the edge which is zero, so $\z = \l B$.  
Thus for any $B\in\G$,
$$
D^2_B(\z)<0 \ \ \text{for any vector $\z$ not tangent to the ray through $B$}.
$$
Therefore the restriction of $f$ to any non-radial  closed interval $I\ss\G$ is strictly concave,
and so the inequality is strict on $I$. Of course, if $I$ is radial, the endpoints are 
positive multiples of each other since they lie in a open ray. \qed

\medskip

\noindent
{\bf Proof of (1) $\Rightarrow$ (2):} Combining (1) with Lemma \Ph.2 we have that
$$
\g\On(\det\, A)\On \ \leq \ {1\over n}\bra {A}{D_B \gg\ON} \qquad\forall\, A>0\ \ {\rm and }\ \ B\in \G.
$$
Hence, 
$$
\g\On \ \leq \   \underset{\det\, A=1} { \underset{A>0}  \inf} {1\over n} \bra{A}{D_B\gg\ON} 
\qquad \forall\ B\in\G.
$$
By Corollary \EE.4, $\D_B \gg\ON\geq 0$. By the classical fact (3) in (\LO.1c), which says that 
$$
 \underset{\det\, A=1} { \underset{A>0}  \inf } {1\over n} \bra PA\ =\ (\det P)\On\qquad{\rm if} \ \ P\geq 0,
\eqno{(\Ph.4)}
$$
this infimum equals $\det(D_B\gg\ON))\On$.  Thus
$$
\g\On \ \leq\ \det(D_B\gg\ON)\On\quad{\rm or}\quad 0\ <\ \g\ \leq\ \det(D_B\gg\ON).
\eqno{(\Ph.5)}
$$
Since $D_B \gg\ON\geq 0$ and has positive determinant, $D_B\gg\ON$ is positive definite, completing
the proof of (2)\qed

\bigskip

\noindent
{\bf Proof (2) $\Rightarrow$ (1):}  By the classical fact (\Ph.4) applied to $P\equiv D_B\gg\ON >0$,  assumption (2) yields
$$
\g\On \ \leq\ (\det D_B\gg)\On \ = \ {1\over n}   \underset{\det\, A=1} { \underset{A>0}  \inf}  \bra{A}{D_B\gg\ON} \ =\ {1\over n} \underset {A>0} \inf{ \bra A {D_B\gg\ON}   \over   (\det A)\On   }.
$$
Euler's Formula with $B=A>0$ says that $\bra{D_A\gg\ON}{A} = \gg\ON(A)$. Thus we have 
$\g\On \leq {1\over n} {  \gg\ON(A)   \over  (\det A)\On   }$  which is (1).  \qed

There is a stronger (equivalent) statement of the Majorization of the Determinant Inequality (1) above.

\Lemma {\Ph.3} {\sl
Under the hypotheses of Proposition \Ph.1,

(1) \ is equivalent to 

(1)$'$   \ \ $n\g^{1\over n} (\det \, A)^{1 \over n} \ \leq \ \gg^{1\over N}(B+A) - \gg^{1\over N}(B) 
\quad \forall\, B\in \G, \  {\rm and}\  A>0.$
}

\noindent
{\bf Proof that (1) $\Rightarrow$ (1)$'$.}
By the concavity of $\gg^{1\over N}$, one has
$$
\half \gg^{1\over N}(A)  + \half \gg^{1\over N}(B) \ \leq \ 
\half \gg^{1\over N} \left({A+B \over 2}\right) \quad \forall\, A,B \in\G.
$$
Since $\gg^{1\over N}$ is homogeneous of degree one, this implies that
$$
 \gg^{1\over N}(A)  +  \gg^{1\over N}(B) \ \leq \ 
 \gg^{1\over N}(A+B) \quad \forall\, A,B \in\G.
$$
Hence, we have
$$
\gg^{1\over N}(A+B) - \gg^{1\over N}(B) \ \geq\ \gg^{1\over N}(A)
\ \geq \ 
n\g^{1\over n} (\det \, A)^{1 \over n} 
$$
by (1) if $A>0$.\qed

\vskip .3in

%\Remark{\Ph.4}

\centerline{\bf \Ph.4 Expanding the Majorization of the Determinant }
\centerline{\bf to a Two-sided Bound}

The inequality in Proposition \Ph.1 (1) was established in [HL$_9$]  (cf. [HL$_8$])
for a wide class of \GD operators $\gg$.  It turns out that this inequality can be expanded to an elegant 
two-sided bound (\Ph.7) given in terms of the determinant and trace of the second derivative.

  The conditions on $\gg$ in [HL$_9$] were
   {\bf The Central Ray Hypothesis}     that $D_I\gg = k I$ for some $k>0$, and a certain coefficient condition.
   Here the constant $\g$ satisfies $n\g\On = \gg(I)\ON$, and the constant $k$ can be computed to be 
   $k=\gg(I)/n$.

Combining (1) with the upper bound of Lemma \Ph.2 we have the two-sided bound
$$
\gg(I)\ON (\det\, A)\On \ \leq \gg(A)\ON\ \leq\ \bra {A}{D_B\gg\ON}.
\eqno{(\Ph.6)}
$$
where the left inequality is for $A>0$ and the right one is for $A,B \in\G$.
Taking $B=I$,  using the Central Ray Hypothesis and the fact that $k=\gg((I)/n$, this becomes
$$
\gg(I)\ON (\det\, A)\On \ \leq \gg(A)\ON \ \leq \gg(I)^{{1\over N}-1} k \D(A) \ = \ \gg(I)\ON {1\over n} \D(A).
$$
where $\D(A) \equiv \tr(A)$.  This gives us the following very nice result.
 
 \Prop{\Ph.5}
 {\sl
For \GD operator $\gg$ as above (as in [HL$_9$]), one has the two-sided bound
 $$
\qquad\  \qquad\ \qquad \ \ \ \ (\det \, A)\On \ \leq \ \left( {\gg(A)\over \gg(I)}    \right)\ON \ \leq\ {1\over n} \D(A) \qquad{\rm requiring}
 \eqno{(\Ph.7)}
 $$  
$$
\text{ $A>0$ for the left inequality \ \ \ and \ \ \  $A\in\G$ for the right inequality.}
$$ }

Note that each term in (\Ph.7) is equal to 1 at $A=I$.

\Remark{\Ph.6}  This result reflects the subequation inclusions 
$$
\cp\ \ss\ \overline\G\ \ss\ \D
 \eqno{(\Ph.8)}
 $$   
In fact, if for $c>0$ we normalize $\gg(I)=1$  and set %normalize $\gg(I)=1$
$$
\begin{aligned}
\cp(c) \ &\equiv \ \{A \in \cp : \det(A)\On  \geq c\}, \\
\overline{\G(c)} \ &\equiv\ \{A\in \overline\G : \gg(A)\ON\geq c \}, \\
\D(c) \ &\equiv \ \{A\in \D :  \tr(A)/n \geq c\},
\end{aligned}
$$
then
$$
\cp(c)\ \ss\ \overline{\G(c)}\ \ss\ \D(c)
 \eqno{(\Ph.9)}
 $$

%\end{document} 
%%%%%%%%%%%%%%%%%%%%%%%%%%%%%%%%%%%%%%%%%%%%%%%%%%%%
%%%%%%%%%%%%%%%%%%%%%%%%%%%%%%%%%%%%%%%%%%%%%%%%%%%%
%%%%%%%%%%%%%%%%%%%%%%%%%%%%%%%%%%%%%%%%%%%%%%%%%%%%
%%%%%%%%%%%%%%%%%%%%%%%%%%%%%%%%%%%%%%%%%%%%%%%%%%%%
%%%%%%%%%%%%%%%%%%%%%%%%%%%%%%%%%%%%%%%%%%%%%%%%%%%%
%%%%%%%%%%%%%%%%%%%%%%%%%%%%%%%%%%%%%%%%%%%%%%%%%%%%
%%%%%%%%%%%%%%%%%%%%%%%%%%%%%%%%%%%%%%%%%%%%%%%%%%%%
%%%%%%%%%%%%%%%%%%%%%%%%%%%%%%%%%%%%%%%%%%%%%%%%%%%%
%%%%%%%%%%%%%%%%%%%%%%%%%%%%%%%%%%%%%%%%%%%%%%%%%%%%
%%%%%%%%%%%%%%%%%%%%%%%%%%%%%%%%%%%%%%%%%%%%%%%%%%%%
%%%%%%%%%%%%%%%%%%%%%%%%%%%%%%%%%%%%%%%%%%%%%%%%%%%%

\def\Aa{A}

\def\Fcal{{\mathcal F}}
\def\cG{{\overline \G}}
\def\bbM{{\mathbb M}}
\def\Mbf{{\mathbf M}}
\def\Nul{{\mathbf N}}
\def\iv{^{-1}}
\def\ggb{\overline {\gg}}
\def\dlog{D_A \log \, \gg}

\vskip.5in

\centerline{\bf \headfont  Appendix \Aa. Some Standard Ellipticity Remarks.}

As background for Section \DD we list some of the standard ways of describing ellipticity for an operator $F\in C(\Fcal)$
with domain $\Fcal \ss \Sn$.
Then we note they hold for a \GD operator $F=\gg$ with admissibility domain $\Fcal \equiv \cG$.

If $F$ is $C^1$, then the weakest form of  {\bf classical ellipticity} means that the {\bf linearization} $L_B$ of $F$ 
at each point $B$ has a {\bf coefficient matrix} $\Mbf\geq 0$.  By definition
$$
L_B(A) \ \equiv\ {d\over dt}\biggr|_{t=0} F(B+tA), 
$$ 
and in turn this equals $\bra {D_BF}A$.  So the coefficient matrix for $L_B$ is $\Mbf\equiv D_BF$.
In order to employ viscosity theory and have consistency of viscosity solutions to $F(D^2u)=0$
with $C^2$ solutions, one must require the weakest form of ellipticity in the viscosity context:

(\Aa.1) (Positivity) \ \ $\Fcal + \cp\ss\Fcal$, and

(\Aa.2) (Monotonicity) \ \ $F(A+P) \geq F(A)\ \ \forall\, A\in \Fcal$ and $P\geq 0$.

An additional compatibility condition between the operator $F$ and its domain $\Fcal$
is required to analyze the equation $F(D^2 u)=0$.

(\Aa.3) \ \ (Compatibility of $F$ and $\Fcal$)\ \ $\Int \Fcal = \{F>0\}$ and 
$\partial \Fcal = \{A\in \Fcal : F(A)=0\}$.

\noindent
If these conditions are met, we say that $F, \Fcal$ is an {\bf elliptic operator/subequation pair}.
In the case where $F\equiv \gg$ and $\Fcal=\G$ are a \Ga\ polynomial and its \Ga\, cone on $\rn$,
the condition (\Aa.1) is required by our definition of a \GD operator.  The Monotonicity Condition (\Aa.2)
follows from the formula (\EE.3) \note{??} for the first derivative of $\gg$.  This same formula (\EE.3) 
was used in Section \HH\ to prove that $D_B\gg \in \G^*$ if $B\in \G$ (Prop.\ \HH.1).  The completeness criterion (6)
in Prop.\ \CC.2 says that $\G^*\ss \Int \cp$.  This proves the following.

\Theorem {\Aa.1} {\sl
If $\gg, \G$ is a complete \GD operator, then the coefficient matrix $\Mbf \equiv D_B\gg$ for the linearization $L_B$
of $\gg$ at each point $B\in\G$, is positive definite, i.e., $L_B$ is {\bf uniformly} elliptic.
}

Compatibility (\Aa.3) of $\gg, \G$ is the statement that $\gg>0$ on $\G$ and $\partial \G = \{A\in \overline \G: \gg(A)=0\}$, which is true.

 %%%%%%%%%%%%%%%%%%%%%%%%%%%%%%%%%%%%%%%%%%%%%%%%%%%%%%%%%%%%%%%%%%%%%%%%%%%%%%%%%%%%%%%%%%%%%%%%%%%%%%%%%%%%%%%%%%%%%%%%%%%%%%%%%%%%%%%%%%%%%%%%%%%%%%%%%%%%%%%%%%%%%%%%%%%%%%%%%%%%%%%%%%%%%%%%%%%%%%%%%%%%%%%%%%%%%%%%%%%%%%%%%%%%%%%%%%%%%%%%%%%%%%%%%%%%%%%%%%%%%%%%%%%%%%%%%%%%%%%%%%%%%%%%%%%%%%%%%%%%%%%%%%%%%%%%%%%%%%%%%%%%%%%%%%%%%%%%%%%%%%%%%%%%%%%%%%%%%

 \def\Aa{A}
 \def\Ab{B}

\vskip.5in

\centerline{\bf \headfont  \Ab: Natural   Linear Transforms which Turn $\gg$ }
  \centerline{\bf   \headfont   and  its Dual  $\gg^*$  into Elliptic Differential Operators}

 Suppose $\gg$ is a \Ga polynomial on $\Sn$  which is $I$-hyperbolic, but whose \Ga cone $\G$ does not contain $\cp$.
 Assume  that 
 $$
 \G \ss \{A : \tr(A)\geq 0\}.
 $$  Then there is a family of linear transformations $L_{t}: \Sn \to \Sn$ which,  as $t \downarrow 0$, widens the
 cone $\G$ until at a certain first  point $t_0$, one has $\cp\ss L_{t_0} \G$, and so at this point $\gg\circ L_{t_0}\iv$ becomes a \GD polynomial. 
  Furthermore, for $0< t <t_0$,  the \Ga-Dirichlet polynomial $\gg\circ L_{t}\iv$
becomes uniformly elliptic. 
 
 It is a major point here that  this discussion also holds for the dual operator $\gg^*$.  That is,
in a siimilar  way $\gg^*$ can also be turned into %an elliptic and then 
a uniformly elliptic differential operator.
This gives {\bf  a natural duality between these differential operators}, which is associated to the basic duality 
of \S \DU.
  
  To start we consider the family of linear transformations $L_t : \Sn\to \Sn$, $t\in\bbr$,  given by 
  $$
  L_t(A)  \ \equiv\ A +( t-1)  (\tr\, A) {1\over n} I, \quad {\rm for} \ t\in \bbr.
  \eqno{(\Ab.1)}
  $$
 Note that
 $$
 \tr \, L_t(A) \ =\  t \, \tr\, A \and  L_1 \ =\ {\rm Id}.
  \eqno{(\Ab.2)}
  $$
%so that 
%$
%L_t\{ \tr = \l\} \ =\  \{ \tr  = t \l\}.
%  $
  
 There is  simple composition law 
 $$
 L_s \circ L_t \ =\ L_{st},
   \eqno{(\Ab.3)}
   $$
    Since $L_1 = {\rm Id}$, this implies that if $t\neq 0$, $L_t$ is invertible with inverse
   $$
   L_t\iv \ =\ L_{1\over t}  \qquad \text{for\  $t\neq 0$.}
   \eqno{(\Ab.4)}
  $$
From (\Ab.1), $L_t$ fixes the hyperplane $\{\tr  = 0\}$  for all $t\in \bbr$.
  
  The law (\Ab.3) follows from a simple calculation. However, it can be better understood under the decomposition 

\centerline{$\phi:\Sn \overset {\cong} \arr  \ \{\tr = 0 \} \oplus (\bbr \cdot  {1\over  n} I)$}

given by 
   $$
   \phi(A) \ \equiv\ \left (A-(\tr\, A) {1\over n}I, \ ( \tr\, A )  {1\over n}I\right )  \ \equiv\ (\s, \tau).
      \eqno{(\Ab.5)}
  $$
  With respect to these coordinates we see that 
  $$
  \boxit
 {$ L_t(\s, \tau) \ =\ (\s, t \tau)$}
     \eqno{(\Ab.6)}
  $$
  which means that $L_t = {\rm Id}_{\{\tr=0\}} \oplus \L_t$ where $\L_t=$ scalar multiplication by $t$.
  From this, the picture is entirely clear.
  
  \Lemma {\Ab.1}
  $$
  \text{ \sl $L_0$ is  projection onto the hyperplane $\{\tr = 0\}$ with kernel $\bbr I$.}
  $$
   $$
   \begin{aligned}
  &\hskip1.in \text{ \sl If $t>0$, $L_t$ fixes $\{\tr=0\}$ and}\\
  &\text{  \sl  if $t \geq 1$ it expands the upper and lower half  spaces by  $t$, and }\\
  &\text{   \sl if $t \leq 1$ it contacts the upper and lower half  spaces by  $t$}\\
  &\hskip 1. in \text{  \sl If $t<0$, then $L_t (\s,\tau)= (\s, -|t|\tau)$, so} \\
  &\text{  \sl the action is the same after interchanging the upper and lower half-spaces.}  
   \end{aligned}
   $$

  \def\cf{{\mathcal F}}
 \def\cb{{\mathcal B}}
 
 \medskip
 \noindent
 {\bf Some Assumptions:} 
 
 We now  suppose that $\cf$ is a closed convex cone in $\{\tr \geq 0\}\ss\Sn$ with vertex at the origin, and which contains $I$.
 We recall that the {\bf edge} of $\cf$ is 
 $$
 E(\cf)  \ \equiv\ \cf \cap (-\cf) \ \ss\ \{\tr=0\}.
 $$
 
 We shall make the additional assumption that 
 $$
 E(\cf) \ =\ \cf\cap \{\tr =0\}.
     \eqno{(\Ab.7)}
  $$
However, this  condition is {\bf always satisfied} if $\cf= \G$, the \Ga cone  for a \Ga polynomial $\gg$ on $\Sn$,
because here the edge equals the nullity of $\gg$ (see (5) in Section 2).

 \Lemma{\Ab.2} {\sl  Under these assumptions the following conditions are equivalent.}

(a)\ \ 
$E(\cf)= 0$.
  
(b)\ \ 
$ \cb \ \equiv \ \{A\in \cf : \tr\, A=1\} \ \ \text{is a compact set (and  so a compact base for $\cf$)}$
  
  (c)\ \ 
$ I \in \cf^* \ \ \text{($\equiv$ the open polar of $\cf$)}$
  
\pf
(a) $\Rightarrow$ (b):  Assume (a).  
Then $\cf\cap \{\tr = 1\}$ is closed since $\cf$ is closed. It is bounded since, if not, there exists
$v_j\in \cf\cap \{\tr=1\}$ with $|v_j|\to\infty$.  We can assume that $v_j/|v_j| $ converges to a non-zero element $v\in \{\tr =0\}$.
Then the ray $\overline{0 v_j}$ converges to the ray $\overline{0 v} \ss \{\tr=0\}\ss \cf$ (since $\cf$ is closed),
which with assumption (\Ab.7) contradicts (a).

\noindent 
(b) $\Rightarrow$ (c):  The statement (b) says that the base $\cb$ for $\cf$ is a compact set.  We have assumed that
$I\in \cf$, so ${1\over n}I\in \cb$.  Thus $\bra I v >0$ for all $v\in \cf-\{0\}$ which means that (c) holds.

\noindent 
(c) $\Rightarrow$ (a):   Note that $I^\perp = \{\tr=0\}$, so if (c) holds, then a neighborhood of $I$ lies in the polar of $\cf$.
This means that for any direction $u$ in $\{\tr=0\}$, there exist $A\in \cf^*$ with $\bra A u < 0$.  Thus $u=0$ and (a) holds.  \qed

   It is easy to see the following.
   
   \Lemma{\Ab.3} 
   {\sl With $\cf$ as above we have:
 $$
 \begin{aligned}
& \text{As $t\uparrow  \infty$, $L_t(\cf)$ contracts monotonically  down to the ray $\{t I: t\geq 0\}$, }  \\
&\hskip.3in   \text{and   $L_{-t}(\cf)$ contracts monotonically down to the ray $\{t I: t\leq 0\}$.} \\
& \text{As $t\downarrow 0$ $L_t(\cf)$ expands monotonically to $\{\tr >0\}\cup \{0\}$,} \\
 &\hskip .3in  \text{and $L_{-t}(\cf)$ expands monotonically to $\{\tr <0\}\cup \{0\}$.}
 \end{aligned}
 $$
 }

\vskip .5in
%\vfill\eject

\centerline{\headfont References} 

\noindent
[AO] S. Abja and G. Olive, {\sl Local regularity for concave homogeneous complex
degenerate elliptic equations dominating  the
Monge-Amp\`ere equation},   Ann.\ Mat.\ Pura Appl.\ (4) {\bf201} (2022) no.\ 2, 561-587.

\noindent
[ADO]  S. Abja, S. Dinew and G. Olive, {\sl Uniform estimates for concave
homogeneous complex degenerate elliptic equations comparable to the Monge-Amp\`ere
equation}, Potential Anal.  {\bf 59} (2023), no. 4, 1507-1524. \newline DOI 10.1007/s11118-022-10009-w

\noindent
[AV]  S. Alesker and M. Verbitsky, {\sl Plurisubharmonic functions on hypercomplex  manifolds and HKT-geometry}, J. Geom. Anal. {\bf 16} (2006), no. 3, 375-399.
 ArXiv:0510140v3

\noindent
[BGLS]   H. H.  Bauschke, O.  G\"uler, A. S.  Lewis, H. S.  Sendov, {\sl Hyperbolic polynomials and convex
analysis},  Canad. J. Math. {\bf 53} (2001), no. 3, 470-488.

\noindent
[CC] L.  Caffarelli and X. Cabr\'e, {\sl Fully nonlinear elliptic equations}, American
Math. Soc. Colloquium Publications, vol. 43, American Math.
Soc., Providence, RI, 1995.

\noindent
[CHLP] M. Cirant, R. Harvey,  H. B. Lawson and K. Payne,  {\sl Comparison Principles for General Potential Theories and PDEs}, 
Annals of Mathematics Studies, No. 218, Princeton University Press, Princeton, NJ, 2023.

\noindent
[CNS] L. Caffarelli, L. Nirenberg, and J. Spruck, {\sl The Dirichlet problem for nonlinear
second-order elliptic equations. III. Functions of the eigenvalues of the Hessian},
Acta Math. {\bf 155} (1985), no. 3-4, 261-301.

\noindent
[D] S. Dinew, {\sl Interior estimates for p-plurisubharmonic functions},
Indiana Univ. Math. J. {\bf  72} (2023), no. 5, 2025-2057.
arXiv:2006.12979.

\noindent
[G\aa r] L. G\aa rding, {\sl An inequality for hyperbolic polynomials}, J. Math. Mech. 8 (1959),
957-965.

\noindent
[G\"ul]  O. G\"uler, {\sl Hyperbolic polynomials and interior point methods for convex programming}, 
Mathematics of Operations Research, {\bf 22} No. 2, 1997, 350-377.

\noindent
[Gur]  l. Gurvits, {\sl Van der Waerden/Schrijver-Valiant like conjectures
and stable (aka hyperbolic) homogeneous
polynomials :
one theorem for all},  Electron. J. Combin. {\bf 15} (2008), no. 1, Research Paper 66, 26 pp.

\noindent
[GPT] B. Guo, D.H. Phong, and F. Tong, {\sl On $L^\infty$-estimates for complex Monge-Amp`ere equations}, Ann. of Math. (2) {\bf 198} (2023), no. 1, 393-418.
arXiv:2106.02224.

\noindent
[GP$_1$] B.Guo and D. H.Phong, {\sl On  $L^\infty$-estimates for fully nonlinear partial differential equations}, Ann. of Math. (2) {\bf 200} (2024), no. 1, 365-398.
  arXiv:2204.12549v1.

\noindent
[GP$_2$]  \ \------------ , {\sl Uniform entropy and energy bounds for fully non-linear equations}, 
Comm. Anal. Geom. {\bf 32} (2024), no. 8, 2305-2325.  arXiv:2207.08983v1.

\noindent
 [HL$_1$] F. R. Harvey and H.B. Lawson, {\sl Dirichlet duality and the nonlinear
Dirichlet problem}, Comm. Pure Appl. Math. 62 (2009), no. 3, 396-443.

\noindent
[HL$_2$] \ \------------ , {\sl Dirichlet duality and the nonlinear Dirichlet problem on Riemannian manifolds},
J. Diff. Geom. 88 (2011), no. 3, 395-482.

\noindent
[HL$_3$] \ \------------ , {\sl  Hyperbolic polynomials and the Dirichlet problem},    ArXiv:0912.5220.

\noindent
[HL$_4$] \ \------------ ,   {\sl  G\aa rding's theory of hyperbolic polynomials},
   {\sl Comm. in Pure  App. Math.}  {\bf 66} no. 7 (2013), 1102-1128.
  
\noindent
[HL$_5$] \ \------------ ,  {\sl Lagrangian potential theory and a Lagrangian equation of Monge-Amp\`ere type}, ,pp. 217- 257 in Surveys in Differential Geometry, edited by H.-D. Cao, J. Li, R. Schoen and S.-T. Yau, vol. 22, International Press, Somerville, MA, 2018.

\noindent
[HL$_6$]    \ \------------ ,    {\sl The inhomogeneous Dirichlet problem for natural operators on manifolds}, 
Annales de l'Institut Fourier, {\bf 69} No. 7 (2019),  3017-3064.
ArXiv:1805.11121.

 \noindent
[HL$_7$]    \ \------------ ,    {\sl Pluriharmonics in general potential theories}, Contemporary Mathematics, {\bf  735} (2019), 145-168,

\noindent
[HL$_8$]    \ \------------ ,  {\sl  Determinant majorization and the work of Guo-Phong-Tong and Abja-Olive},
   Calc. of Var. and PDE's  {\bf 62} article number: 153 (2023). ArXiv:2207.01729

\noindent
[HL$_9$]    \ \------------ ,  {\sl  A definitive determinant  majorization  result for nonlinear operators}, June, 2024, Duke J, 2025.
ArXiv:2407.05408.

\noindent
[HL$_{10}$]    \ \------------ ,  {\sl Some results for G\aa rding-Dirichlet operators on compact manifolds}, to appear.

   \noindent
[K]   N. V. Krylov,    {\sl  On the general notion of fully nonlinear second-order elliptic equations},    Trans. Amer. Math. Soc. (3)
 {\bf  347}  (1979), 30-34.

\noindent
[Ro]  R. T. Rockafellar, Convex Analysis, Princeton University Press, Princeton, NJ, 1970.

\noindent
[W]  Y. Wang, {\sl On the $C^{2,\a}$-regularity of the complex Monge-Amp\`ere equation}, Math.
Res. Lett. 19 (2012), no. 4, 939-946.

\noindent
[Y] S.T. Yau, {\sl On the Ricci curvature of a compact K\"ahler manifold and the complex Monge-Amp\`ere
equation. I}, Comm. Pure Appl. Math. {\bf  31} (1978), 339-411.

\vskip .3in

\noindent
F. Reese Harvey\\Rice University\\Houston,Texas, U.S.A.\\harvey@rice.edu

\noindent
H. Blaine Lawson, Jr. \\ Stony Brook University\\Stony Brook, New York, U.S.A.\\blaine@math.sunysb.edu

\end{document}